\documentclass{amsart}
\usepackage[utf8]{inputenc}
\usepackage[english]{babel}
\usepackage{amsmath,amsthm, amsfonts, amssymb, setspace, enumerate, graphicx, mathtools, booktabs, float, caption, subcaption, multicol}
\usepackage[dvipsnames]{xcolor}

\theoremstyle{plain}
\newtheorem{thm}{Theorem}
\newtheorem{cor}{Corollary}

\newtheorem{lem}{Lemma}

\theoremstyle{remark}
\newtheorem{rem}{Remark}

\newtheorem{exm}{Example}
\newtheorem{ques}{Question}

\theoremstyle{definition}
\newtheorem{defn}{Definition}

\newcommand{\R}{\mathbb{R}}
\newcommand{\s}{\mathbb{S}}
\newcommand{\D}{\mathbb{D}}

\newcommand{\dd}{\partial}

\begin{document}
\title{On Round Surgery on framed links and their Diagrams for 3-manifolds}

\author{Prerak Deep}
\address{Department of Mathematics, Indian Institute of Science Education and Research Bhopal}
\curraddr{}
\email{prerakdg@gmail.com}
\thanks{}

\author{Dheeraj Kulkarni}
\address{Department of Mathematics,  Indian Institute of Science Education and Research Bhopal}
\curraddr{}
\email{dheeraj@iiserb.ac.in}
\thanks{}



\date{}

\dedicatory{}

\begin{abstract}
We introduce a novel method for performing round surgery on framed links in a 3-manifold. In $\mathbb{S}^3$, these operations give rise to round surgery diagrams for 3-manifolds, directly analogous to classical Dehn surgery diagrams. We establish an explicit correspondence between a specific class of these round surgery diagrams and integral Dehn surgery diagrams. As a consequence, we prove a round surgery analogue of the Lickorish--Wallace theorem, demonstrating that any closed, connected, oriented 3-manifold can be obtained via round surgery on a framed link in $\mathbb{S}^3$. This result recovers Asimov's foundational theorem for oriented 3-manifolds within a purely diagrammatic framework. 
   
Since distinct round surgery diagrams can represent homeomorphic 3-manifolds, we naturally investigate whether a version of Kirby calculus exists in this setting. To this end, we define four types of diagrammatic moves and prove that any two round surgery diagrams representing the same 3-manifold are related by a finite sequence of these moves, thereby establishing a ``round version" of Kirby calculus. 

Finally, we provide two applications of this framework. First, we outline a method to construct standard Kirby diagrams directly from round surgery diagrams. Second, we establish the existence of taut foliations—and consequently tight contact structures—on 3-manifolds obtained via index-1 round surgeries on two-component fibered links in $\mathbb{S}^3$.
    % The bridge correspondence becomes one-to-one and onto up to the equivalence defined by the moves.

\end{abstract}

\maketitle

\section{Introduction}

In 1975, Asimov \cite{Asimov} introduced the notion of round handles as a generalisation of standard handles. A round handle is defined as the product of a circle with an ordinary handle, and a manifold $M$ is said to admit a round handle decomposition if it can be realised as a union of round handles. Asimov demonstrated that any compact flow manifold $M$ of dimension at least four admits a round handle decomposition (see Main Theorem in \cite{Asimov}). Shortly thereafter, in 1976, Thurston \cite{Thurs} introduced round Morse functions while investigating codimension-1 foliations, proving that the existence of a round Morse function on a manifold $M$ is equivalent to the existence of a round handle decomposition of $M$.

Attaching a round handle to an $n$-manifold along its boundary produces a new $n$-manifold with boundary. Analogous to ordinary handle attachments, the modification of the boundary resulting from a round handle attachment is called a \emph{round surgery} on the boundary manifold. Miyoshi effectively utilized round surgeries in \cite{Miyoshi} to characterize foliated cobordisms.

In dimension three, Dehn surgery is a thoroughly explored and foundational topic. Dehn surgery on a 3-manifold involves removing a tubular neighbourhood of a link and gluing back solid tori along the toral boundary components. The resulting closed 3-manifold is uniquely determined by the isotopy classes of the homeomorphisms used for the glueing. A cornerstone of this theory is the Dehn--Lickorish--Wallace Theorem \cite{L, W}, which guarantees a surgery presentation for every closed, connected, oriented 3-manifold. Consequently, any such 3-manifold can be described via a generic link projection in $\mathbb{R}^2$ adorned with an integer framing coefficient assigned to each component. While a 3-manifold admits infinitely many surgery presentations, Kirby's theorem \cite{RK, Book-GompfStip} establishes that any two such presentations are related by a finite sequence of diagrammatic moves, known as Kirby moves. Topologically, integral surgeries on 3-manifolds correspond to handle attachments on 4-manifolds; Kirby proved his theorem by analyzing handle decompositions of 4-manifolds and the deformations of their associated Morse functions.

In analogy with the Dehn--Lickorish--Wallace Theorem, Asimov showed that any compact, connected 3-manifold can be obtained via round surgeries on $\mathbb{S}^3$ (see Theorem Y in \cite{Asimov}). However, unlike the classical theorem, Asimov did not provide a systematic presentation of a 3-manifold via a framed link in $\mathbb{S}^3$. In the contact topology setting, Adachi \cite{Jiro, Jiro19} introduced contact round surgeries of index 1 on Legendrian and transverse links in contact 3-manifolds. However, a purely topological treatment of round surgery remains largely unexamined. In this article, we investigate a more general version of index-1 and index-2 round surgeries on links in smooth oriented 3-manifolds, addressing the following fundamental questions:

\begin{enumerate}
    \item Can we define round surgery on framed links in a 3-manifold, thereby establishing a formal round surgery presentation?
    \item If so, what is the precise relationship between round surgery on framed links and Asimov's original round surgery?
\end{enumerate}

Answering these questions allows us to draw robust parallels between Dehn surgery and round surgery, motivating whether 3-manifolds admit round surgery presentations akin to classical Dehn surgery diagrams via link projections in $\mathbb{R}^2$.

In Section \ref{Sec: RSDs}, we introduce a novel method for performing round surgeries on framed links in $\mathbb{S}^3$. We prove that on $\mathbb{S}^3$, this definition is equivalent to Asimov's round surgery. Specifically, we establish the following result:

\begin{thm}\label{thm: ARSequalsRS}
On $\mathbb{S}^3$, Asimov's round surgery is equivalent to performing a round surgery on a framed link.
\end{thm}

This equivalence allows us to propose a formal notion of a \emph{round surgery presentation} (or \emph{round surgery diagram}) to describe any closed, connected, oriented 3-manifold (see Section \ref{Sec: RSDs}). As a consequence, we obtain a round surgery analogue of the Lickorish--Wallace theorem:

\begin{cor}[Lickorish--Wallace Theorem for round surgery]\label{cor: AsimovResult}
Any closed, connected, oriented 3-manifold can be represented by a round surgery diagram consisting of index-1 and index-2 components.
\end{cor}

On $\s^3$, we note that a round 2-surgery on a knot yields a disconnected manifold unless it is performed in tandem with a round 1-surgery link. We call such a configuration a \emph{joint pair}, and a diagram composed exclusively of these pairs is a \emph{round surgery diagram of joint pairs}.% (We note that in Proposition 4.2 of \cite{Jiro19}, the author discussed an ordered pair of contact round surgeries mimicking a $(\pm 1)$-Dehn surgery, which can be viewed as a specific contact geometric realization of the topological joint pairs defined in Subsection \ref{subsec: Joint pair}.)

Furthermore, we investigate the problem of relating different round surgery diagrams that describe homeomorphic 3-manifolds, in the spirit of Kirby's theorem. We define four types of diagrammatic equivalence moves and establish the following result:

\begin{thm}[Equivalence theorem]\label{thm: EqThm}
Suppose $R_1$ and $R_2$ are two round surgery diagrams of joint pairs that yield the same closed, connected, oriented 3-manifold. Then $R_2$ can be obtained from $R_1$ by a finite sequence of equivalence moves 1--4.
\end{thm}

To achieve this, we establish a formal correspondence between classical Dehn surgery diagrams and a specific class of round surgery diagrams:

\begin{thm}[The Bridge Theorem]\label{BridgeTheorem}
The following two statements establish a bridge between integral Dehn surgery diagrams and round surgery diagrams:
\begin{enumerate}
    \item[(a)] Let $L=\bigcup_{i=1}^{n} (L_{i1} \cup L_{i2})$ be a round surgery diagram of joint pairs in $\mathbb{S}^3$ for a closed, connected, oriented smooth 3-manifold $M_{L}$, where $n_{ij}\in \mathbb{Z}$ is the round 1-surgery coefficient on $L_{ij}$ and $m_i\in \mathbb{Z}$ is the round 2-surgery coefficient on $L_{i2}$ for $1\leq i\leq n$ and $1 \leq j \leq 2$. Then there exists a corresponding Dehn surgery diagram on $L$ where $n_{i1}-n_{i2}+m_i$ and $m_i$ are the surgery coefficients on $L_{i1}$ and $L_{i2}$, respectively.
    \item[(b)] Conversely, given an integral Dehn surgery diagram $L= L_1 \cup \cdots \cup L_n$ of a closed, connected, oriented smooth 3-manifold with surgery coefficients $n_i \in \mathbb{Z}$ on each component $L_i$, there exists an equivalent round surgery diagram of joint pairs.
\end{enumerate}
\end{thm}

Using Theorem \ref{BridgeTheorem}, we provide a set of diagrammatic moves that relate any two round surgery diagrams representing the same 3-manifold (see Section \ref{Sec : KCforRSDJP}), thereby extending Kirby calculus to the round surgery framework. Additionally, Theorem \ref{BridgeTheorem} allows us to view Asimov's result (Corollary \ref{cor: AsimovResult}) as a direct manifestation of the classical Dehn--Lickorish--Wallace Theorem within the round surgery setting.

The article is organized into six sections. In Section \ref{sec: Prelim}, we review the preliminaries of Dehn surgery, Kirby's theorem, Asimov's round surgery, and round handles. In Section \ref{Sec: RSDs}, we formally define round surgery on framed links and explore its foundational properties. Section \ref{sec: BridgeThmC4} details the correspondence between Dehn surgery and round surgery, culminating in the proof of the Bridge Theorem. In Section \ref{Sec : KCforRSDJP}, we introduce the equivalence moves on round surgery diagrams of joint pairs and prove the Equivalence Theorem. Finally, Section \ref{Sec: Application} explores two applications of round surgery diagrams.

\section*{Acknowledgement}
This research is part of the PhD thesis of the first Author.
The first author was supported by the Prime Minister Research Fellowship-0400216, Ministry of Education, Government of India, during his PhD.

\section{Preliminaries}\label{sec: Prelim}
We say that a manifold is closed if it is a compact manifold without a boundary. 
By Dehn--Lickorish--Wallace theorem \cite{L} and \cite{W}, every closed connected oriented 3-manifold can be obtained by Dehn $(\pm1)$-surgeries on some link in $\s^3$. 
Thus, any closed connected oriented 3-manifold can be presented in terms of a diagram in $\R^2$ of a link projection with each component having $\pm1$ coefficient on top of it. 
These diagrams on links are called {\it Dehn surgery diagrams} of the corresponding 3-manifold, and the coefficients on the link components are called {\it Dehn surgery coefficients}. 
%A link with a fixed framing is called a {\it framed link}.
%Throughout this article, we refer to surgery diagrams mentioned above as \emph{Dehn surgery diagrams}. 

For a given knot $K$ in an oriented 3-manifold $M$, we denote by $N(K)$ a tubular neighbourhood of $K$, which is homeomorphic to the solid torus $\s^1 \times \D^2$.
We can specify a \emph{framing} by choosing a {\it longitude} $l$ on $\partial N(K)$. 
%We call this longitude a {\it chosen longitude}.
We say a knot $K$ is a {\it framed knot} if a framing is specified, i.e., a longitude is chosen on the boundary. 
We say a link $L$ in $M$ is a {\it framed link} if each of its components is a framed knot.
We denote the meridian, an essential simple closed curve that bounds a disk in $N(K)$, by a letter $m$.
For a knot $K\subset \s^3$, there is a natural choice of framing by the Seifert surface of $K$. In this case, the chosen longitude is called {\it canonical longitude}. 
We call the pair of meridian $m$ and canonical longitude $l$, the \textit{canonical coordinate system} on the boundary torus $\dd N(K)$. 
Throughout this article, the Dehn surgery coefficient is given with respect to the canonical coordinate system. 
For details about the Dehn surgery, the reader may refer to \cite{Book-GompfStip}, \cite{Book-Rolfsen}, \cite{Book-S}.

%Consider a knot $K\subset \s^3$. Let $S_{K}$ denote a Seifert surface of $K$. We define $l:= S_{K} \cap \dd N(K)$ to be the {\it canonical longitude} on $\dd N(K)$. 
%The curve $l$ is isotopic to the knot $K$. 
%Denote  by $\hat{n}$, a normal vector on $\dd N(K)$ pointing out of  $\overline{\s^3 \setminus N(K)}$. 
%We may choose an orientation on curves $m$ and $l$ such that $\left<m, l, \hat{n}\right>$ determines the orientation of $\s^3$.
%

Recall that an integral Dehn surgery on a link in a 3-manifold $M$ can be realised as the effect on one of the boundaries of $M \times [0,1]$ when a 2-handle $\D^2 \times \D^2$ is attached to it.
An $n$-ball $\D^k \times \D^{n-k}$ is called a $k$-handle of dimension $n$. 
It is known that any compact manifold can be decomposed into handles. A decomposition of a manifold into handles of various indices is called a {\it handle-body decomposition}.
In particular, a compact connected 4-manifold can be described on paper by drawing the attaching spheres of handles of indices 1 and 2, and mentioning other handles. 
These diagrams are known as the {\it Kirby diagram}.
%\textcolor{red}{Kirby diagrams that consist of only attaching spheres of 2-handles are the same as the integral Dehn surgery diagrams.} 

\begin{figure}[ht]
    \centering
    \includegraphics[scale=0.4]{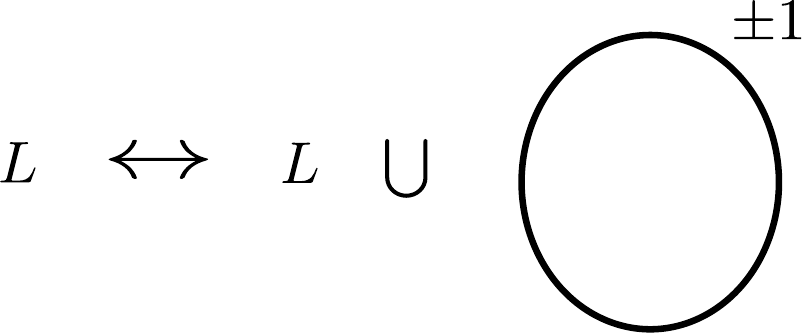}
    \caption{Kirby move of type 1}
    \label{fig:FirstKM}
\end{figure}

In \cite{RK}, Kirby discussed two types of moves on Dehn surgery diagrams having only integer coefficients. The moves introduced by Kirby do not change the resultant 3-manifold. 
The Kirby move of type 1 is an addition (or deletion) of an unknot with $\pm1$ framing to (or from) a surgery diagram (see Figure \ref{fig:FirstKM}). 

\begin{figure}[ht]
    \centering
    \includegraphics[scale=0.4]{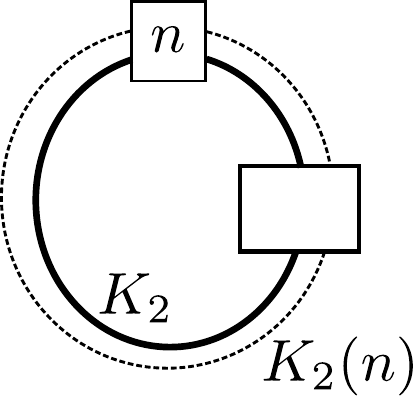}
    \caption{Framing curve $K_2(n)$ of $K_2$. The box on the right denotes the knotting of $K$, and the box in the middle denotes the linking of $K_2(n)$ with $K_2$.}
    \label{fig: FramedCurve}
\end{figure}

\begin{figure}[ht]
    \centering
    \includegraphics[scale=0.4]{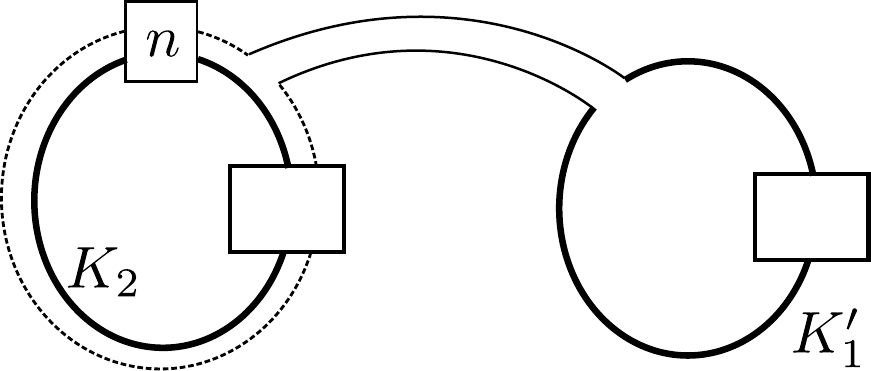}
    \caption{Band Connected sum of $K_1$ with $K_2$.}
    \label{fig: BandConnSum}
\end{figure}

Before we discuss the Kirby move of type 2. 
We recall the notion of band-connected sum. 
Suppose $K_1$ and $K_2$ are two framed knots with framing $m$ and $n$.
We denote $K_2(n)$ to be the framing curve of $K_2$ such that $lk(K_2, K_2(n))=n$ (see Figure \ref{fig: FramedCurve}).
Recall that the band connected sum of $K_1$ and $K_2(n)$, denoted by $K'_1:=K_1 \#_b K_2(n)$, is the attachment of a band as shown in Figure \ref{fig: BandConnSum}.

\begin{figure}[ht]
    \centering
    \includegraphics[scale=0.4]{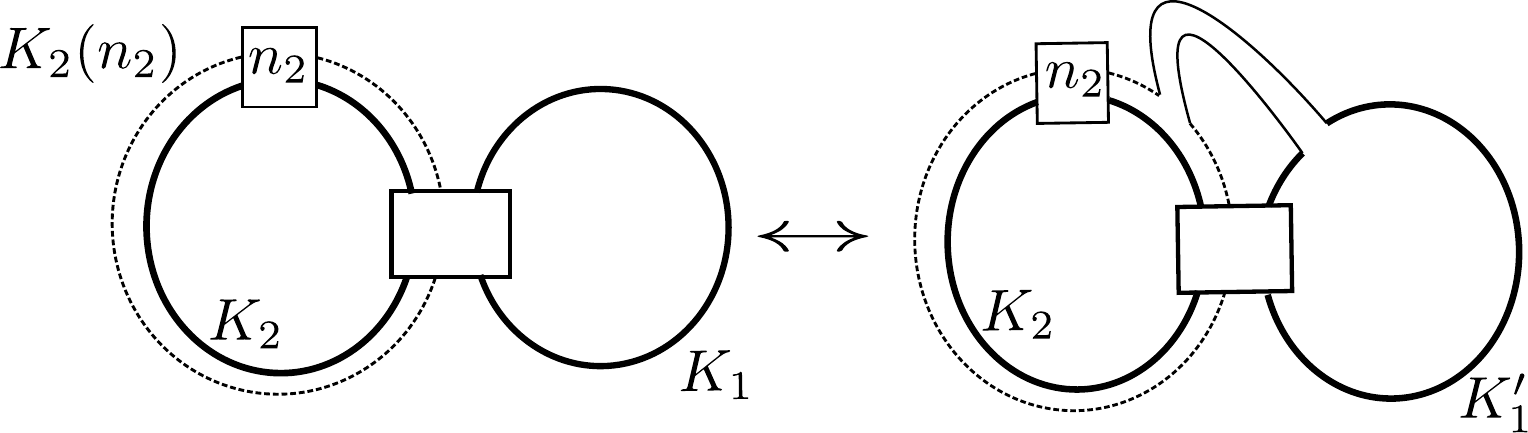}
    \caption{Kirby move of type 2. On both sides, the box in the middle denotes the knotting of each component and its linking.}
    \label{fig:SecondKM}
\end{figure}

The Kirby move of type 2 is the sliding of one component of the link $L$ over another. In particular, given two components $K_1$ and $K_2$ in $L$ with framing $n_1$ and $n_2$, respectively. 
We may slide $K_1$ over $K_2$ along the framing curve $K_2(n_2)$.  
We replace $K_1 \cup K_2$ by $K'_1 \cup K_2$, where $K'_1= K_1 \#_b K_2(n_2)$ and $K'_1$ gets a new framing given by $n_1 + n_2 + 2lk(K_1, K_2)$ (see Figure \ref{fig:SecondKM}). 
This move corresponds to the sliding of a 2-handle corresponding to $K_1$ over the 2-handle corresponding to $K_2$. 
It is easy to see that both moves do not change the resultant 3-manifold. Kirby proved that given a Dehn surgery diagram for a 3-manifold, the moves of type 1 and 2 are sufficient to describe all possible Dehn surgery diagrams for the same 3-manifold. 

%Every handle-body decomposition of $M$ is associated with a real-valued function on $M$ called a {\it Morse function}.  
In \cite{Asimov}, Asimov introduced the notion of a round handle as a generalisation of the ordinary handles. 
The round handles are defined as the product of a circle with an ordinary handle. More precisely, we have the following definition.
\begin{defn}
A product manifold $\s^1\times \D^k \times \D^{n-k-1}$ is called an {\it $n$-dimesional round handle of index $k$}. 
\end{defn}

The boundary of the round handle of index $k$ has two components given by $\s^1\times \s^{k-1} \times \D^{n-k-1}$ and $\s^1 \times \D^k \times \s^{n-k-2}$. 
We call $\s^1\times \s^{k-1} \times \D^{n-k-1}$ to be an \emph{attaching region} of the round handle of index $k$. 

\begin{defn}
     An $n$-manifold $X^{\prime}$ is obtained by {\it attaching an $n$-dimensional round handle of index $k$} to $X$ if $X^{\prime}= X \bigcup_{\phi}( \s^1 \times \D^k \times \D^{n-k-1})$, where the attaching region of the round handle is embedded in $\partial X$ by a specified embedding $\phi: \s^1 \times \s^{k-1}\times \D^{n-k-1} \to \partial X$. 
\end{defn}

We say that $n$-manifold $X$ admits a round handle decomposition if $$X=  R_0  \cup_{\phi_1} R_1 \cup_{\phi_2} \cdots \cup_{\phi_l} R_l,$$ where each $R_k$ denotes a round handle of some index 
$0\leq k\leq n-1$. %The manifold $X$ is obtained inductively by attaching $R_i$ to \textcolor{red}{$\partial_{+}$} of intermediate manifold obtained by preceeding round handles. 

In a handle-body decomposition of a manifold, we say two handles $h_k$ and $h_{k+1}$ are {\it independently attached} if their attaching regions are disjoint.
The following result from \cite{Asimov} proves that two independently attached handles of indices $k$ and $k+1$ can be isotoped to a round handle of index $k$. 
 
\begin{thm}[Fundamental Lemma of Round Handles, \cite{Asimov}]\label{Thm: FRH}
Let $(X, \dd_{-}X, \dd_{+}X)$ be a flow manifold with non-empty boundary $\partial X$.
Let $X^{\prime} = X \cup_{\psi_1} h_k \cup_{\psi_2} h_{k+1}$, where $h_k$ and $h_{k+1}$ are attached independently via attaching maps $\psi_1$ and $\psi_2$ to $\partial_{+} X$. Then if $k \geq 1, \, X^{\prime} = X \bigcup_{\phi} R_k$ via some attaching map $\phi$.
\end{thm}

Conversely, any attached round $k$-handle can be realised as an attached pair of a $(k-1)$-handle attachment and a $k$-handle, as discussed in the following lemma. 

\begin{lem}[\cite{MFS}]\label{Lem: RHasOH}
    Any round $k-$handle can be realised as a pair of ordinary handles of indices $k$ and $k+1$.  
\end{lem}
 
The effect of attaching a $(n+1)$-dimensional round $k$-handle to an $(n+1)$-manifold $X$ on its boundary is a {\it round $k$-surgery} on $\dd X$ for $0\leq k\leq n-1$. More precisely, we have the following definition of the round surgery.
\begin{defn}
Let $N$ be $n$-dimensional manifold. Let $\phi : \s^1 \times \s^{k-1} \times \D^{n-k} \to N $ be an embedding. 
A round $k$-surgery on $n$-manifold $N$ is the operation of removing the embedded region $\phi\left(\s^1 \times \s^{k-1} \times \D^{n-k}\right) $ from $N$ and glueing $\s^1 \times \D^{k} \times \s^{n-k-1}$ to get a new $n$-manifold $M$. More precisely, we have, 
$$M := \overline{N \setminus \phi(\s^1 \times \s^{k-1}\times \D^{n-k})} \bigcup_{\text{id}} \s^1 \times \D^{k} \times \s^{n-k-1}$$
for $0\leq k\leq n-1$. The manifold $M$ is said to be obtained from performing round $k$-surgery (or round surgery of index $k$) on $N$ along the embedding $\phi$. 
\end{defn}

\begin{exm}\label{ex: R1S-3-torus}
We consider a round 1-surgery on $\s^3$ along the embedding $\phi: \s^1\times \s^0 \times \D^2 \to \s^3$ defined as follows: $\phi (\s^1 \times \s^0 \times \D^2)$ is a tubular neighbourhood $N$ of the Hopf link, and $\phi(\s^1 \times \{p\} \times \{-\}) \mapsto l_{p}$ and $\phi(\{-\} \times \{p\} \times \s^1) \mapsto m_{p}$, where $(m_p, l_p)$ denotes the pair of simple closed curves defining the canonical coordinate system on boundary tori $\dd N$, for each $p \in \s^0$. 
Following the above mapping, we glue a thickened torus $\s^1 \times \D^1\times \s^1$. After glueing, we obtain a circle bundle over a torus. %$3$-torus. +
Thus, this round 1-surgery on $\s^3$ along embedded solid tori $\phi (\s^1 \times \s^0 \times \D^2)$ produces a circle bundle over a torus. %$\mathbb{T}^3$.  
\end{exm}

\begin{exm}\label{ex: R2S-s3U(s1xs2)}
We consider a round 2-surgery along the embedded thickened torus inside a tubular neighbourhood of the unknot $N(U)$ such that the $\phi(\s^1 \times \{-\} \times \{-\})$ maps to the meridian of the solid torus. 
Then the resultant $3$-manifold is $\s^3 \sqcup (\s^1 \times \s^2)$.
\end{exm}

\section{Round Surgery on framed links}\label{Sec: RSDs}

In this section, we define round surgery of indices 1 and 2 on framed links in an oriented 3-manifold $M$. Further, we prove that on $\mathbb{S}^3$, this notion of round surgery is equivalent to Asimov's round surgery. Thus, we can associate a framed link to a sequence of round surgeries on $\mathbb{S}^3$, and we call this associated framed link a \emph{round surgery diagram}. We can think of it as an analogue to a Dehn surgery diagram for an oriented smooth 3-manifold.

Below, we define round surgeries of indices 1 and 2 on framed links in $M$ with rational round surgery coefficients. 
We know that any simple closed curve $c$ on a torus can be expressed with respect to the meridian $m$ and the chosen longitude $l$ as $c= p\cdot m + q \cdot l$. Thus, we get a fraction $\dfrac{p}{q} \in \mathbb{Q}\cup \{\infty\}$. 
The round surgery coefficient on each component of the framed link is such a fraction $\dfrac{p}{q}$ describing a simple closed curve on the boundary torus. 
In the case of a knot $K\subset \s^3$, there is a canonical choice of longitude $l$ given by a Seifert surface of $K$ on $\dd N(K)$.   
%Thus, we can express any simple closed curve $c$ with respect to the meridian $m$ and the canonical longitude $l$ as follows: 
%$$ c= p\cdot m + q\cdot l.$$ 
%Moreover, the canonical longitude $l$ and the meridian $m$ can be used to identify the boundary torus with $\R^2/\mathbb{Z}^2$ by identifying  $m$ to $(1,0)$ and $l$ to $(0,1)$.  
Therefore, the meridian $m$ and canonical longitude $l$ determine the \textit{canonical coordinate system} on $\dd N(K)$.

\subsection{Round 1-surgery on a link with two components}
Suppose $L_1\cup L_2$ denotes a link with two components in $M$. Let integers $n_1$ and $n_2$ be round 1-surgery coefficients on $L_1$ and $L_2$, respectively.  
We remove the tubular neighbourhood $N(L_i)$ of each component from $M$. We obtain the link-complement $$N:=\overline{M\setminus(N(L_1)\cup N(L_2))}.$$ 
Now, we glue a thickened torus $\s^1 \times \D^1\times \s^1$ to $N$ by the glueing map $\phi: \dd(\s^1 \times \D^1 \times \s^1) \to \dd N$.  
The thickened torus has two boundary torus $\s^1 \times \{\pm1\}\times \s^1$, here we treat $\s^0 = \{\pm1\}$.
The boundary torus $\s^1 \times \{(-1)^i\} \times \s^1$ glues to $\partial N(L_i)$. 
The glueing map $\phi$ maps the curve $\{q\} \times \{(-1)^i\} \times \s^1$ to the meridian $m_i$ of $L_i$ and the curve $\s^1 \times \{(-1)^i\} \times \{p\}$ is mapped to the curve $n_i \cdot m_i + l_i$, where $l_i$ denotes the chosen longitude of the link component $L_i$ and for some $p,q \in \s^1$.

We denote by $M_{(L^{n_1}_1 \cup L^{n_2}_2)}$ the $3$-manifold obtained by performing a round $1$-surgery on $L_1 \cup L_2$ with round $1$-surgery coefficients $n_1$ and $n_2$ on $L_1$ and $L_2$, respectively.

%The following diagrams are some examples of round surgery diagrams on $\s^3$.
\begin{exm}\label{ex: R1SD-3torus}
Example \ref{ex: R1S-3-torus} of round 1-surgery on $\s^3$ can be presented by the Hopf link with round 1-surgery coefficient $0$ on each component as shown in the Figure~\ref{fig:3-torus}. %$\mathbb{T}^3$ can be presented by round $1-$surgery diagram as a Hopf link with $0$ framing on both components. 
    \begin{figure}[ht]
    \centering
    \includegraphics[scale=0.4]{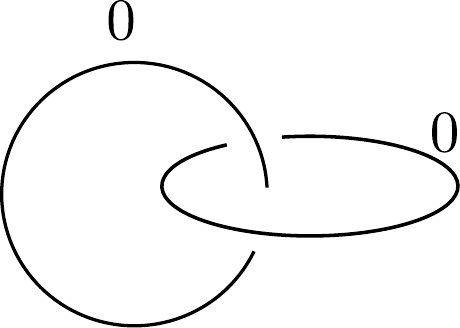}
    \caption{A circle bundle over a torus.}
    \label{fig:3-torus}
\end{figure}
\end{exm}

\begin{thm}\label{thm: SameClassOfR1SD}
Suppose some integers $n_1, n_2, n^{\prime}_1$ and $n^{\prime}_2$ satisfies $n_1-n_2 = n_1^{\prime}-n^{\prime}_2=k$. Then,
$$M_{(L_1^{n_1}\cup L_2^{n_2})}\text{ is diffeomorphic to } M_{(L_1^{n^{\prime}_1}\cup L_2^{n^{\prime}_2})}.$$ \end{thm}

\begin{proof}
   Let $m_i$ and $l_i$ denote the meridian and longitude on the boundary tori $\dd N(L_i)$ for each $i=1,2$.
    Suppose $x= \{q\} \times \{(-1)^i\} \times \s^1 $ and $y= \s^1 \times \{(-1)^i\} \times \{p\}$ denote the coordinate system on the boundary of the glued thickened torus $\mathbb{T}^2 \times [-1,1]$, where $i=1,2$.
    Then, we glue $\mathbb{T}^2 \times \{(-1)^i\}$ to $\dd N(L_i)$ as 
    \begin{eqnarray*}
        x &\mapsto &  m_i\\
        y &\mapsto & n_i\cdot m_i + l_i.
    \end{eqnarray*}
    %for each $j=1,2$. 
    Suppose $\phi_i$ denote the above homomorphism from $\mathbb{T}^2 \times \{(-1)^i\}$ to $\dd N(L_i)$. 
    Then we can realise glueing of the thickened torus as a homeomorphism $\phi_2\circ \phi_1^{-1}$ from $\dd N(L_1)$ to $\dd N(L_2)$. 
    It maps $\phi_2 \circ \phi_1^{-1}(m_1)=m_2$ and  
    \begin{eqnarray*}
    \phi_2\circ\phi_1^{-1}(l_1)= \phi_2(y-n_1 \cdot x)= (n_2-n_1)\cdot m_2 + l_2 = k \cdot m_2 + l_2.
    \end{eqnarray*}
    Thus, the homeomorphism corresponds to $\begin{pmatrix}
         1 & k \\
         0 & 1
    \end{pmatrix} \in SL(2, \mathbb{Z})$ with respect to the given basis on $\dd (N(L_1)\cup N(L_2))$.
    Therefore, the glueing homeomorphism is independent of the coefficients $n_1$ and $n_2$ but depends only on their difference $k$. 
    %Notice that, for fixed $k\in \mathbb{Z}$ and different values of $n\in \mathbb{Z}$, we get the same matrix for the homeomorphism. 
    Hence, the resultant 3-manifolds, after performing round 1-surgeries on $L$ with coefficient $n$ on $L_{1}$ and $n+k$ on $L_2$, for $n\in \mathbb{Z}$ and a fixed integer $k$, are diffeomorphic.  
\end{proof}

\begin{rem}\label{rem: SameR1SC}
    %For $k=0$, we get the same round 1-surgery coefficient on both of the components of $L$. Thus, by the above Lemma \ref{lem: SameResultAfterR1S}, round 1-surgeries on a link $L_{1}\cup L_{2}$ with the same coefficients on both of the components correspond to the diffeomorphic 3-manifold. 

As a consequence, round 1-surgeries on the same link with round 1-surgery coefficients on both components to be the same integer produce the diffeomorphic 3-manifolds.
\end{rem}

%Recall that round 1-surgery on $\mathbb{S}^3$ is an operation of removing embedded $\s^1 \times \s^0 \times \D^2$ from $\mathbb{S}^3 $ and glueing $\s^1 \times \D^1 \times \s^1$ by the identity map. 
%Thus, $$N:= \overline{\mathbb{S}^3 \setminus \phi(\s^1\times \s^0 \times \D^2)} \bigcup_{\text{id}} \s^1 \times \D^1 \times \s^1; $$
%where $\phi :\partial (\s^1\times \s^0 \times \D^2) \to \mathbb{S}^3$ is an embedding. 
 
\begin{lem} \label{Lem:R1Diagram}
Asimov's round 1-surgery on $\s^3$ can be completely determined by a two-component framed link in $\s^3$. 
\end{lem}

\begin{proof}
In round 1-surgery, we remove two solid tori and glue in a thickened torus. 
Recall that any embedded solid torus $\s^1 \times \D^2$ in $\mathbb{S}^3$ can be realised as a tubular neighbourhood of its core curve $\s^1 \times \{0\}$. 
Thus, a link $L= L_1 \cup L_2$ with two components determines the embedding of $\s^1 \times \s^0 \times \D^2$ in $\s^3$ for a round 1-surgery. 
In particular, $\phi(\s^1 \times \s^0 \times \D^2)= N(L_1) \cup N(L_2)$ where $N({L_j})$ denote the tubular neighbourhood of the link component $L_j$. 
Next, we prove that the glueing homeomorphism of a thickened torus is completely determined by an integer surgery coefficient on each component of the link. 

On $\partial N(L_j)$, we consider a canonical coordinate system given by meridian $m_j$ and longitude $l_j$. 

We want to express the round 1-surgery coordinates curves with respect to $m_{j}$ and $l_{j}$.
Denote the round 1-surgery coordinates curves by $r_{j}=\phi(\s^1 \times \{(-1)^j\} \times \{-\})$ and $ s_{j}= \phi(\{-\} \times \{(-1)^j\} \times \s^1)$. 
We have
\begin{eqnarray*}
    r_{j} =&  p_{j}\cdot m_{j} + q_{j}\cdot  l_{j}\\
   \text{and } \quad s_{j} =& u_{j}\cdot  m_{j} + v_{j} \cdot l_{j},
\end{eqnarray*}
for some $p_{j}, q_{j}, u_{j} \text{ and } v_{j}$ in $\mathbb{Z}$. 
Moreover, for $j=1,2$, the curves $\s^1 \times \{(-1)^j\} \times \{p\}$ and $\{q\} \times \{(-1)^j\} \times \s^1$ form a basis of the first homology of the boundary torus $\s^1 \times \{(-1)^j\}\times \s^1$ of the thickened torus $\s^1 \times [-1,1]\times \s^1$. 

Since we glue $\s^1\times \D^1 \times \s^1$ by the identity homeomorphism, the images of the curves $\s^1 \times \{(-1)^j\}\times \{p\} \mapsto  r_{j} = \phi(\s^1 \times \{(-1)^j\} \times \{p\}) \text{ and } \{q\} \times \{(-1)^j\}\times \s^1 \mapsto s_{j} =\phi(\{q\} \times \{(-1)^j\} \times \partial \D^2)$. 

Since $\phi(\{q\} \times \{(-1)^j\} \times \partial \D^2)$ is a disk bounding curve in $N(L_j)$, it is a meridian on the boundary torus $\partial N(L_j)$.
The meridian curve is unique up to isotopy; we get $s_{j} = m_{j}$.
Therefore, $u_{j}= 1$ and $ v_{j}= 0$ for each $j=1,2$. 

Further, the curve $r_{j}$ is isotopic to the core curve $L_j (=\s^1 \times \{(-1)^j\} \times \{0\})$ of $N(L_j)$. 
We get $ r_{j} =  n_{j} \cdot m_{j} + 1 \cdot l_{j}$, i.e. $p_{j}=n_{j}$ and $q_{j}=1$ for some $n_j \in \mathbb{Z}$. 
Therefore, the round 1-surgery is completely determined by $(\frac{n_j}{1}, \frac{1}{0})$.
Since the image $s_j$ is always fixed, embedding can be completely determined by a pair of integers $n_1$ on $L_1$ and $n_2$ on $L_2$. 
Hence, a round surgery can be presented by an integral framed link with two components. 
\end{proof} 

%\begin{rem}The integers $n_1$ on $L_1$ and $n_2$ on $L_2$ are called the {\it round 1-surgery coefficients} on the respective components.\end{rem}

\begin{rem}
Suppose $M$ is obtained by round 1-surgery on a framed link $L$ in $\s^3$. 
Then $\pi_1(M)$ is not trivial. 
This is because the interval $I$ of the glued thickened torus can be joined with another interval $I_0\subset \s^3\setminus N(L)$ to form a circle that is non-contractible in the new 3-manifold, where $N(L)$ is a tubular neighbourhood of the round 1-surgery link $L$.
\end{rem}

%\begin{lem}\label{lem: SameResultAfterR1S}
 %   Let $L= L_{1}\cup L_{2}  \subset \s^3$ denote a round 1-surgery link  with round 1-surgery coefficients $n_1=n$ on $L_1$ and $n_2=n+k$ on $L_2$. Then, we obtain diffeomorphic 3-manifolds after performing round 1-surgeries on $L$ for any $n\in \mathbb{Z}$ and a fixed $k\in \mathbb{Z}$.
%\end{lem}

\subsection{Round 2-surgery on a framed knot}
Suppose $K$ is a knot in $M$ and $\dfrac{p}{q}\in \mathbb{Q}\cup \{\infty\}$ is the round 2-surgery coefficient. 
Let $N_d(K)$ denote a tubular neighbourhood of $K$ of some radius $d$. 
For some $0<d_1 < d_2 < d$, we remove the thickened torus $N_{d_2} \setminus N_{d_1}$ from $M$ and obtain
$$N : = \overline{M \setminus \left(N_{d_2} \setminus N_{d_1}\right)}.$$
Now, we glue two copies of a solid torus, i.e., $\s^1 \times \D^2 \times \s^0$. 
For some $q \in \s^1$, we map the meridians $\{q\}\times \partial \D^2 \times \{(-1)^i\}$ to $p \cdot m_i + q \cdot l_i$, where $m_i$ and $l_i$ denote the meridian and chosen longitude on $\partial N_{d_i}$ for $i=1,2$.

\begin{exm}\label{ex: R2SD-s3U(s1xs2)}
The Example~\ref{ex: R2S-s3U(s1xs2)} of round 2-surgery on $\s^3$ producing $\s^3 \sqcup (\s^1 \times \s^2)$ can be presented by a $0$-framed unknot as shown in Figure \ref{fig:0-unknot}. 
\begin{figure}[ht]
    \centering
    \includegraphics[scale=0.4]{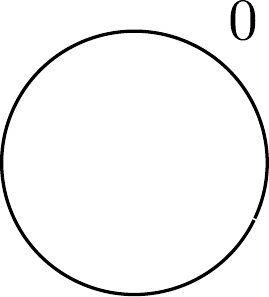}
    \caption{$\s^3 \sqcup (\s^1\times \s^2) $}
    \label{fig:0-unknot}
\end{figure}
\end{exm}

\begin{exm}
    $\s^3 \sqcup \s^3$ can be presented by a round 2-surgery diagram on a $(\pm1)$-framed unknot.
\end{exm}

We are interested in associating knot diagrams to resultant 3-manifolds. For that, we study round surgery of index 2 on $\s^3$. 
The following lemma is required to define the round surgery diagram for the round surgery of index 2. 

\begin{lem}\label{Lem:EmThickTorus}
Any embedded thickened torus $\mathbb{T}^2\times [1,2]$ in $\s^3$ is homeomorphic to $N_2(K)\setminus N_1(K)$ for some knot $K \in \s^3$, where $N_r$ is a tubular neighbourhood homeomorphic to a solid torus $\s^1 \times \D^2_r :=\{(\theta, (x,y)) \in \s^1 \times \R^2 | x^2 +y^2 = r^2 \text{ and } \theta \in \s^1\}$.

\end{lem}
\begin{proof}
Suppose $\phi: \mathbb{T}^2 \times [1,2] \to \s^3$ is an embedding. 
Then by the solid torus theorem \cite{STT}, there exists a solid torus $N_1$ bounded by either the torus $\mathbb{T}^2 \times \{1\}$ or $\mathbb{T}^2 \times \{2\}$. 
Without loss of generality, we assume that $\dd N_1 = \mathbb{T}^2 \times \{1\}$. 
Further, we can glue $N_1$ with the $\mathbb{T}^2 \times [1,2]$ to produce a new solid torus $N^{\prime} := N_1 \bigcup_{\mathbb{T}^2\times \{1\}} \mathbb{T}^2 \times [1,2]$ with boundary $\mathbb{T}^2 \times \{2\}$.
Suppose $c_1$ denotes the core curve of the solid torus. 
Thus, $N'$ and $N_1$ are tubular neighbourhoods of the knot $c_1$. 
Hence, $\phi(\mathbb{T}^2 \times [1,2])$ is homeomorphic to $N'\setminus N_1 \subset N'(c_1)$ for some knot $c_1$.  

We would like to claim that the choice of knot $c_1$ is unique up to isotopy.
Suppose there is another embedded solid torus $N_2$ bounding $\mathbb{T}^2 \times \{2\}$ given by an embedding $\psi: \s^1 \times \D^2 \to N_2 \subset M$. 
We denote $c_2$ to be the core curve of the solid torus $N_2$. 
Since $N_2 \neq N'$, we claim that $N_2 \cap N' = \mathbb{T}^2 \times \{2\}$.
Suppose $p \in N_2 \cap N'$ such that $p\notin \mathbb{T}^2 \times \{2\}$. 
The point $p$ is in the interiors of $N'$ and $N_2$. 
Let $s_1$ and $s_2$ be straight line segments lying in the interiors of $N'$ and $N_2$ joining $p$ with some point $q_1 \in c_1$ and $q_2 \in c_2$, respectively. 
Thus, the path $s_1 \cup s_2$ is a path connecting the interior of $N'$ with $N_2$ without intersecting $\mathbb{T}^2\times \{2\}$. 
We remove the embedded torus $\mathbb{T}^2 \times \{2\}$ from $\s^3$ and get a 3-manifold with two connected components. 
The union of the interiors is contained in one of the components. 
By the solid torus theorem, the interior of $N'$ is one of the components of the complement. 
Thus, the interior of $N_2$ is a subset of the interior of $N'$. 
Since $\dd N' =\dd N_2$, we get $N'= N_2$, which is a contradiction to $N' \neq N_2$. 
Therefore,  $N_2 \cap N' = \mathbb{T}^2 \times \{2\}$. 
In particular, the solid torus $N_2 $ is the knot complement of $c_1$ in $\s^3$. 
By Heegard splitting of genus 1 of $\s^3$, we know that the knot complement of a knot $K$ in $\s^3$ is a solid torus; then the knot is isotopic to the unknot.
%knots in $\s^3$ are isotopic if their knot complements are homeomorphic by an orientation-preserving homeomorphism. 
Since the knot complement is a solid torus, it implies that the core curves $c_1$ and $c_2$ are unknots. Therefore, $\mathbb{T}^2 \times \{2\}$ bounding another solid torus is only possible in the case of an unknot. 
Hence, the knot $c_1$ is unique up to isotopy.
\end{proof}

\begin{rem}
In the case of an embedded thickened torus in $\s^3$ determined by an unknot, there is a choice between core curves $c_1$ and $c_2$ of the solid tori $N_1$ and $N_2$ in the above proof. 
In the context of round 2-surgery, suppose $N_1$ is a solid torus that is disjoint from the rest of the round surgery link. We take the unknot $c_1$ to be the round 2-surgery knot. 
\end{rem}

For a round 2-sugery, we remove embedded thickened torus $\phi(\s^1 \times \s^1 \times \D^1)$ from $\s^3$ and glue two solid tori $\s^1 \times \D^2 \times \s^0$ using the identity homeomorphism. 

\begin{lem}\label{Lem:R2Diagram}
A round 2-surgery on $\mathbb{S}^3$ can be determined by a knot in $\s^3$ with a rational coefficient. 
\end{lem}
\begin{proof}
Suppose $\phi : \s^1 \times \s^1 \times [-1,1] \to \s^3$ is an embedding. 
By Lemma \ref{Lem:EmThickTorus}, the embedded thickened torus $\phi(\s^1 \times \s^1 \times [-1,1])$ is completely determined by a knot $K \subset \mathbb{S}^3$.
In particular, the embedded thickened torus $\phi(\s^1 \times \s^1 \times [-1,1]) = N_2 (K) \setminus N_1(K)$, where $N_i(K)=\{(\theta, (x,y)) \in \s^1 \times \R^2 \, |\,  \theta \in S^1, \,  x^2 + y^2\leq i^2 \text{ for }x, y \in \mathbb{R}\}$ is a tubular neighbourhood of the round 2-surgery knot $K$. 
Since $\phi$ is a diffeomorphism onto its image, the boundary of $\mathbb{T}^2 \times [-1,1]$ maps to the boundary $\dd (N_2(K)\setminus N_1(K))$. Without loss of generality, we can assume that $\phi\left(\mathbb{T}^2 \times \{(-1)^j\} = \dd N_j(K)\right)$ where $j=1,2$. 
On $\dd N_j$, we have the canonical coordinate system of curves $m_j$ and $l_j$. 
Moreover, the meridian $m_1$ is isotopic to $m_2$, and $l\_1$ is isotopic to $l_2$ in $N_2(K)\setminus N_1(K)$.
On the boundary of $\s^1 \times \s^1 \times [-1,1]$, curves $r_{j}=\s^1\times \{p\} \times \{(-1)^j\} $ and $s_{j} = \{q\}  \times \s^1\times \{(-1)^j\}$ determine the basis of the first homology group of $\mathbb{T}^2 \times \{(-1)^j\}$ for some fixed $p, q \in \s^1$ and $j=1,2$.
Clearly, the basis curves $r_1$ is isotopic to $r_2$ and $s_1$ is isotopic to $s_2$ in $\mathbb{T}^2 \times [-1,1]$. 
We can express diffeomorphism $\mathbb{T}^2 \times \{1\} \to \dd N_{1}$ as a matrix $A= \begin{pmatrix}
    p & r\\
    q & s\\
\end{pmatrix}$ with respect to the basis $r_1$, $s_1$ on $\mathbb{T}^2 \times \{-1\}$ and canonical coordinates $m_1$ and $l_1$ on $\dd N_1$. 

Now, we can think of $\phi: \mathbb{T}^2 \times [-1,1] \to \s^3$ as an isotopy between diffeomorphisms $\phi|_{\mathbb{T}^2 \times \{(-1)^j\}}$ onto their images for $j=1,2$. 
Recall that any isotopy class of diffeomorphism is given by the same matrix $A \in SL(2, \mathbb{Z})$.
Thus, for $j=1,2$, $\phi|_{\mathbb{T}^2 \times \{(-1)^j\}}$ is given by the same matrix $A$.

In round 2-surgery, we remove the embedded thickened torus $\phi(\s^1 \times \s^1 \times [-1,1])$ from $\mathbb{S}^3$. 
Denote $N :=\partial \overline{(\s^3\setminus \phi(\s^1 \times \s^1 \times [-1,1]))}$.
We observe that $N$ has two connected components. 
We glue two solid tori $\s^1 \times \D^2 \times \s^0 $ to each connected component of $N$ along their boundary tori by the identity diffeomorphism.
Denote $\s^0 =\{-1,1\}$. Under the identity glueing map, 
the longitude $\s^1 \times \{p\} \times \{(-1)^j\}$ maps to $\phi(r_j)$ and the meridian $\{q\} \times \dd \D^2 \times \{(-1)^j\}$ maps to $\phi(s_j)$, for $j=1,2$. 

Recall that the glueing map of a solid torus can be completely determined by the image of the meridian of the boundary torus. 
Thus, in our case, it is sufficient to know the coefficients of $\phi(s_j) = A \begin{pmatrix}
    1 & 0
\end{pmatrix}^T =p \cdot m_j + q \cdot l_j$, where $p, q \in \mathbb{Z}$ such that $(p,q)=1$. 
The quotient $\frac{p}{q}$ denotes the round 2-surgery coefficient on a $K$. For $\frac{1}{0}$-coefficient on $K$, we get $\phi(s_j)= m_j$, for $j=1,2$, yielding the resultant 3-manifold $\mathbb{S}^3 \sqcup (\s^1 \times \s^2)$. 
%Notice that the round 2-surgery on any knot in $\mathbb{S}^3$ with surgery coefficient $\frac{1}{0}$ corresponds to the $\mathbb{S}^3 \sqcup \s^1 \times \s^2$.
We can avoid this surgery coefficient because we can obtain $\mathbb{S}^3 \sqcup (\s^1 \times \s^2)$ by doing a round 2-surgery on an unknot with round 2-surgery coefficient ${0}$. 
Therefore, we get a rational number with a knot that completely determines the round surgery of index 2. 
This rational number is the {round 2-surgery coefficient}. \end{proof}

%In particular, we can get a round 2-surgery diagram consisting of an integral framed knot.
%Now, we need to understand the possible image $\phi(s_j)$ under the embedding.  
%Notice that there are two choices. 
%In the embedded $\mathbb{T}^2\times [1,2]$, the $\phi(s_j)$ is a first homology-generating curve, i.e. either $m_j$ or a longitude $f_j$ isotopic to the core $K$, for $j=1,2$.
%Thus, we have either $\phi(s_{j})= m_j$ or $\phi(s_{j})=f_j= n_{j}\cdot m_j + l_j$, for $j=1,2$ respectively. 
%Further, the images $\phi(s_{j})$ are isotopic curves on the boundaries $\partial N_2\setminus N_1$ by a linear isotopy given on the embedded thicken torus $\phi(\s^1 \times \s^1 \times \D^1)$, which implies $n_1= n_2=n $.  
%In the first choice, we get $\frac{1}{0}$ coefficient and $\frac{n}{1}$ in the second choice. 
%. 
%We can obtain this by doing a round 2-surgery on an unknot with round 2-surgery coefficient $n=0$ in the second choice. 

\begin{rem}{(Round $2$-surgery on $\mathbb{S}^3$ gives disconnected manifold.)}\label{rem: R2GivesDisconnected} 
In round 2-surgery, we remove an embedded $\mathbb{T}^2\times [0,1]$ from $\mathbb{S}^3$ and glue two solid tori $S^1 \times S^0 \times D^2$. From Lemma \ref{Lem:EmThickTorus}, we get a 3-manifold with two components after removing the $\mathbb{T}^2 \times [0,1]$. Thus, we obtain a disconnected 3-manifold after gluing solid tori. Below, we define the way to perform round 1-surgery and round 2-surgery together to get a connected 3-manifold. 
\end{rem}

\begin{rem}
    By the definition of round surgery, it is clear that a round surgery on a closed manifold produces a closed manifold. Thus, we get a closed 3-manifold after performing round 1 and 2 surgeries on $\mathbb{S}^3$. Moreover, we fix an orientation on $\s^3$ to define the coordinate curves. Thus, we get an oriented smooth 3-manifold  after round surgeries on $\s^3$. 
\end{rem}

\begin{proof}[proof of the theorem \ref{thm: ARSequalsRS}]
   From the above Lemmas \ref{Lem:R1Diagram} and \ref{Lem:R2Diagram}, it is clear that Asimov's round surgery can be described using our description of round surgery. Moreover, it is also easy to see that the round surgery on framed links in $S^3$ can be realised as the round surgery performed on a pair of solid tori or on a thickened torus, depending on the index of the surgery under consideration. As discussed in the previous Lemmas 2 and 4, the glueing maps in both cases of the round surgery on links are the identity homeomorphisms. Thus, both of the round surgery descriptions become equivalent on $\s^3$.
\end{proof}

The preceding theorem establishes that on $\s^3$, Asimov's round surgery is equivalent to surgery performed on a framed link. Consequently, this yields combinatorial descriptions of 3-manifolds obtained via round surgeries on $\s^3$. It is natural to inquire whether this equivalence between Asimov's round surgery and surgery on framed links extends to an arbitrary 3-manifold $M$. In particular, we pose the following question:

\begin{ques}
    Let $N$ be the 3-manifold obtained by a round surgery of index-1 on an embedded $\s^1 \times \s^1 \times \s^0$ in a 3-manifold $M$. Can we realise this round surgery of index-1 as being performed on some framed link $L \subset M$?
\end{ques}

A natural candidate for the link $L$ consists of the core curves of the solid tori. However, this approach is obstructed by the absence of a canonical reference framing. To circumvent this difficulty, we may restrict our attention to cases where the ambient 3-manifold is simply connected (which is $\s^3$), or consider round surgeries performed on tubular neighbourhoods of null-homologous links. In particular, we want to know the necessary conditions on the ambient 3-manifold such that the round surgery of index 1 can be performed on a framed link.

A similar inquiry naturally arises for round surgeries of index-2:

\begin{ques}
    Suppose $N$ is the 3-manifold obtained by performing a round surgery of index-2 on a thickened torus in a 3-manifold $M$. Under what conditions on $M$ does there exist a (rationally) framed knot $K \subset M$ such that a round surgery of index-2 on $K$ yields $N$?   
\end{ques}

In the rest of the article, we study the round surgery diagrams on $\s^3$. Thus, we restrict ourselves to round surgeries on $\s^3$.  

\subsection{Joint pair of round surgeries of indices 1 and 2}\label{subsec: Joint pair}
By Remark \ref{rem: R2GivesDisconnected}, we know that round surgery of index 2 on a knot produces a 3-manifold with two connected components. 
In particular, removing an embedded thickened torus from $\s^3$ gives a two-component 3-manifold. 
By Lemma \ref{Lem:EmThickTorus}, we know that one component is the solid torus, and the other component is the knot complement of the core curve of the solid torus.   
We need to know possible embeddings for round surgeries on $\s^3$ such that we don't get two components after performing surgeries or a combination of surgeries. 
From Lemma \ref{Lem:EmThickTorus}, we know that every thickened torus is part of some solid torus.
If this solid torus is a part of some round 1-surgery, then we get a connected manifold after performing round 1-surgery along that solid torus and removing that embedded thickened torus. 
In the following Lemma \ref{lem: RemovedTorusNbd}, we prove it. 

Now, let $L=L_{11} \cup L_{12}$ in $\s^3$ be a framed link for round surgery of index $1$. Let $\phi(\mathbb{T}^{2} \times I) $ be an embedded thickened torus such that $\phi(\mathbb{T}^{2} \times\{0\})$ is parallel to $\overline{\partial\left(\s^{3} \backslash N\left(L_{12}\right)\right)}$, i.e., $\phi(\mathbb{T}^2 \times I)= N_2(L_{12}) \setminus N_1(L_{12})$. Then perform a round 1-surgery on $L=L_{11} \cup L_{12}$, i.e.,
$$
M_{L}:=\overline{\s^{3} \backslash \psi\left(\s^{1} \times \s^{0} \times \D^{2}\right)} \bigcup_{\s^{1} \times \s^{0} \times \s^{1}} \s^1 \times \D^{1} \times \s^{1}
$$
where $\psi : \s^1 \times \s^0 \times \D^2 \to \s^3$ is an embedding. 
Since $\phi(\mathbb{T}^{2} \times I) \subset \s^{3} \backslash N(L),$ we get $\phi(\mathbb{T}^{2} \times I )\subset M_{L}$. 
We want to remove embedded $\phi(\mathbb{T}^{2} \times I)$ from $M_{L}$. 
Denote $M_{L}^{\prime}=\overline{M_{L} \backslash\left(\phi(\mathbb{T}^{2} \times I)\right)}$.
Clearly, $\partial M_{L}^{\prime}=\phi(\mathbb{T}^2 \times \s^0) =\phi(\mathbb{T}^{2} \times \{0\}) \cup \phi(\mathbb{T}^{2} \times \{1\})$.

\begin{lem}\label{lem: RemovedTorusNbd}
Removing the above embedded $\phi(\mathbb{T}^{2} \times I)$ (parallel to the glued torus in round 1-surgery) nullifies the effect of the glued thickened torus in round 1-surgery, i.e.
$$
\overline{M_{L} \backslash\phi\left(\mathbb{T}^{2} \times I\right)} \cong \overline{\s^{3} \backslash N(L)},
$$    
where $L=L_{11} \cup L_{12}$ and $M_{L}$ described above.
\end{lem}

\begin{proof}
We have given that the embedded $\mathbb{T}^{2} \times I$ is parallel to the glued torus in round 1-surgery. 
In $M_{L}^{\prime}$, the boundary torus $\phi (\mathbb{T}^2 \times \{0\})$ is given to be parallel to the $\dd N(L_{12})$ and the other boundary torus $\phi (\mathbb{T}^2 \times \{1\})$ can be isotoped along the glued thickened torus of round 1-surgery to realise it as parallel to the $\dd N(L_{11})$.   
Thus, $ M_{L}^{\prime} \cong \overline{M_{L} \backslash \phi\left(\mathbb{T}^{2} \times I\right)} \cong \overline{\s^{3} \backslash N(L)}$.
\end{proof}

We observe that performing a round 2-surgery on a thickened torus parallel to the glued torus of round 1-surgery produces a connected 3-manifold. 
Notice that we are performing a round 2-surgery on the thickened torus $N_2(L_{12})\setminus N_1(L_{11})$, which is a part of the tubular neighbourhood of a link component $L_{12}$ of the round 1-surgery link $L$. 
Thus, the component $L_{12}$ is a round 2-surgery knot as well. 
We call such a pair of round 1-surgery and round 2-surgery on the round 1-surgery link a {\it joint pair of round surgeries of indices 1 and 2}. 
Formally, we define it as follows. 

\begin{figure}[ht]
    \centering
    \includegraphics[scale=0.4]{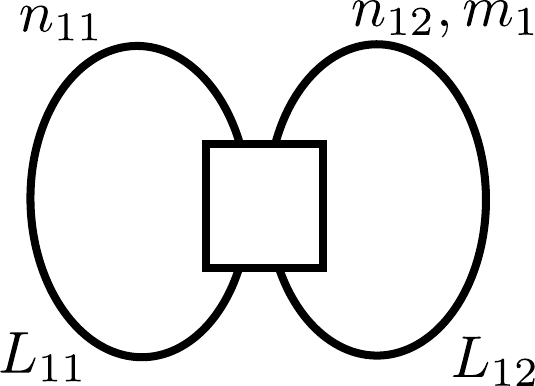}
    \caption{A diagram of a Joint pair $L_{11} \cup L_{12}$ with round 1-surgery coefficient $n_{11}$ on $L_{11}$ and $n_{12}$ on $L_{12}$, and round 2-surgery coefficient $m_1$ on $L_{12}$. The box in the middle represents the linking between the components and the knotting of $L_{11}$ and $L_{12}$.}
    \label{fig: JP}
\end{figure}

\begin{defn}
    A round 1-surgery link $L_{11} \cup L_{12}$ is said to be a {\it joint pair of round surgeries of indices 1 and 2} if one of the components of $L_{11} \cup L_{12}$ is treated as a round 2-surgery knot. 
    We denote the coefficient of the round 2-surgery knot on the top of that component next to the round 1-surgery coefficient as shown in Figure \ref{fig: JP}.
\end{defn}

\begin{rem}\label{rem: StdPos}
%    We say a joint pair $L_{11} \cup L_{12}$ in a {\it standard position} if the second component $L_{12}$ of round 1-surgery link $L_{11} \cup L_{12}$ is also a round 2-surgery knot.
   % From this point onward, for a given joint pair, we will assume that we are performing round 2-surgery on the second component $L_{i2}$
%We fix a convention of indexing the link component as $L_{12}$ round 2-surgery knot for a joint pair of round surgery diagrams $L=L_{11} \cup L_{12}$. 
%Further, for a joint pair, we denote the round 1-surgery coefficient by $n_{11}$ and $n_{12}$, and the round 2-surgery coefficient $m_1$ next to the $n_{12}$ on top of $L_{12}$ as shown in the figure \ref{fig: JP}. 

We fix a convention to index a joint pair $L$ as $L_{i1} \cup L_{i2}$ such that $L_{i1} \cup L_{i2}$ is a round 1-surgery link with round 1-surgery coefficient $n_{i1}$ on $L_{i1}$ and $n_{i2}$ on $L_{i2}$, and $L_{i2}$ is also a round 2-surgery knot with round 2-surgery coefficient $m_i$, for some $i\in \mathbb{N}$. 
\end{rem}

\subsection{A round surgery diagram in $\mathbb{S}^3$ of a connected closed 3-manifolds}

%We know that round 2-surgery produces a 3-manifold with two components. Thus, a round surgery diagram for a closed connected 3-manifold has a round 2-surgery knot in a joint pair with some round 1-surgery link. 

Suppose $M$ is a 3-manifold obtained by a sequence of round surgeries of indices 1 and 2 on $\s^3$.
For such a sequence of round surgeries, we can associate a framed link corresponding to the round surgeries of indices 1 and 2. We call such a framed link a {\it round surgery diagram}. 
\begin{defn}
    A link $L=L_1 \cup \cdots \cup L_{n}$ is said to be a round surgery diagram of a 3-manifold $M$, if 
    \begin{enumerate}
        \item for each $i$, $L_i$ is either a round 1-surgery link or a round 2-surgery knot. 
        \item for some round 1-surgery link $L_i$, $L_{i2}$ is also a round 2-surgery knot, i.e. a $L_i$ is a joint pair. 
    \end{enumerate}
\end{defn}

We know that round 2-surgery produces a $3$-manifold with two components. 
If $M$ is connected, then a round 2-surgery knot must be in a joint pair with a round 1-surgery link; otherwise, we get a disconnected resultant manifold.  
Below, we define a round surgery diagram for a closed, connected, oriented 3-manifold.  
\begin{defn}
    A round surgery diagram in $\mathbb{S}^3$ of a closed connected oriented 3-manifold is a link $L=L_{1} \cup \ldots \cup L_{k}$ such that 
    \begin{enumerate}
        \item for each $1\leq i\leq n$, $L_i$ is a link with two components with an integral framing corresponding to the round 1-surgery,
        \item and for some of the $1 \leq i \leq n$, $L_i$ is a joint pair of the round surgeries of indices 1 and 2. 
    \end{enumerate}
\end{defn}

\subsection{Trivial round surgery diagram}
In the proof of Lemma \ref{Lem:R2Diagram}, we have a choice of $1/0$-round 2-surgery coefficient. Let us consider a joint pair with a $1/0$-round 2-surgery coefficient.
Suppose $L= L_{11} \cup L_{12}$ is joint pair with round 1-surgery coefficient $n_{11}$ on $L_{11}$ and $n_{12}$ on $L_{12}$ and round 2-surgery coefficient $1/0$ on $L_{12}$. 
After round 1-surgery and removing the thickened torus for round 2-surgery from $\s^3$, we are left with $\overline{\s^3 \setminus N (L)}$ as proven in the Lemma \ref{lem: RemovedTorusNbd}. 
Now, gluing two solid tori with the $1/0$ round 2-surgery coefficient is the same as performing two trivial Dehn surgeries together. 
Thus, the resultant manifold is again $\s^3$. 
We call such a round surgery a trivial round surgery. 
In particular, we say a round surgery diagram {\it trivial} in $\mathbb{S}^3$ if the round 2-surgery coefficient on each joint pair $L_{j2}$ is $1/0$. 

\subsection{Connected sum of 3-manifolds} 

In this subsection, we observe that two round surgery links contained in disjoint 3-balls in $\s^3$ yield the connected sum of the manifolds obtained by performing round surgery on each link individually. 

Suppose two round surgery links $L_1$ and $L_2$ are contained in two disjoint $3$-balls in $\s^3$. 
For each $1\leq i\leq 2$, we denote by $M_{L_i}$ the 3-manifold obtained after performing round surgery on $L_i$. 
We can think of $\s^3$ as a connected sum of its two copies, i.e., $\s^3_1\#\s^3_2$, such that $L_1 \subset \s^3_1$ and $L_2\subset \s^3_2$. 
Then, it is easy to see that we obtain the connected sum $M_{L_1}\#M_{L_2}$ after performing round surgeries on $L_1$ and $L_2$. 
Moreover, let $K_{ij}$ denote a round 2-surgery knot in the link $L_i$ such that $K_{ij}$ is not in a joint pair with any round 1-surgery link.
Since round surgery on a round 2-surgery knot yields a manifold with two components, $M_{L_i}$ has disconnected components $M_{K_{ij}}$ obtained from glueing tubular neighbourhoods of $K_{ij}$ with a solid torus and $M'_{L_j}$ obtained by the rest of the round surgeries. 
Then, the resultant 3-manifold can be expressed as follows: 
$$M_{L_1\sqcup L_2}= \bigsqcup_{j=1}^{k_1} M_{K_{1j}} \bigsqcup \left(M'_{L_1} \#M'_{L_2}\right) \bigsqcup_{j=1}^{k_2} M_{K_{2j}}. $$

%We know that the complement of the tubular neighbourhood of the surgery link remains unchanged after a round surgery over a link. Thus, we get a connected sum of the manifolds after performing round surgery on the given round surgery links. 

%Suppose two round surgery links are contained in two disjoint 3-balls in $\mathbb{S}^3$. We know that the complement of the tubular neighbourhood of the surgery link remains unchanged after a round of surgery over a link. Thus, we get a connected sum of the manifolds after performing round surgery on the given round surgery links. 

\section{Bridge between Dehn surgery diagrams and round surgery diagrams}\label{sec: BridgeThmC4}

We aim to convert a given round surgery diagram for round surgery of indices 1 and 2 into a Dehn surgery diagram. 
In particular, we want to establish a bridge between round surgery diagrams and Dehn surgery diagrams.

\subsection{Correspondence between round surgery diagrams and surgery diagrams.}
Suppose the link $L=L_{11} \cup L_{12}$ is a joint pair of round surgeries of indices 1 and 2. 
Then we perform round 1-surgery on the tubular neighbourhood $N(L)(:= N(L_{11}) \cup N(L_{12}))$ of the link $L$ and round 2-surgery on a thickened torus $N_2(L_{12})\setminus N_1(L_{12})$ parallel to the $\dd N(L_{12})$. 
From Lemma \ref{lem: RemovedTorusNbd}, we know that after performing round 1-surgery and removing the thickened torus, we get $\overline{\s^3 \setminus N(L)}$. 
We glue two solid tori to complete the round 2-surgery. 
Notice that, glueing $\s^1 \times \D^2 \times \s^0$ to $\overline{\s^3\setminus N(L)}$ is same as performing Dehn surgery on $L$. 
In Lemma \ref{lem: JPtoOSD}, we prove this and calculate the Dehn surgery coefficients. 

\begin{figure}[ht]
        \centering
        \includegraphics[scale=0.4]{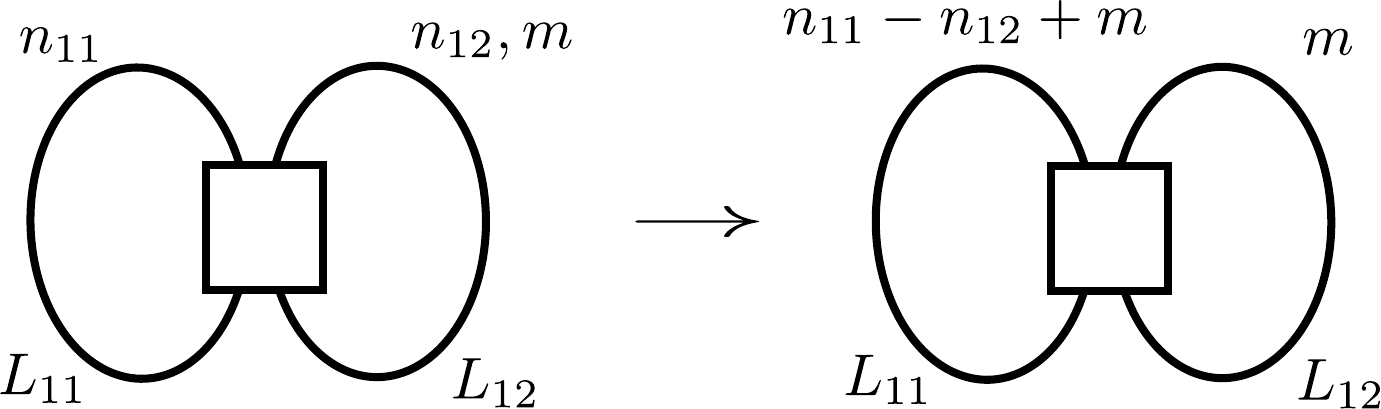}
        \caption{A joint pair to a Dehn surgery diagram.}
        \label{fig:JPtoOSD}
    \end{figure}

Moreover, any Dehn surgery on some two-component link can be viewed as performing round surgeries on a joint pair.
This observation provides a correspondence between the set of round surgery diagrams of joint pairs and integral Dehn surgery diagrams. 
\begin{lem}[{\bf Joint pair to a Dehn surgery diagram}]\label{lem: JPtoOSD}
    Let $L=L_{11} \cup L_{12}$ be a round surgery diagram of a closed connected oriented 3-manifold such that $L$ is a joint pair of round surgeries of indices 1 and 2. Then, it corresponds to the Dehn surgery link $L$ with the surgery coefficient shown in Figure \ref{fig:JPtoOSD}.
\end{lem}
\begin{proof}
Suppose we perform round 1-surgery on $L_{11} \cup L_{12}$ with round surgery coefficient $n_{11}$ on $L_{11}$ and $n_{12}$ on $L_{12}$. 
For round 1-surgery, the solid tori are embedded as the tubular neighbourhoods $N(L_{11}) \sqcup N_1(L_{12})$ of the $L_{1j}$ in $\s^3$. 
We consider a larger radius tubular neighbourhood $N_2(L_{12})$ of $L_{12}$ such that $N_1 (L_{12}) \subset N_2(L_{12})$. 
We also perform a round 2-surgery on $L_{12}$ on the embedded thickened torus $N_2(L_{12}) \setminus N_1 (L_{12})$ in $\s^3$.

On the boundary torus $\partial N(L_{1j})$, we have the canonical coordinate system given by the curves meridian $m_j$ and longitude $l_{j}$ for $j=1,2$. 
%We use the notation from the proof of Lemma \ref{Lem:R1Diagram} to describe glueing. 
%In round 1-surgery, we glue $\s^1 \times \s^1 \times [1,2]$ to $\dd N(L_{1j})$ such that 
%\begin{eqnarray*}
 %     & r_{j} \mapsto m_{j}\\
 %      \text{ and } &  s_{j} \mapsto n_{1j} \cdot m_{j} + l_{j}.
%\end{eqnarray*}
We get $N: = \overline{\s^3\setminus \psi(\s^1\times \D^2\times \s^0)} \cup \s^1 \times \s^1 \times [1,2]$ after round 1-surgery, where $\psi : \s^1 \times \D^2 \times  \s^0 \to \s^3$ is an embedding.
   
On $L_{12}$, we perform a round 2-surgery with respect to framing $m$.
After removing the embedded thickened torus from $N$, we get $N^{\prime}= N\setminus (N_2 (L_{12}) \setminus N_1 (L_{12}))$. 
On $\partial N^{\prime} =\dd N_2 (L_{12}) \sqcup \dd N_1 (L_{12})$, we have coordinate curves $m_{2}$ and $l_{2}$. 
The round 2-surgery framing curve on each $\dd N_i (L_{12})$ can be expressed as $f_i= p\cdot m_{2} + q \cdot l_{2}$ for $i=1,2$, because $m = \frac{p}{q} \in \mathbb{Q}$. 

Notice that the boundary component $\dd N_1(L_{12})$ %$ T^2 \times \{2\}$ 
can be realised as the boundary of a tubular neighbourhood of $L_{11}$. 
We observe that gluing a solid torus to $\dd N_2(L_{12})$ is the same as performing a Dehn surgery on $L_{12}$ by surgery coefficient $m$. 
Thus, to realise this round 2-surgery as Dehn surgeries, we also need to express the framing curve $f_1$ in terms of the coordinate curves $m_1$ and $l_1$ of $\partial N(L_{11})$.
Suppose $\mu$ and $\lambda$ are the coordinate curves on the glued thickened torus. 
Then on $\partial N_1(L_{12})$, $\lambda \mapsto n_{12} \cdot m_2 + l_2$ and $\mu \mapsto m_2$, and on $\partial N(L_{11})$, $\lambda \mapsto n_{11} \cdot m_1 + l_1$ and $\mu \mapsto m_1$. 
Thus, we conclude that $m_1$ and $m_2$ are isotopic and  
\begin{eqnarray*}
        l_2&= &\lambda - n_{12} \cdot m_2\\
        &= & n_{11} \cdot m_1 + l_1  + n_{12}\cdot m_2\\
        &= & (n_{11}-n_{12}) \cdot m_1 + l_1, \quad \because \, m_1=m_2. 
\end{eqnarray*}
Therefore, 
\begin{eqnarray*}
        f_1 &= & p\cdot m_2 + q \cdot\{(n_{11}-n_{12})\cdot m_1 +l_1\}\\
        &=& (q(n_{11}-n_{12}) +p) \cdot m_1 + q \cdot l_1.
\end{eqnarray*}
Thus, the Dehn surgery coefficient is $(n_{11}-n_{12}+p/q) = (n_{11}-n_{12}+m)$ on $L_{11}$ and $m$ on $L_{12}$ for a given joint pair of round surgery diagram. 
\end{proof}

\begin{defn}
A round surgery diagram $L = L_{1} \cup \ldots\cup L_{n}$ is said to be {\it round surgery diagram of joint pairs} if for each $1\leq i\leq n$
\begin{enumerate}
    \item $L_i$ is a link with two components, i.e. $L_i = L_{i1} \cup L_{i2}$
    \item $L_{i1} \cup L_{i2}$ is a joint pair of round surgeries of indices 1 and 2 such that $n_{ij}$ is a round 1-surgery coefficient on $L_{ij}$ and $m_j$ is a round 2-surgery coefficient on $L_{i2}$, where $1\leq j\leq 2$. 
\end{enumerate} 
\end{defn}

\begin{rem}{(Restriction to integer coefficients to round 2-surgery knots.)}It is not always possible to convert a given Dehn surgery diagram with rational coefficients to a joint pair. 
    For example, we take a pair of unlinked unknots $U_1 \cup U_2$ with rational coefficients $\frac{p}{q}$ on $U_1$ and $+1$ on $U_2$. We claim that this rational Dehn surgery diagram can not be converted to a joint pair for round surgery. To prove the claim, assume that it is possible to convert it into a joint pair of round surgery diagram with round 1-surgery coefficient $n_1$ on $U_1$ and $n_2$ on $U_2$ for some integers $n_1,n_2$ and round 2-surgery coefficient $m$ on $U_2$. By Lemma \ref{lem: JPtoOSD}, we can convert this joint pair back to the Dehn surgery diagram with the Dehn surgery coefficients $n_1 -n_2 +m$ on $U_1$ and $m$ on $U_2$. 
    It implies, $m=1$ and $n_1 -n_2 +m = \frac{p}{q}$, which is contradiction because $n_1, n_2 \in \mathbb{Z}$.
\end{rem}

 %We are interested in studying the equivalence moves on the round surgery diagrams similar to Kirby moves on integral Dehn surgery diagrams. 
In view of the above remark, we will only consider integral round 2-surgery coefficients from now on.  
 
\begin{lem}[{Ordinary surgery diagram to a joint pair of round surgeries}]\label{lem: OSDtoJPRS}
Suppose $L=L_1 \cup \cdots \cup L_n$ is a Dehn surgery diagram such that component $L_i$ has a surgery coefficient $n_i$ for $1 \leq i\leq n$. 
If $n$ is even, say $n =2m$, then $(L_1 \cup L_2)\cup  \cdots \cup (L_{2m-1} \cup L_{2m})$ is a round surgery diagram of joint pairs such that each pair $L_{2i-1} \cup L_{2i}$ in the parentheses is a joint pair with round 1-surgery coefficients $n_i -n_{i+1}+k $ on $L_{i}$ and $k$ on $L_{i+1}$, and round 2-surgery coefficient $n_{i+1}$ on $L_{i+1}$, for $1 \leq i \leq m$ and $k\in \mathbb{Z}$. 

If $n$ is odd, say $n=2m-1$, then we add an unlinked unknot $L_{2m}$ with framing $(\pm1)$ to convert $L$ into a new Dehn surgery diagram with an even number of components and apply the previous statement.
\end{lem}

\begin{proof}
Suppose we have given a surgery diagram consisting of a link $L= L_1 \cup L_2$ with two components such that $n_j$ is a surgery coefficient on $L_j$ for $1\leq j\leq 2$. 
\begin{figure}[ht]
    \centering
    \includegraphics[scale=0.4]{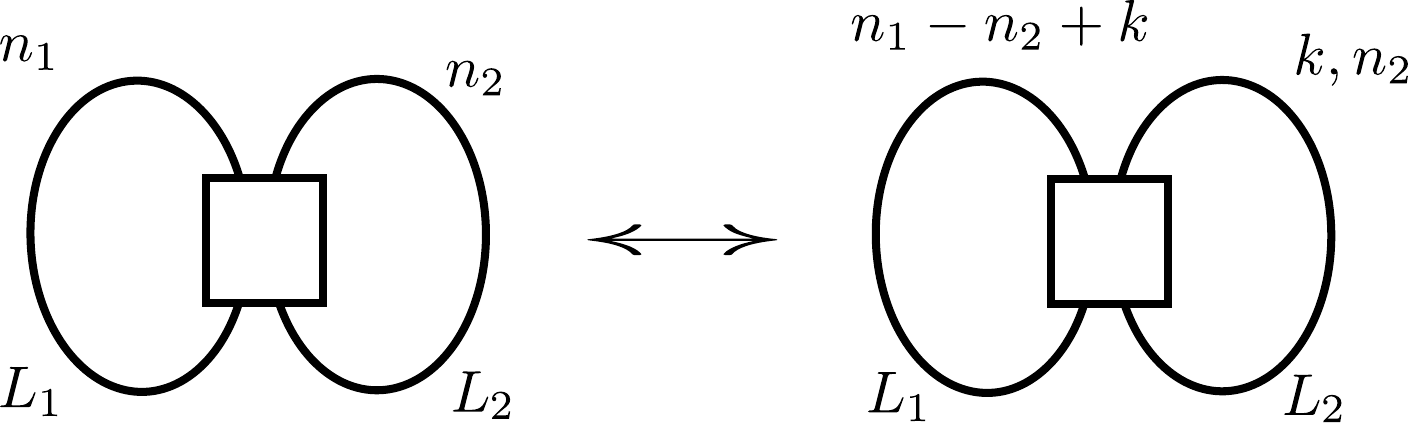}
    \caption{An ordinary surgery diagram consists of two components to a joint pair round surgery diagram on the left.}
    \label{fig: 2ComOSDtoJPRSD}
\end{figure}
Then, we can realise them as a joint pair $L_1 \cup L_2$ such that the pair has the round 1-surgery coefficients $n_1-n_2+k$ on $L_1$ and $k$ on $L_2$, and round 2-surgery coefficient $n_2$ on $L_2$, where $k$ can be any integer (see Figure \ref{fig: 2ComOSDtoJPRSD}).
We observe that we can convert the joint pair $L_1 \cup L_2$ to the Dehn surgery diagram using Lemma \ref{lem: JPtoOSD}. 
We obtain the same Dehn surgery diagram $L_1 \cup L_2$ with surgery coefficient $n_j$ on $L_j$ for $j=1,2$.
Thus, the joint pair and Dehn surgery diagram correspond to the same 3-manifold. 

Suppose we have a surgery diagram consisting of a single component $K$ with surgery coefficient $n$. 
We need two components to convert a Dehn surgery diagram into a joint pair.
Thus, we perform the Kirby move of type 1, i.e., addition of an unlinked disjoint unknot with $\pm1$ surgery coefficient to $K$ to get $K \cup U$. 
Now, we convert the new Dehn surgery diagram to the joint pair as discussed above (see Figure \ref{fig: OSDtoJP}). 

\begin{figure}[ht]
    \centering
    \includegraphics[scale=0.4]{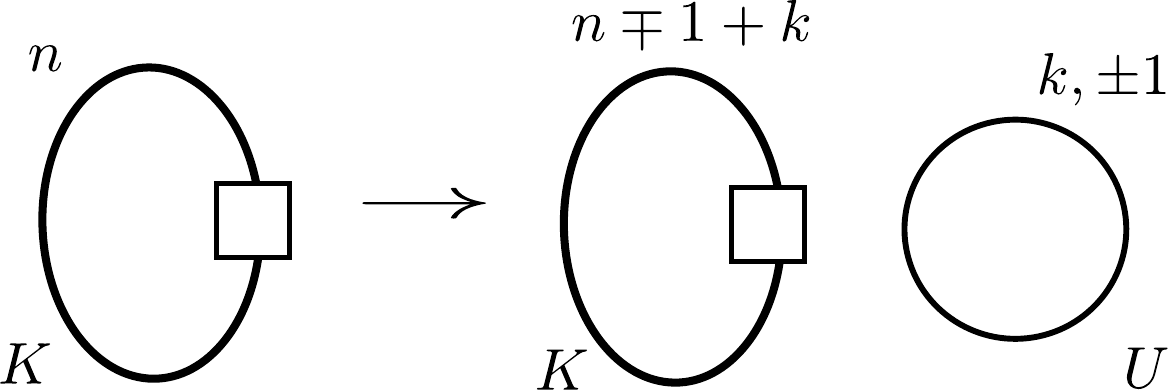}
    \caption{A joint pair from a Dehn surgery diagram consisting of a single component.}
    \label{fig: OSDtoJP}
\end{figure}
Suppose we have given a Dehn surgery diagram $L= L_1 \cup \cdots \cup L_n$. 
We use the order of the link components to combine two components to form a candidate for a joint pair. 
For each pair of link components, we treat them as a Dehn surgery diagram consisting of two components and then convert them into a joint pair, as mentioned above.  
If $n$ is odd, then we treat the last component $L_n$ as a Dehn surgery diagram consisting of a single component and convert it into a joint pair.
In conclusion, we obtain a round surgery diagram of joint pairs from the given Dehn surgery diagram. 
\end{proof}

\begin{rem}\label{rem: AmbiguityInOSDtoJP}
    There is a choice involved in the proof of the above lemma about which component to assign a round 2-surgery curve. 
    However, it does not matter because the resultant 3-manifolds will be the same (from Lemma \ref{lem: JPtoOSD}). 
    In Section \ref{Sec : KCforRSDJP}, we give an equivalence relation to address this issue. 
\end{rem} 

From Lemma \ref{lem: JPtoOSD} and Lemma \ref{lem: OSDtoJPRS}, we have proved Theorem \ref{BridgeTheorem}: the Bridge Theorem. 

\iffalse \begin{thm}[{The Bridge Theorem}]\label{BridgeTheorem}
The following two statements establish a bridge between a set of integral Dehn surgery diagrams and round surgery diagrams. 

\begin{enumerate}
    \item[(a)] Given a round surgery diagram of joint pairs $L=\bigcup_{i=1}^{n} (L_{i1} \cup L_{i2})$ in $\s^3$ of closed connected oriented smooth 3-manifold $M_{L}$ such that $n_{ij}\in \mathbb{Z}$ is a round 1-surgery coefficient on $L_{ij}$ and $m_i\in \mathbb{Z}$ is a round 2-surgery coefficient on $L_{i2}$ for $1\leq i\leq n$ and $1 \leq j \leq 2$. Then there is a Dehn surgery diagram given on $L$ such that $n_{i1}-n_{i2}+m_i$ and $m_i$ are surgery coefficients on $L_{i1}$ and $L_{i2}$, respectively.
    \item[(b)] Given an integral Dehn surgery diagram $L= L_1 \cup \cdots \cup L_n$ of closed connected oriented smooth 3-manifold such that the integer $n_i$ is a surgery coefficient on $L_i$ for each $1\leq i\leq n$. Then there is a round surgery diagram of joint pairs.
\end{enumerate}
\end{thm}\fi

As a corollary, we have Asimov's theorem for connected, oriented 3-manifolds. 

\begin{proof}[Proof of Corollary \ref{cor: AsimovResult}]
    By Lickorish--Wallace theorem, given a closed connected oriented smooth 3-manifold $M$, there is a Dehn surgery diagram $L$ consisting of a framed link $L$. 
   There is a round surgery diagram of joint pairs for this surgery diagram $L$ from Theorem \ref{BridgeTheorem}. 
\end{proof}

\begin{cor}%[Lickorish--Wallace theorem for round surgery diagrams] \label{}
Given a round surgery diagram of joint pairs $L=\bigcup_{i=1}^{n}(L_{i1} \cup L_{i2})$, we can obtain a new round surgery diagram of joint pairs $L'$ from $L$ such that all the round surgery coefficients are either of $\pm1$. 
\end{cor}
\begin{proof}
    Suppose we have given a round surgery diagram of joint pairs $L= \bigcup_{i=1}^{n}(L_{i1}\cup L_{i2})$ with round 1-surgery coefficient $n_{ij}$ on $L_{ij}$ and round 2-surgery coefficient $m_i$ on $L_{i2}$ for $1 \leq i\leq n$ and $1\leq j\leq2$. 
    By Theorem \ref{BridgeTheorem}, we convert $L$ to a Dehn surgery diagram $L$ with Dehn surgery $n_{i1}-n_{i2}+ m_i$ and $m_i$ on $L_{i1}$ and $L_{i2}$ for $ 1\leq i \leq n$. 
    Recall that any Dehn surgery diagram with an integral surgery coefficient can be converted into a Dehn surgery diagram $L'$, with each component having a surgery coefficient either $\pm1$. 
    If $L'$ has an odd number of components, then we add an unlinked disjoint unknot with either of $\pm1$ surgery coefficient. 
    Suppose $L' = \bigcup_{i=1}^{m}(L'_{i1} \cup L'_{i2})$ has even components with $n'_{ij}$ is a  suregry coefficient on $L'_{ij}$ such that $n'_{ij} \in \{\pm1\}$ for $1 \leq i\leq n$ and $1\leq j\leq 2$. 
    By Theorem \ref{BridgeTheorem}, we convert $L'$ to a round surgery diagram $\bigcup_{i=1}^{m}(L'_{i1}\cup L'_{i2})$ with round 1-surgery coefficient $n'_{i1}-n'_{i2}+k_i$ on $L'_{i1}$ and $k_i$ on $L'_{12}$ and round 2-surgery coefficient $n'_{i2}$ on $L'_{i2}$ where $k_i\in \mathbb{Z}$ for $1\leq i\leq n$. 
    For each $i$, we take $k_i=n'_{i2}$. 
    Hence, we get a round surgery diagram of joint pairs such that round surgery coefficients are either of $\pm1$.
\end{proof}

We end this subsection with the following definition of the \textit{round surgery diagram of index 1}. 

\begin{defn}
A round surgery diagram \( L = L_1 \cup \ldots \cup L_n \subset S^3 \) for a 3-manifold is called a \emph{round 1-surgery diagram} if the following conditions hold for each \( 1 \leq i \leq n \) and \( 1 \leq j \leq 2 \):
\begin{enumerate}
    \item Each \( L_i \) consists of exactly two components, \( L_{i1} \cup L_{i2} \), forming only a \emph{round 1-surgery link}.
    \item Each component \( L_{ij} \) is assigned a \emph{round 1-surgery coefficient} \( n_{ij} \in \mathbb{Z} \).
\end{enumerate}
\end{defn}

\subsection{From a Kirby diagram to a round 1-surgery diagram}
Suppose two ordinary handles $h_j$ and $h_{j+1}$ are independently attached to an $n$-manifold $M$.
Thus, the attaching regions and the attaching spheres of handles are disjoint from each other. 
In particular, in the case of dimension 4, suppose $X$ is a 4-manifold obtained by independently attaching a 1-handle $h_1$ and a 2-handle $h_2$ to $\D^4$. 
Then, the Kirby diagram of $X$ consists of a disjoint pair of $\D^3$ and a framed knot $K$.
In the lemma below, we will use the proof of the fundamental lemma of a round handle by Asimov in \cite{Asimov} to obtain a round 1-surgery diagram from a given Kirby diagram consisting of a framed knot disjoint from a pair of 3-balls.   

\begin{lem}[{Kirby diagrams to round 1-surgery diagrams}]\label{Lem: OSDtoR1SD}
    Suppose $\mathcal{A}$ is a Kirby diagram consisting of a pair of $\D^3$ and an integral framed knot $K$ such that $K$ is disjoint from the pair of $\D^3$. 
    Then $\mathcal{A}$ can be isotoped to a round 1-surgery diagram. 
\end{lem}

\begin{proof}
Suppose $\dd_{+} h_{j}$ and $\dd_{-}h_j$ denote the belt region and attaching region of the handle $h_j$, respectively. 
In $\mathcal{A}$, a pair of $\D^3$ determines the attaching sphere of a 1-handle, and a framed knot $K$ determines the attaching sphere of a 2-handle on $\s^3$.
We may choose a tubular neighbourhood $N(K)$ of $K$ such that $N(K)\cap \partial_{-}h_1 = \emptyset$. 
Since $\partial_{-}h_2= N(K)$, we see that the attaching regions are disjoint, i.e. $\dd_{-}h_1\cap \dd_{-}h_2= \emptyset$. 
Therefore, we have a pair of independently attached handles $h_1$ and $h_2$. 
Thus, we may isotope the attaching map of $h_2$ as discussed in the proof of a fundamental theorem of the round handle in \cite{Asimov}. 

In particular, we have a smooth map $F: \s^1 \times [0,1] \to \s^3 \cup \partial_+ h_1$ such that $F_0(\s^1)$ is disjoint from $h_1$ and $F_1(\s^1) \cap h_1 = \D^1 \times \s^0$. 
Moreover, this isotopy can be extended to the attaching region, i.e.,  $F' : (\s^1\times \D^2)\times [0,1]\to \s^3 \times \partial_+h_1$ such that $F'_0(\s^1\times \D^2) \cap h_1 = \emptyset $ and $F'_1 (\s^1 \times \D^2)\cap h_1 = (\D^1 \times \D^2) \times \s^0$. 

We know that the framing of the attaching sphere $\s^1 \times \{0\} $ of $h_2$ is a choice of longitude on the boundary of $\s^1 \times \D^2$ with respect to the trivialisation defined by the Seifert surface. 
This framing can be determined by an integer $n$, and it is also the surgery coefficient of the 2-surgery.  

Since the isotopy $F$ perturbs an arc of the $\s^1$ and keeps the rest of the arc identical, we may assume that the perturbing arc has Seifert trivialisation. 
After the extended isotopy $F'$, we can realise $F'(\s^1\times \D^2)= ({\D_1}^1\times \D^2) \cup ({\D_2}^1\times \D^2 \times \s^0) \cup ({\D_4}\times \D^2)$. 
By the extended isotopy $F'$, the thickened arcs $\D_1^1 \times \D^2 $ and $\D_4^1 \times \D^2$ lies in $\s^3$ and $\D^{1}_2\times \D^2 \times \s^0 $ lies on the 1-handle. 

Further, the arc $\D_1^1 \times \{0\} $ is a part of the perturbed arc of the 2-surgery knot $K$, and the arc $\D_j^1 \times \{0\} $ is the unperturbed arc. 
While realising the pair of handles $h_1$ and $h_2$ as a round handle, we join the end of the arc $\D_1^1 \times \{0\} $ by $\D^1 \in \s^3$ to get a knot $L_1$ in $\s^3$. Notice that $L_1$ is an unknot with 0-framing. 
Similarly, we join the end of the other arc $\D_1^1 \times \{0\}$ by $\D^1 \in \s^3$ to get another knot $L_2$. Clearly, the knot types of $L_2$ and $K$ are the same, and the framing of $L_{2}$ is the same as that of $K$.  
As a result, we get a link with two unlinked components $L_{1}$ and $L_{2}$ with framing $0$ on $ L_{1} $ and $n$ on $L_{2}$.
Thus, we get the corresponding round 2-surgery diagram, a two-component framed link with framing $0$ and $n$. 
\end{proof}

\begin{rem}(Round surgery diagrams of joint pairs are the optimal choice to define the Equivalence theorem).
By Lickorish--Wallace theorem, there is an integral surgery diagram $\mathcal{K}_0$ for a circle bundle over a torus $M$ obtained in Example \ref{ex: R1S-3-torus} via round 1-surgery. 
As discussed in Example~\ref{ex: R1SD-3torus}, we also have a round 1-surgery diagram of $M$ on a Hopf link.  
By Lemma~\ref{Lem: R1SDtoOSD}, we get the Kirby diagram with a 1-handle and two 2-handles $\mathcal{K}$ for the 4-manifold whose boundary is $M$. 
Since a surgery diagram can be thought of as a Kirby diagram of a 4-manifold with $M$ as its boundary, we have two Kirby diagrams $\mathcal{K}_0$ and $\mathcal{K}$ of some 4-manifold with the circle bundle over a torus as boundary given in the Figure \ref{fig:KirbyDiagramFor3-torus}. 
%These two Kirby diagrams $\mathcal{K}$ and $ \mathcal{K}_0 $ correspond to 4-manifolds whose boundary is the 3-torus.
If we want to convert $\mathcal{K}_0$ %(see Figure \ref{fig:BR3-torus}) 
to the other Kirby diagram $\mathcal{K}$, then we need to insert a 1-handle.% in the Boromean rings. 
The only way to do it is by the addition of a cancelling pair of handles of indices 1 and 2. The newly created cancelling pair of handles can not be made into a non-cancelling pair of handles by any Kirby move. 
In other words, after performing any Kirby move, a cancelling pair continues to be a cancelling pair.    
Hence, a 1-handle can not be added to $\mathcal{K}_0$ to obtain the Kirby diagram $\mathcal{K}$. 
\iffalse\begin{figure}[ht]
    \centering
    \includegraphics[scale=0.4]{Figures/BoromeanRing.pdf}
    \caption{Surgery diagram of 3-torus}
    \label{fig:BR3-torus}
\end{figure}\fi
Thus, we can achieve the Kirby calculus only for the round surgery diagrams consisting of joint pairs. 
Therefore, the subset of the round surgery diagram of joint pairs is the optimal choice to have the bridge theorem.
\end{rem}

\section{Equivalence Moves on round surgery diagrams of joint pairs}
\label{Sec : KCforRSDJP}

In this section, we define equivalence moves on round surgery diagrams of joint pairs. 
We prove that these moves are sufficient to establish a relation between two round surgery diagrams of joint pairs corresponding to the same 3-manifold. 

The equivalence moves are defined as follows.

\begin{enumerate}
\item {\bf Equivalence Move 1.}

The round surgery coefficients of a given joint pair $L_{11} \cup L_{12}$ with round 1-surgery coefficient $n_{11}$ and $n_{12}$, and round 2-surgery coefficient $m_1$ can be changed to new round surgery coefficient as shown below in Figure \ref{fig: EQM1} for any integer $k\in \mathbb{Z}$. 
In particular, if we take $k= n_{12}$, we get the original round surgery coefficients. 
\begin{figure}[ht]
    \centering
    \includegraphics[scale=0.3]{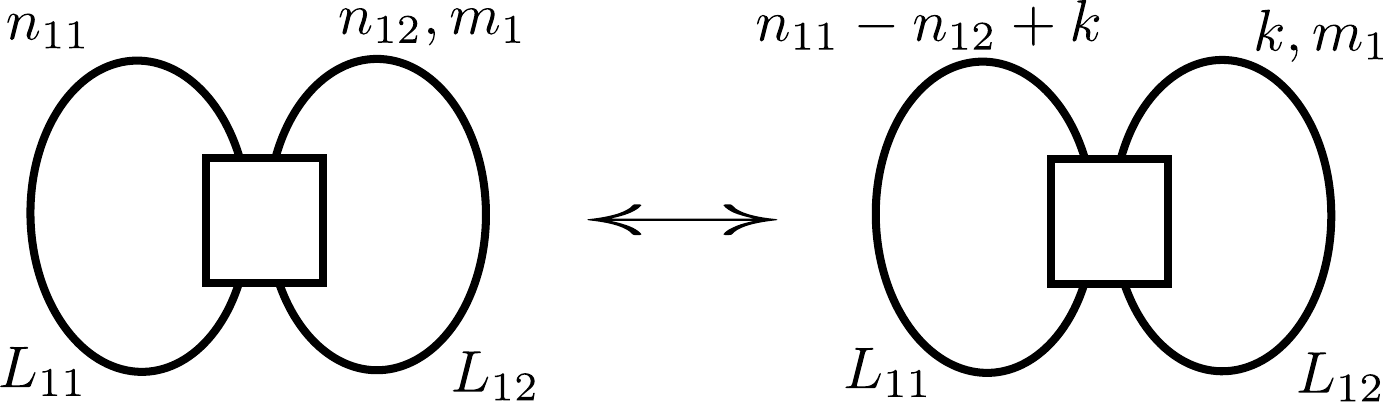}
    \caption{Equivalence move 1}
    \label{fig: EQM1}
\end{figure}

\item {\bf Equivalence Move 2 (Shuffle moves).} 
%Suppose we have given an ordinary surgery diagram consisting of a two-component link. 
%We apply Lemma \ref{lem: OSDtoJPRS}, and there is an ambiguity while assigning the round 2-suregry knot to a component. 
%As a result, we may get two different joint pairs, but both of them correspond to the same 3-manifold because their Dehn surgery diagram is the same (from Lemma \ref{lem: JPtoOSD}).
%Thus, we get some kind of equivalence between the round surgery diagrams of joint pairs. 
We define this equivalence move to address the ambiguity pointed out in Remark \ref{rem: AmbiguityInOSDtoJP}.

Given a round surgery diagram of joint pairs $L= L_1 \cup \ldots \cup L_{n}$ such that for each $1\leq i\leq n, \, L_{i}= L_{i1}\cup L_{i2}$ is a joint pair of round surgeries of index 1 and 2. 
In order to overcome the ambiguity observed in Remark \ref{rem: AmbiguityInOSDtoJP}, we define two types of shuffling on the round 2-surgery knot, namely {\it shuffle moves of type A and B}.

\paragraph{\bf Shuffle move of type A}
Suppose $L_{11} \cup L_{12}$ is a joint pair with round 1-surgery coefficients $n_{11}$ on $L_{11}$ and $n_{12}$ on $L_{12}$, and round 2-surgery coefficient $m_1$ on $L_{12}$.
The shuffle move of type A changes the round 2-surgery knot from $L_{12}$ to $L_{11}$ as shown in Figure \ref{fig:SMTA}. 
In particular, we have converted $L_{11}\cup L_{12} \to L^{\prime}_{11}\cup L^\prime_{12}$ where $L^{\prime}_{11}= L_{12}$ and $L^{\prime}_{12}= L_{11}$. 
Further, $L^{\prime}$ has a new round surgery coefficient. 
The round 1-surgery coefficients on $L^{\prime}_{11}$ is $k-n_{11}+n_{12}$ and $L^{\prime}_{12}$ is $k$, and round 2-surgery coefficient on $L^{\prime}_{12}$ is $n_{11}+m_1-n_{12}$ for any $k\in \mathbb{Z}$.  

\begin{figure}[ht]
        \centering
        \includegraphics[scale=0.3]{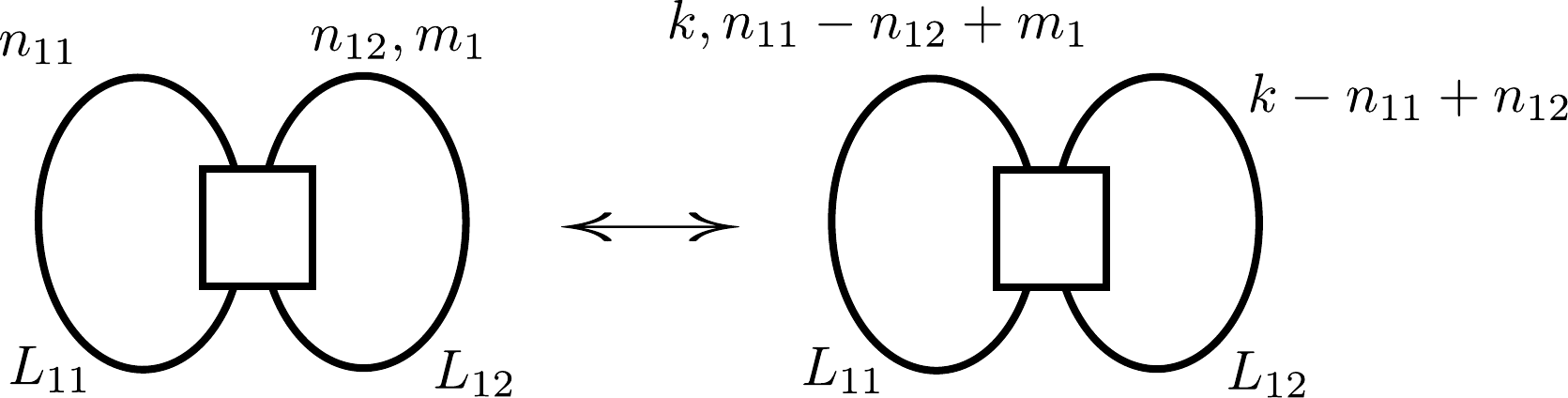}
        \caption{Shuffle move of type A}
        \label{fig:SMTA}
    \end{figure}

\paragraph{\bf Shuffle move of type B}
Suppose we have given a round surgery diagram consisting of two joint pairs $L= \bigcup_{i=1}^{2} L_{i1}\cup L_{i2}$ such that $L_{ij}$ has round 1-surgery coefficients $n_{ij}$ and $L_{i2}$ has round 2-surgery coefficient $m_i$ for $1\leq i,j\leq 2$. 
 \begin{figure}[ht]
        \centering
        \includegraphics[scale=0.3]{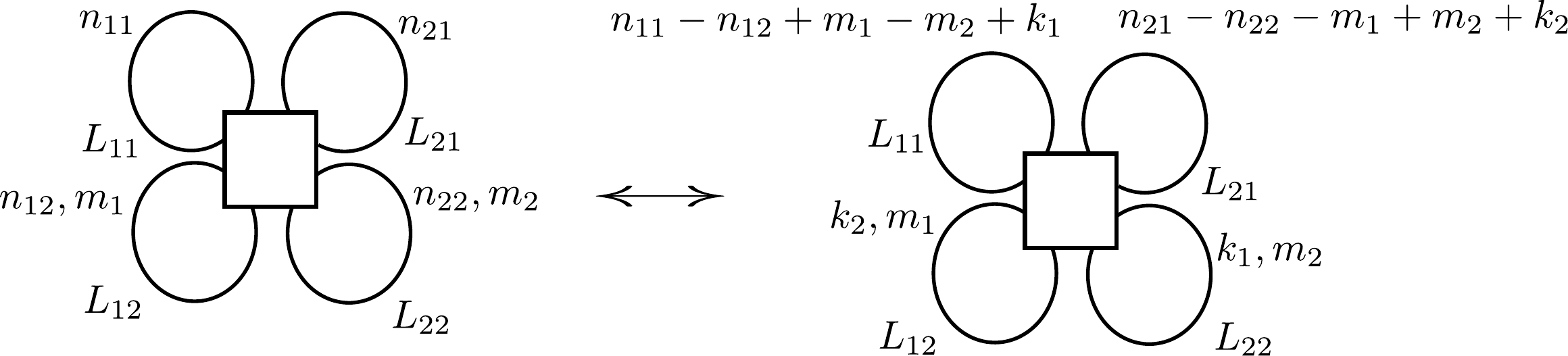}
        \caption{Shuffle move of type `B'}
        \label{fig:SMTB}
    \end{figure}
The shuffle move of type B interchanges round 2-surgery knots $L_{12}$ with $L_{22}$ as shown in Figure \ref{fig:SMTB}.
In particular, we have converted $L \to L^{\prime} (:= \bigcup_{i=1}^{2} L^{\prime}_{i1}\cup L^\prime_{i2})$ where $L^{\prime}_{12}= L_{22}$ and $L^\prime_{22}= L_{12}$ such that the round 1-surgery coefficients on $L^{\prime}_{11}$ is $n_{11}-n_{12}+m_1-m_2+k_1$ and $L^{\prime}_{12}$ is $k_2$, and round 2-surgery coefficient on $L^{\prime}_{12}$ is $m_1$; and  $L^{\prime}_{21}$ is $n_{21}-n_{22}+m_2-m_1+k_2$ and $L^{\prime}_{22}$ is $k_1$, and round 2-surgery coefficient on $L^{\prime}_{22}$ is $m_2$ for any $k_1, k_2\in \mathbb{Z}$.

%Any shuffling of a round 2-surgery knot in a round surgery diagram of joint pairs can be realised as a finitely many times application of the shuffle moves of type A and B on it. 
%We state it as a lemma below. 

\item {\bf Equivalence Move 3.} 

Addition or deletion of joint pair $U=U_1 \cup U_2$, consisting of an unlinked pair of unknots with round 1-surgery coefficient $k_1$ on $U_1$ and $k \in \mathbb{Z}$ on $U_2$ where $k_1 = k $ or $k\pm2$, and round 2-surgery coefficient $\pm1$ on $U_2$, to the given round surgery diagram $R$ as shown in the Figure \ref{fig:FRKM}.
\begin{figure}[ht]
        \centering
        \includegraphics[scale=0.4]{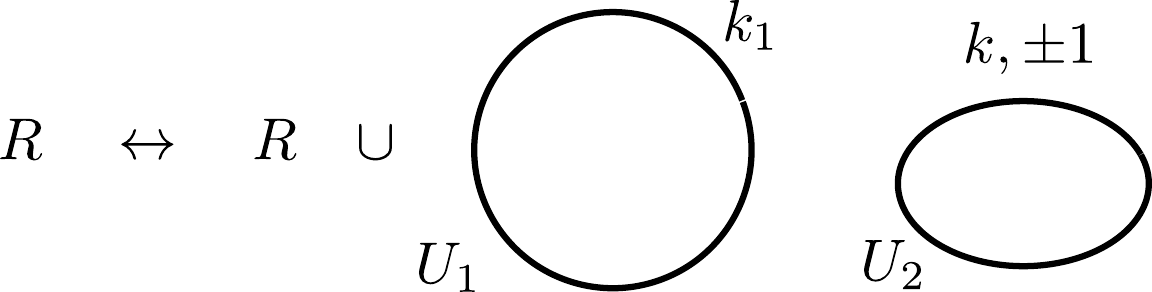}
        \caption{Equivalence move of type 3.}
        \label{fig:FRKM}
\end{figure}
\iffalse

\begin{figure}[h]
        \centering
        \includegraphics[scale=0.3]{RKM Images/FRKMAddition.pdf}
        \caption{First round Kirby move as addition}
        \label{fig:FRKMA}
    \end{figure}

    \begin{figure}[h]
        \centering
        \includegraphics[scale=0.3]{RKM Images/RFKMDeletion1.pdf}
        \caption{No deletion in single unknot component in First round Kirby move}
        \label{fig:FRKMD1}
    \end{figure}

    \begin{figure}[h]
        \centering
        \includegraphics[scale=0.3]{RKM Images/RFKMDeletion2.pdf}
        \caption{A pair of unknots deleted from the round surgery diagram}
        \label{fig:FRKMD2}
    \end{figure}
\fi 
\item {\bf Equivalence Move 4.} 
Suppose we have given a round surgery diagram $R$ consisting of at least one joint pair $L_{11}\cup L_{12}$. 
We define the following two moves on a single joint pair. We denote the linking number between $L_{ij}$ and $L_{i'j'}$ by $l_{iji'j'}$, and $L_{ij}(n)$ denotes an isotopic curve on the boundary of the tubular neighbourhood $N(L_{ij})$ that goes $n$ times around the oriented core curve $L_{ij}$.  

In Figure \ref{fig: EM4-1112}, we describe the first type of equivalence move 4 on a joint pair $L_{11}\cup L_{12}$. 
This equivalence move 4 on a joint pair $L_{11}\cup L_{12}$  replaces the $L_{11}$ by the band connected sum $L^{\prime}_{11} =  L_{11}\#_b L_{12}(m_1)$. 
Under this move, the new component $L^{\prime}_{11}$ gets a new round 1-surgery coefficient $n_{11}^{\prime}= n_{11}-n_{12}+m_1+ 2l_{1112}+k$ and the other coefficients remain the same.

\begin{figure}[ht]
    \centering
    \includegraphics[scale=0.3]{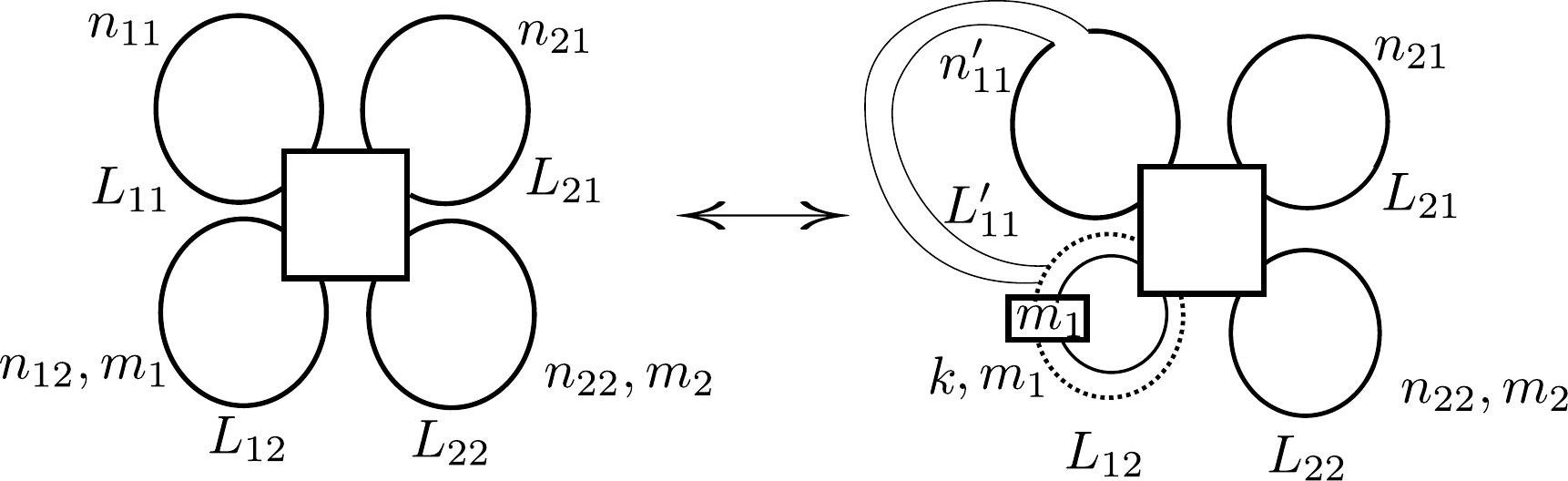}
    \caption{This equivalence move 4 on a Joint pair $L_{11}\cup L_{12}$  replaces the $L_{11}$ by the band connected sum $L_{11}\#_b L_{12}(m_1)$ and changes the round surgery coefficients to $n_{11}^{\prime}$ on $L'_{11}$ and $n_{12}, m_1$ on $L_{12}$.}
    \label{fig: EM4-1112}
\end{figure}

In Figure \ref{fig: EM4-1211}, we describe the second type of equivalence move 4 on a joint pair $L_{11}\cup L_{12}$.
This equivalence move 4 on a joint pair $L_{11}\cup L_{12}$ replaces the $L_{12}$ by the band connected sum $L_{12}\#_b L_{11}(f)$, where $f=n_{11}-n_{12}+m_1$. 
Under this move, the new joint pair $L_{11}\cup L^{\prime}_{12}$ gets a new round 1-surgery coefficient $n^{\prime}_{11}=-m_1-2l_{1112}+k$ on $L_{11}$ and $k$ on $L^{\prime}_{12}$, for $k\in \mathbb{Z}$ and $m_1=2m_1+n_{11}-n_{12}+2l_{1112}$ on $L^{\prime}_{12}$. 
Further, the remaining joint pairs and their round surgery coefficients remain unchanged under this move. 

\begin{figure}[ht]
    \centering
    \includegraphics[scale=0.3]{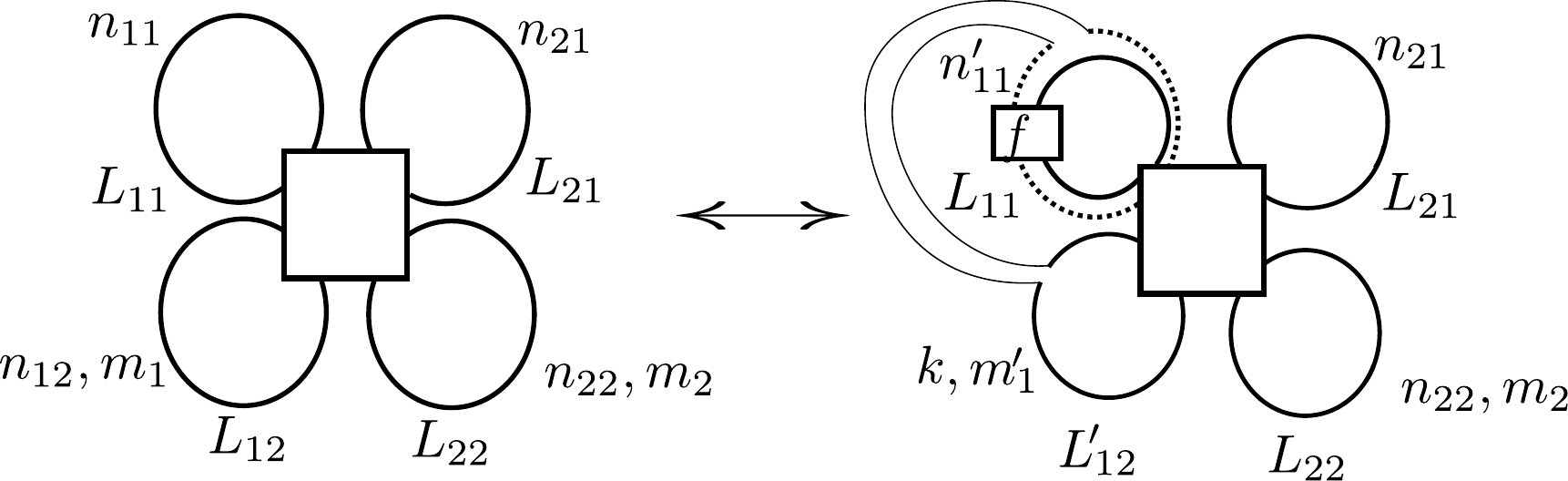}
    \caption{This equivalence move 4 on a Joint pair $L_{11}\cup L_{12}$  replaces the $L_{12}$ by $L_{12}\#_b L_{11}(f)$ and changes the round 1-surgery coefficients to $n_{11}^{\prime}$ on $L_{11}$ and $k$ on $L^{\prime}_{12}$, for $k\in \mathbb{Z}$ and round 2-surgery coefficient on $L^{\prime}_{12}$ to $m^{\prime}_1$.}
    \label{fig: EM4-1211}
\end{figure}

Now we assume that round surgery diagram $R$ has at least two joint pairs $(L_{11}\cup L_{12}) \cup (L_{21} \cup L_{22})$. 
We define four more types of equivalence move 4. 

In Figure \ref{fig: EM4-1121}, we describe the first type of equivalence move 4 on two joint pairs $(L_{11}\cup L_{12}) \cup (L_{21}\cup L_{22})$.
This equivalence move 4 on $(L_{11}\cup L_{12}) \cup (L_{21}\cup L_{22})$ replaces the $L_{11}$ by the band connected sum $L^{\prime}_{11} = L_{11} \#_b L_{21}(f)$, where $f=n_{21}-n_{22}+m_2$. 
Under this move, the new joint pair $L^{\prime}_{11}\cup L_{12}$ gets a new round 1-surgery coefficient $n^{\prime}_{11}= n_{11}+n_{21}- (n_{12}+n_{22})+m_2+ 2 l_{1121}+k$ on $L^{\prime}_{11}$ and $k$ on $L_{12}$, for $k\in \mathbb{Z}$ and $m_1$ on $L_{12}$. 
Further, the remaining joint pairs and their round surgery coefficients remain unchanged under this move. 

\begin{figure}[ht]
    \centering
    \includegraphics[scale=0.3]{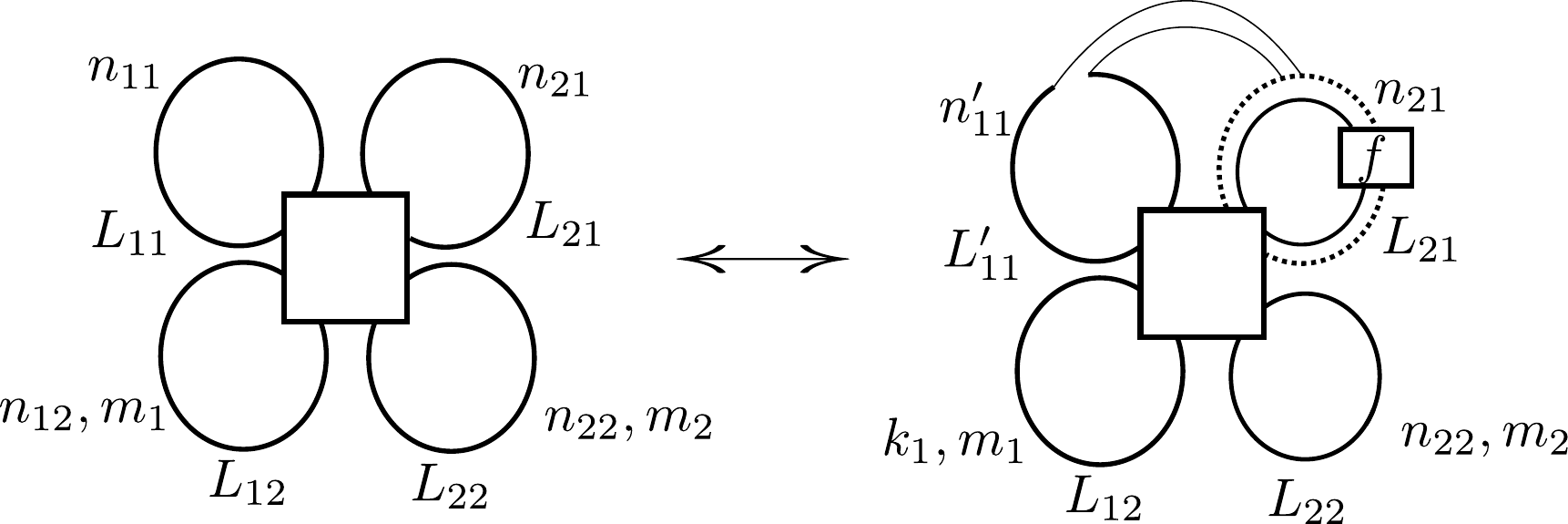}
    \caption{This equivalence move 4 on $(L_{11}\cup L_{12}) \cup (L_{21}\cup L_{22})$ replaces the $L_{11}$ by $L^{\prime}_{11} = L_{11} \#_b L_{21}(f)$ and changes the round 1-surgery coefficients to $n_{11}^{\prime}$ on $L^{\prime}_{11}$ and $k$ on $L_{12}$, for $k\in \mathbb{Z}$.}
    \label{fig: EM4-1121}
\end{figure}

In Figure \ref{fig: EM4-1122}, we describe the second type of equivalence move 4 on two joint pairs $(L_{11}\cup L_{12}) \cup (L_{21}\cup L_{22})$.
This equivalence move 4 on $(L_{11}\cup L_{12}) \cup (L_{21}\cup L_{22})$ replaces the $L_{11}$ by the band connected sum $L^{\prime}_{11} = L_{11} \#_b L_{22}(m_2)$. 
Under this move, the new joint pair $L^{\prime}_{11}\cup L_{12}$ gets a new round 1-surgery coefficient $n^{\prime}_{11}= n_{11}-n_{12}+m_2+2l_{1122}+k$ on $L^{\prime}_{11}$ and $k$ on $L_{12}$, for $k\in \mathbb{Z}$ and round 2-surgry coefficient $m_1$ on $L^{\prime}_{12}$. 
Further, the remaining joint pairs and their round surgery coefficients remain unchanged under this move.

\begin{figure}[ht]
    \centering
    \includegraphics[scale=0.3]{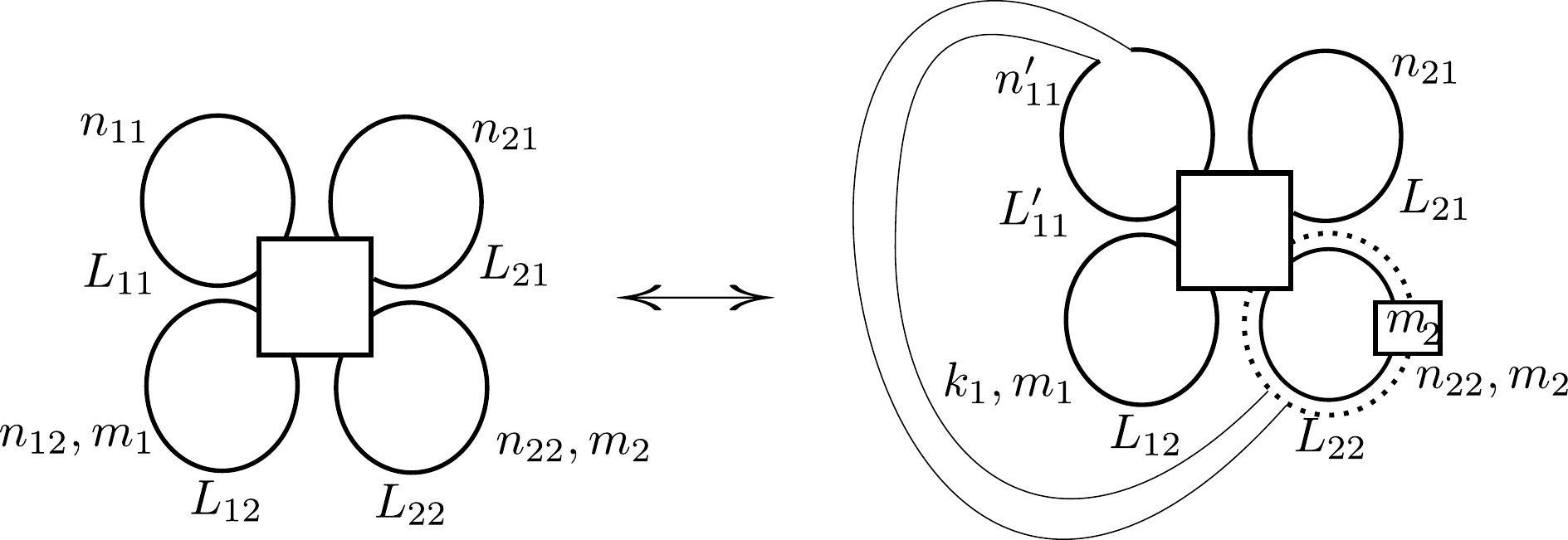}
    \caption{This equivalence move 4 on $(L_{11}\cup L_{12}) \cup (L_{21}\cup L_{22})$ replaces the $L_{11}$ by $L^{\prime}_{11} = L_{11} \#_b L_{22}(m_2)$ and changes the round 1-surgery coefficients to $n_{11}^{\prime}$ on $L^{\prime}_{11}$ and $k$ on $L_{12}$, for $k\in \mathbb{Z}$.}
    \label{fig: EM4-1122}
\end{figure}

In Figure \ref{fig: EM4-1221}, we describe the third type of equivalence move 4 on two joint pairs $(L_{11}\cup L_{12}) \cup (L_{21}\cup L_{22})$.
This equivalence move 4 on $(L_{11}\cup L_{12}) \cup (L_{21}\cup L_{22})$ replaces the $L_{12}$ by the band connected sum $L^{\prime}_{12} = L_{12} \#_b L_{21}(f)$, where $f= n_{21}-n_{22}+m_2$. 
Under this move, the new joint pair $L_{11}\cup L^{\prime}_{12}$ gets a new round 1-surgery coefficient $n^{\prime}_{11}=(n_{11}-n_{12})-(n_{21}-n_{22}) -m_2 - 2l_{1221}+k$ on $L_{11}$ and $k$ on $L^{\prime}_{12}$, for $k\in \mathbb{Z}$ and a new round 2-surgery coefficient $m^{\prime}_{1}= m_1+m_2+n_{21}-n_{22}+2l_{1221}$ on $L^{\prime}_{12}$. 
Further, the remaining joint pairs and their round surgery coefficients remain unchanged under this move.

\begin{figure}[ht]
    \centering
    \includegraphics[scale=0.3]{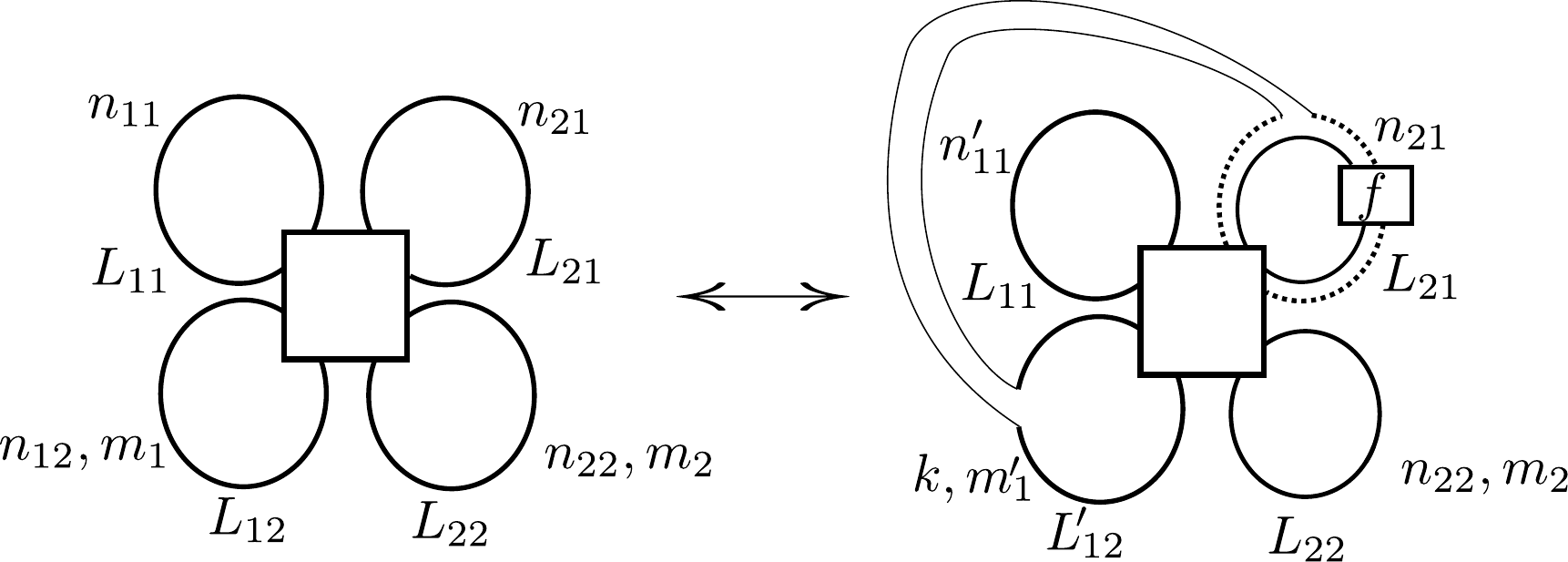}
    \caption{This equivalence move 4 on $(L_{11}\cup L_{12}) \cup (L_{21}\cup L_{22})$ replaces the $L_{12}$ by $L^{\prime}_{12} = L_{12} \#_b L_{21}(f)$ and changes the round 1-surgery coefficients to $n_{11}^{\prime}$ on $L^{\prime}_{11}$ and $k$ on $L_{12}$, for $k\in \mathbb{Z}$ and round 2-surgery coefficient to $m^{\prime}_1$.}
    \label{fig: EM4-1221}
\end{figure}

In Figure \ref{fig: EM4-1222}, we describe the fourth type of equivalence move 4 on two joint pairs $(L_{11}\cup L_{12}) \cup (L_{21}\cup L_{22})$.
This equivalence move 4 on $(L_{11}\cup L_{12}) \cup (L_{21}\cup L_{22}$ replaces the $L_{12}$ by the band connected sum $L^{\prime}_{12} = L_{12} \#_b L_{22}(m_2)$. 
Under this move, the new joint pair $L_{11}\cup L^{\prime}_{12}$ gets a new round 1-surgery coefficient $n^{\prime}_{11}= n_{11}-n_{12}-m_2- 2 l_{1222} +k$ on $L_{11}$ and $k$ on $L^{\prime}_{12}$, for $k\in \mathbb{Z}$ and a new round 2-surgery coefficient $m_{1}^{\prime}= m_1+m_2 +2l_{1222}$ on $L^{\prime}_{12}$. 
Further, the remaining joint pairs and their round surgery coefficients remain unchanged under this move.

\begin{figure}[ht]
    \centering
    \includegraphics[scale=0.3]{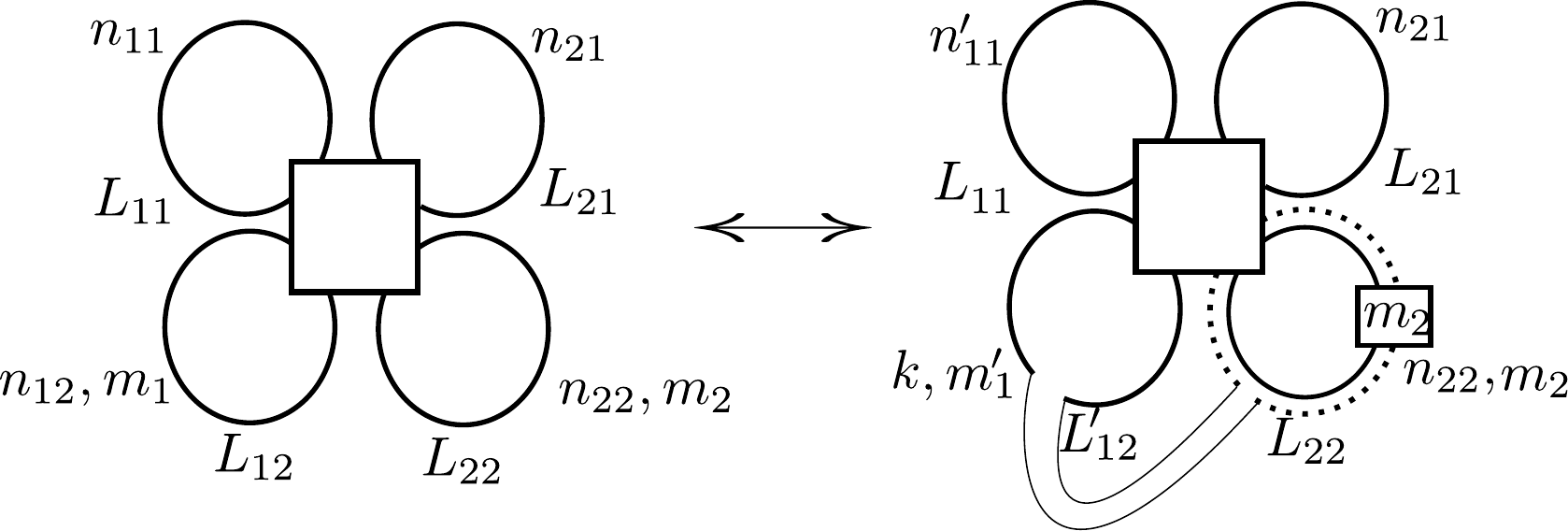}
    \caption{This equivalence move 4 on $(L_{11}\cup L_{12}) \cup (L_{21}\cup L_{22})$ replaces the $L_{12}$ by $L^{\prime}_{12} = L_{12} \#_b L_{22}(m_2)$ and changes the round 1-surgery coefficients to $n_{11}^{\prime}$ on $L^{\prime}_{11}$ and $k$ on $L_{12}$, for $k\in \mathbb{Z}$ and round 2-surgery coefficient to $m^{\prime}_1$.}
    \label{fig: EM4-1222}
\end{figure}
\end{enumerate}

\begin{rem}
    Equivalence move 3 can be thought of as a round version of the Kirby move of type 1, and equivalence move 4 can be realised as a round version of the Kirby move of type 2. 
    In other words, we may also call equivalence move 3 as {\it round Kirby move of type 1, or round first Kirby move} and equivalence move 4 as {\it round Kirby move of type 2, or round second Kirby move}. 
\end{rem}

\begin{lem}
    The application of equivalence moves 1-4 on a given round surgery diagram $L$ does not change the resulting 3-manifold. 
\end{lem}
\begin{proof}
Suppose $L= \bigcup_{i=1}^{n} (L_{i1}\cup L_{i2})$ is a round surgery diagram of joint pairs with round 1-surgery coefficient $n_{ij}$ on component $L_{ij}$ and round 2-surgery coefficient $m_i$ on $L_{i2}$ for $1\leq i\leq n$ and $1\leq j\leq 2$. 
Below, we see that the effect of each equivalence move on $L$ does not change the resultant 3-manifold. 

{\bf Equivalence Move 1.} We apply equivalence move 1 on a joint pair of the given round surgery diagram $L$. 
Without loss of generality, we apply the equivalence move 1 on the first joint pair $L_{11} \cup L_{12}$ and obtain new round surgery coefficients on $L_{11}\cup L_{12}$.
In particular, we get a new joint pair. 
Both of the joint pairs can be converted to the same Dehn surgery diagram $D$ by Lemma \ref{lem: JPtoOSD} as shown in Figure \ref{fig: EM1Proof}. 
\begin{figure}[ht]
        \centering
        \includegraphics[scale=0.3]{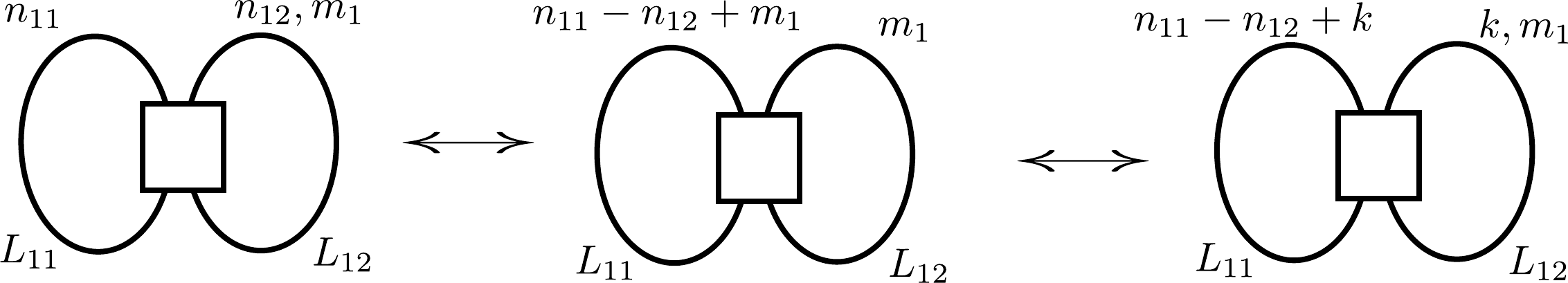}
        \caption{Diagrams involved in the equivalence move 1, i.e., diagrams on the left and right, correspond to the same Dehn surgery diagram in the middle.}
        \label{fig: EM1Proof}
\end{figure}

{\bf Equivalence Move 2.} The equivalence move 2 has two types of separate moves: The shuffle move of type A and type B. 
Since we apply a shuffle move of a certain type once at a time, we discuss their effects separately. 

{\it Shuffle Move of type A.}  Since the shuffle move of type A applies to a joint pair, without loss of generality, we apply the shuffle move A on the first joint pair $L_{11} \cup L_{12}$ of the given round surgery diagram $L$. 
\begin{figure}[ht]
        \centering
        \includegraphics[scale=0.3]{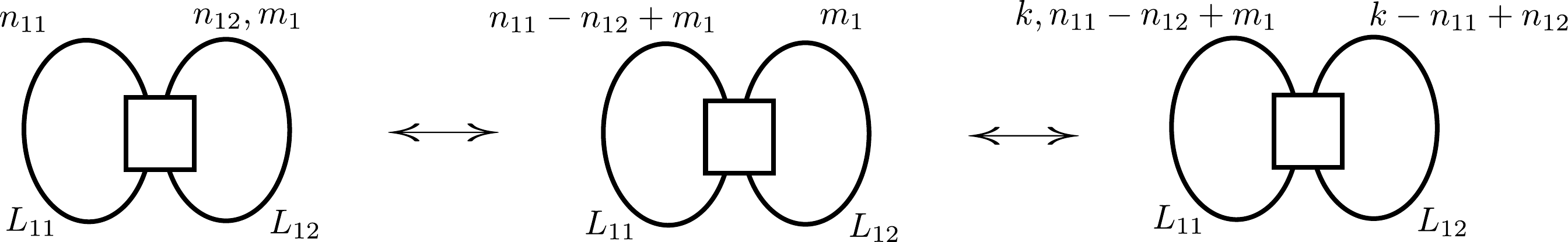}
        \caption{Diagrams involved in the shuffle move of type A, i.e., diagrams on the left and right, correspond to the same Dehn surgery diagram in the middle.}
        \label{fig:SMTAProof}
\end{figure}
Suppose $L^{\prime}_{11} \cup L_{12}^{\prime}$ is the new joint pair obtained after performing a shuffle move of type A on $L_{11} \cup L_{12}$. 
By Theorem \ref{BridgeTheorem}, the Dehn surgery diagram corresponding to both joint pairs is the same (see Figure~\ref{fig:SMTAProof}).
Thus, the resultant 3-manifold remains unchanged.

{\it Shuffle move of type B.} The shuffle move of type B applies to two joint pairs. In the given round surgery diagram $L$, we assume that we apply the type B move on the first two joint pairs $(L_{11} \cup L_{12}) \cup (L_{21} \cup L_{22})=: L_1 \cup L_2$.
After performing a shuffle move of type B on $L_1 \cup L_2$, we get a new round surgery diagram $L^{\prime}= \bigcup_{i=1}^{2} (L^{\prime}_{i1}\cup L^{\prime}_{12}) \bigcup_{i=3}^{n} (L_{i1}\cup L_{i2})$. We denote the new joint pairs by $L^{\prime}_i:= (L^{\prime}_{i1}\cup L^{\prime}_{12})$ for each $i=1,2$. 
\begin{figure}[ht]
    \includegraphics[scale=0.25]{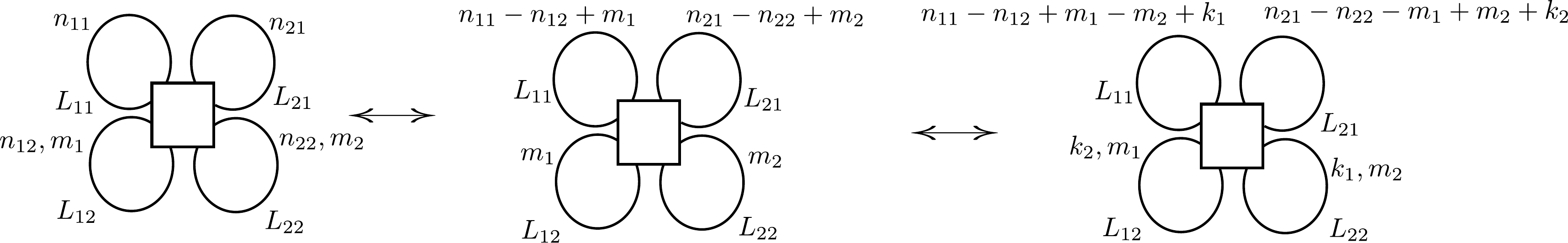}
    \caption{The shuffle move of type B applies to the diagram on the left to obtain the diagram on the right, and both of the diagrams correspond to the same surgery diagram in the middle.}
    \label{fig:SMTBProof}
\end{figure}
By Theorem \ref{BridgeTheorem}, we may convert both round surgery diagrams $L$ and $L^{\prime}$ to the same Dehn surgery diagram as shown in Figure \ref{fig:SMTBProof}. 
Thus, the shuffle move of type B does not change the resultant 3-manifold.

{\bf Equivalence Move 3.}
This move corresponds to the addition or deletion of a joint pair $U= U_1 \cup U_2$ with round 1-surgery coefficient $k, k\pm2$ on  $U_1$ and $k$ on $U_2$, for some $k\in \mathbb{Z}$, and round 2-surgery coefficient $\pm1$ on $U_2$ to the given round surgery diagram $L$.
In particular, we get a new round surgery diagram $L^{\prime}= L \bigsqcup U$. 
By Lemma \ref{lem: OSDtoJPRS}, the joint pair $U$ corresponds to a pair of unknots $U_1 \cup U_2$ with each having Dehn surgery coefficients $n \in \{\pm1\}$. 
Recall that any unlinked disjoint $\pm1$ framed unknot corresponds to a connected sum of the resultant manifold with $\s^3$. In particular, the resultant 3-manifold is unchanged after an addition or deletion of $U$. 

{\bf Equivalence Move 4.}
We have given six types of equivalence moves for a round surgery diagram of joint pairs. 
\begin{figure}[ht]
    \centering
    \includegraphics[scale=0.3]{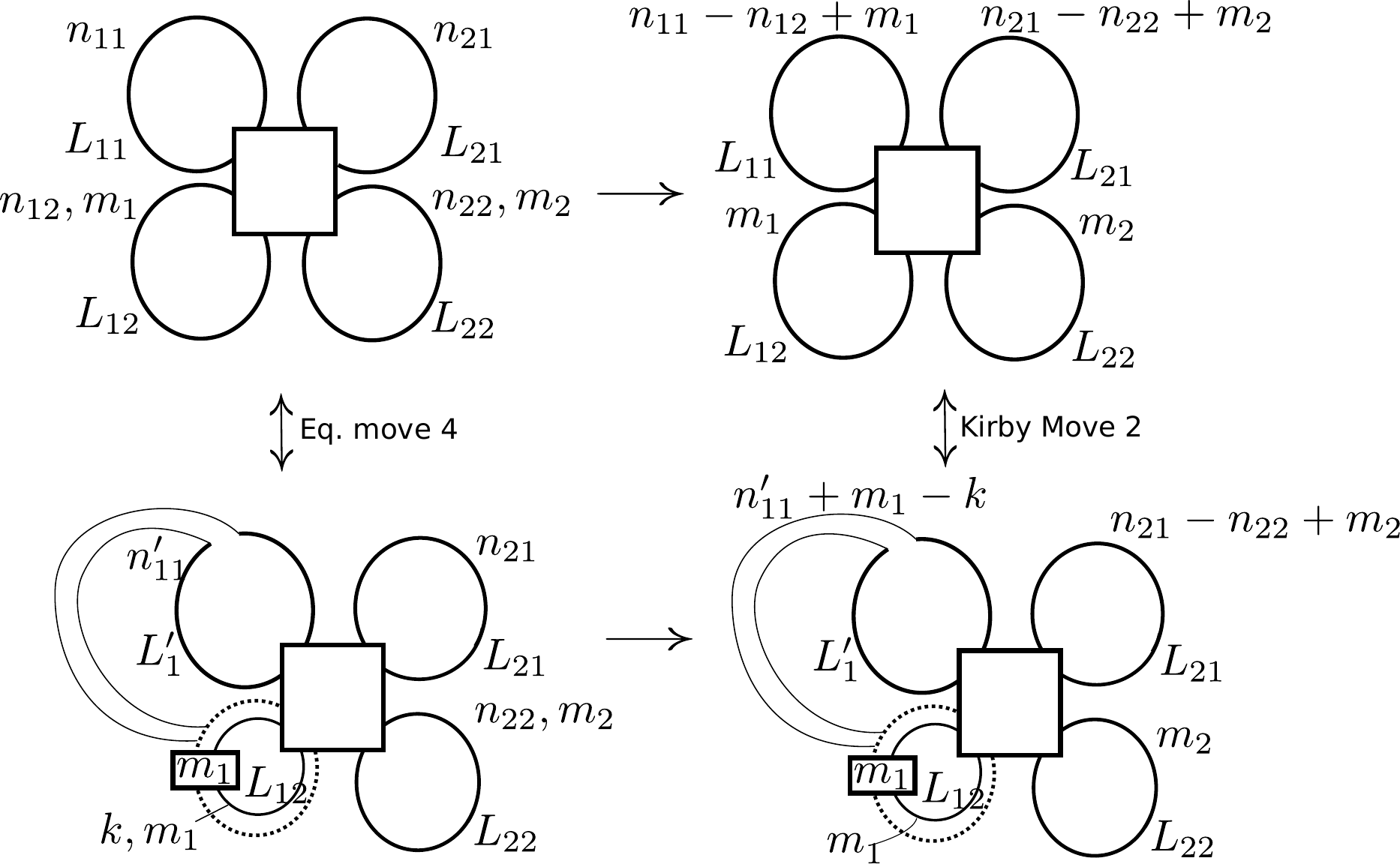}
    \caption{On the left column, the round surgery diagrams are in equivalence move 4, where $L_{11}$ is replaced with $L_{11}\#_b L_{12}(m_1)$, and on the right, their corresponding surgery diagrams in Kirby move 2.}
    \label{fig:EM4-1112-Proof}
\end{figure}

\begin{figure}[ht]
    \centering
    \includegraphics[scale=0.3]{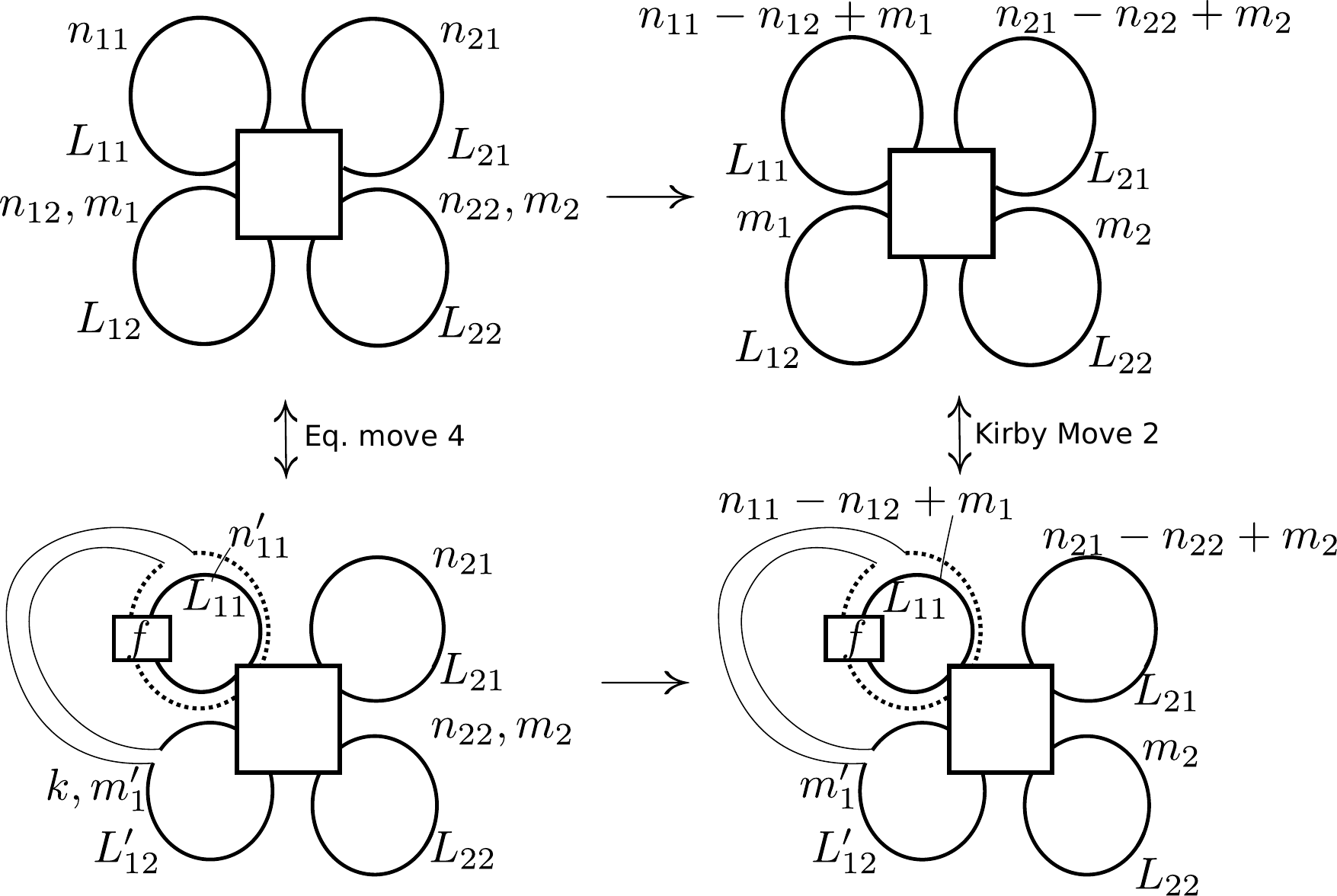}
    \caption{On the left column, the round surgery diagrams are in equivalence move 4, where $L_{12}$ replaces with $L_{12}\#_b L_{11}(f)$, and on the right, their corresponding surgery diagrams in Kirby move 2, where framing $f= n_{11}-n_{12}+m_1$.}
    \label{fig:EM4-1211-Proof}
\end{figure}
Each of the moves does not change the resultant 3-manifold. 
For each type of move, we convert the left and right sides of the equivalence move 4 diagrams to the Dehn surgery diagrams. 
We realise that the effect of the equivalence move 4 on the round surgery diagram is equivalent to the effect of Kirby move 2 on the corresponding Dehn surgery diagrams. 
\begin{figure}[ht]
    \centering
    \includegraphics[scale=0.3]{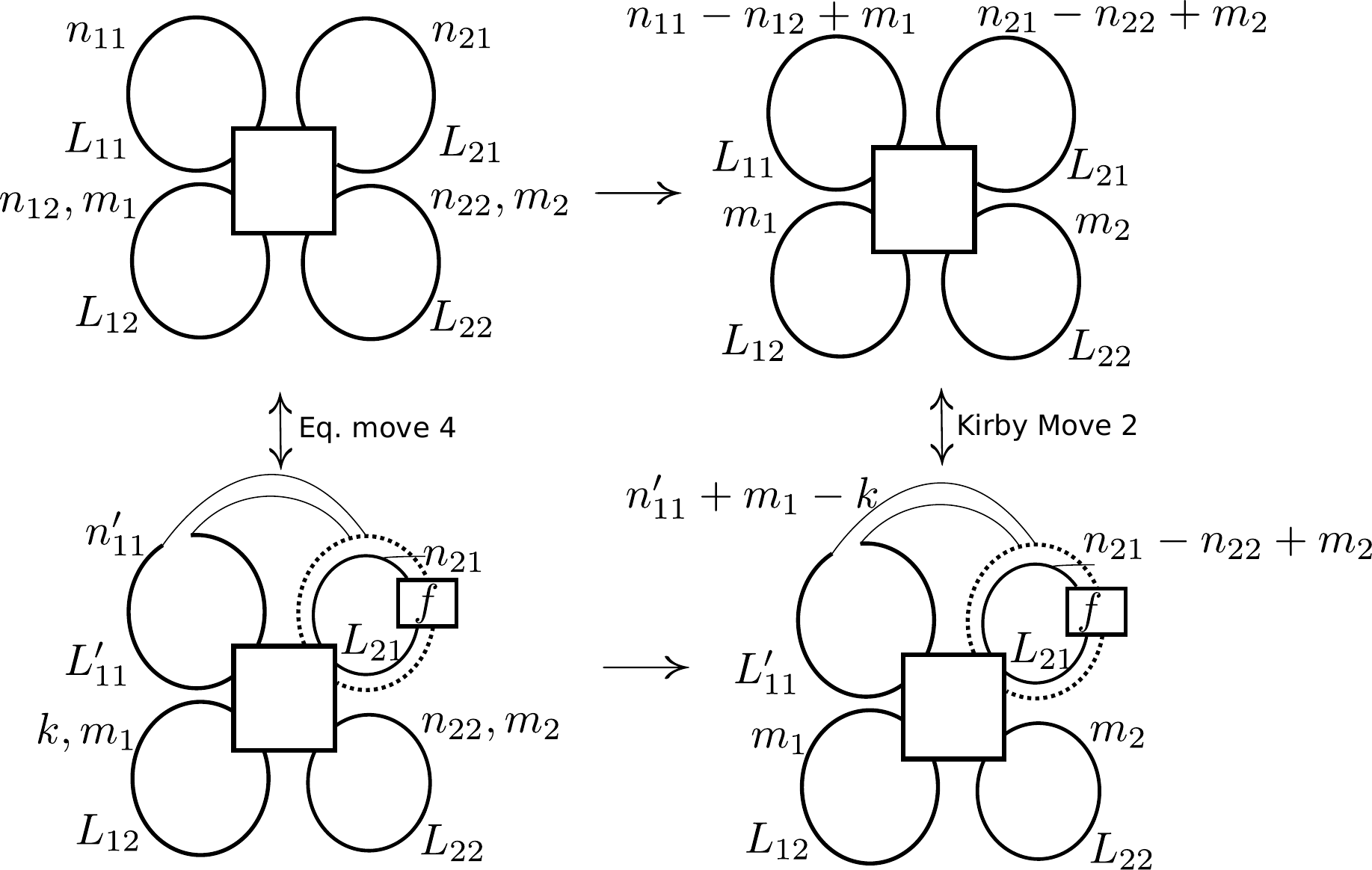}
    \caption{On the left column, the round surgery diagrams are in equivalence move 4, where $L_{11}$ replaces with $L_{11}\#_b L_{21}(f)$, and on the right, their corresponding surgery diagrams in Kirby move 2, where $f= n_{21}-n_{22}+m_2$}
    \label{fig:EM4-1121-Proof}
\end{figure}

\begin{figure}[ht]
    \centering
    \includegraphics[scale=0.3]{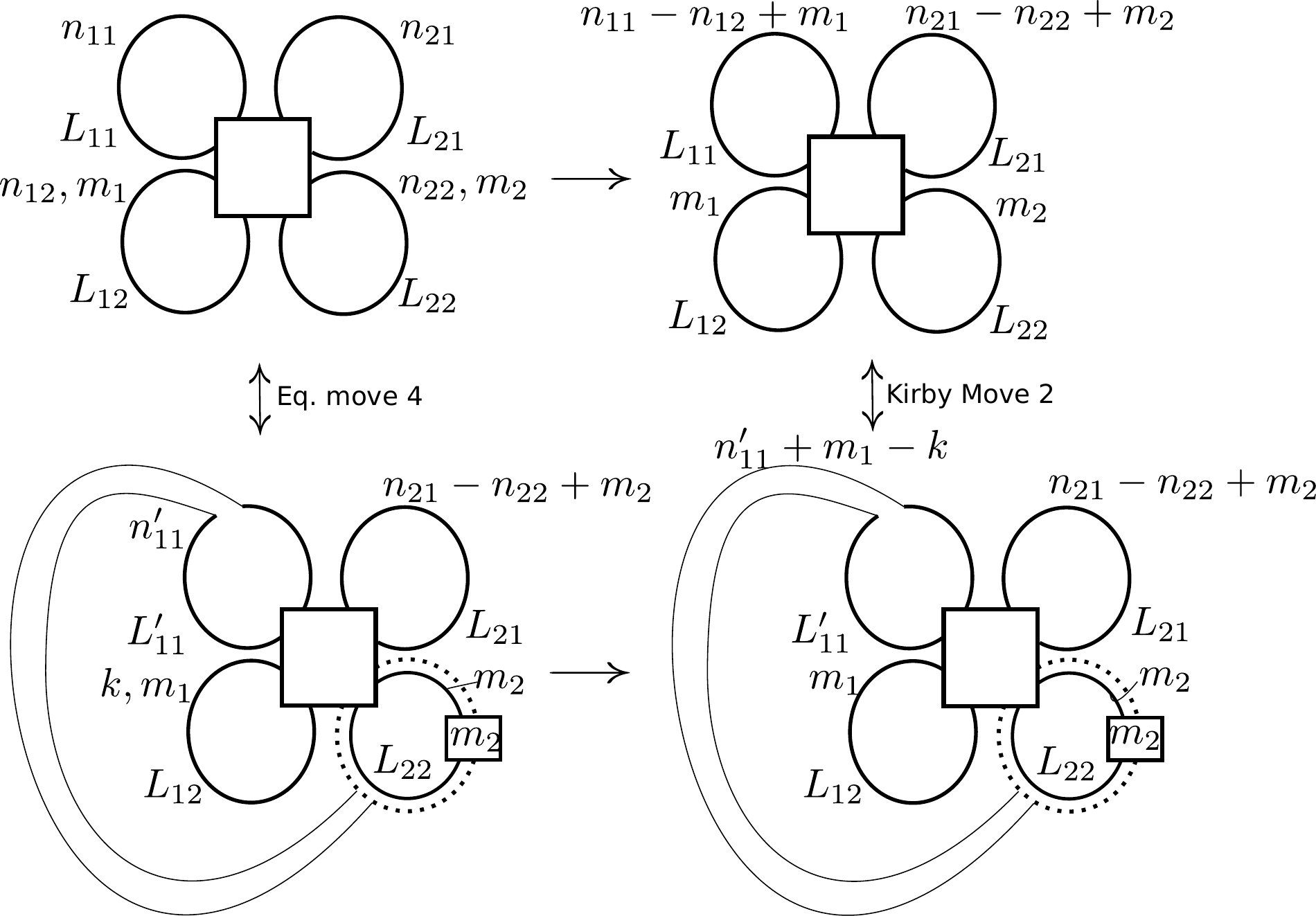}
    \caption{On the left column, the round surgery diagrams are in equivalence move 4, where $L_{11}$ is replaced with $L_{11}\#_b L_{22}(m_2)$, and on the right, their corresponding surgery diagrams in Kirby move 2.}
    \label{fig:EM4-1122-Proof}
\end{figure}

In particular, in Figures \ref{fig:EM4-1112-Proof} and \ref{fig:EM4-1211-Proof}, we realise equivalence move 4 on a single joint pair as performing Kirby move 2 on the corresponding Dehn's surgery diagrams.
And in Figures \ref{fig:EM4-1121-Proof}, \ref{fig:EM4-1122-Proof}, \ref{fig:EM4-1221-Proof} and \ref{fig:EM4-1222-Proof}, we realise the equivalence move 4 defined on two joint pairs as performing Kirby move 2 on the corresponding Dehn's surgery diagrams. 

We know the Kirby move 2 does not change the resultant 3-manifold. Hence, the equivalence move 4 also does not change the resultant 3-manifold. 

\begin{figure}[ht]
    \centering
    \includegraphics[scale=0.3]{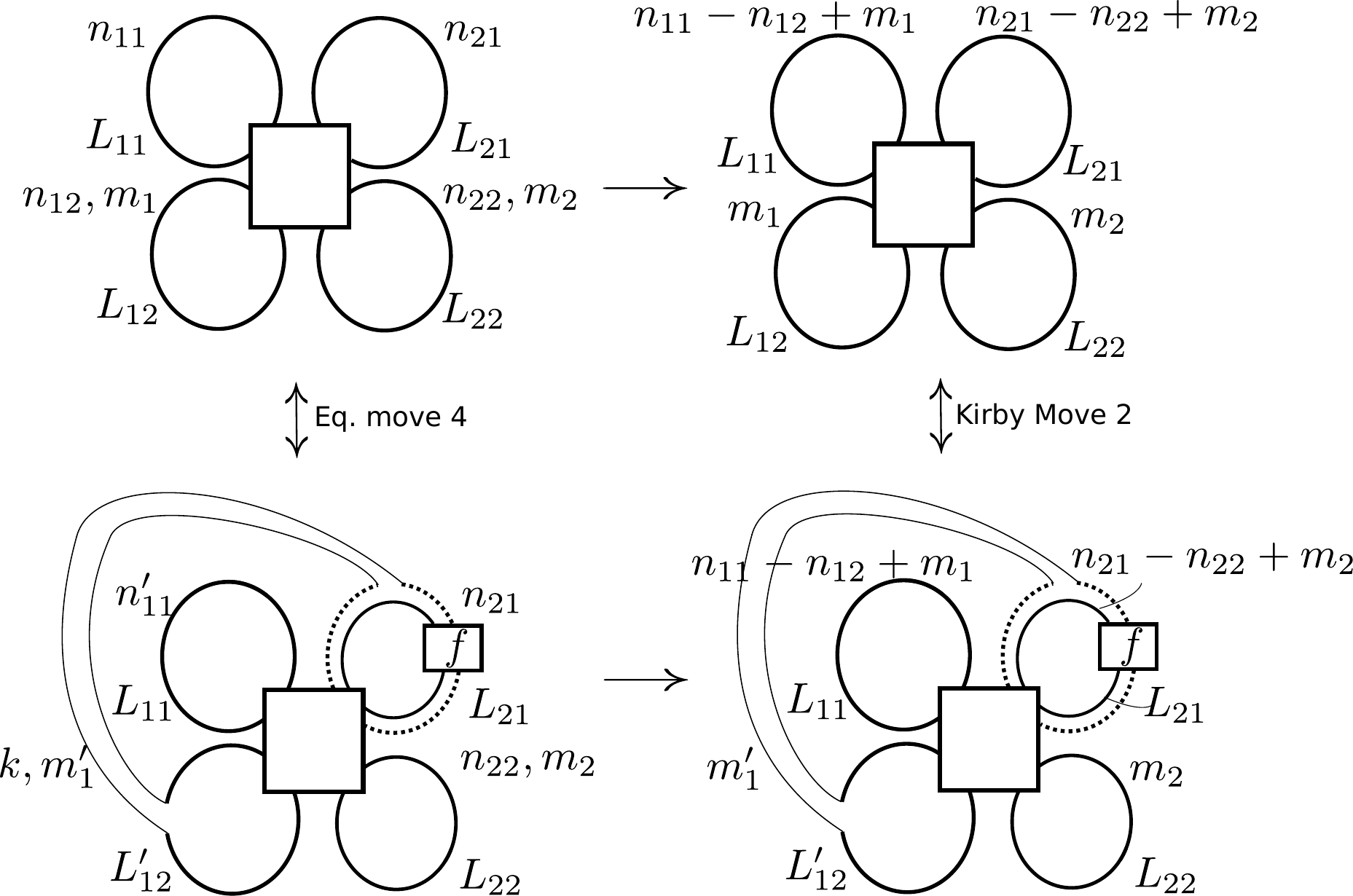}
    \caption{On the left column, the round surgery diagrams are in equivalence move 4, where $L_{12}$ replaces with $L_{12}\#_b L_{21}(f)$, and on the right, their corresponding surgery diagrams in Kirby move 2, where $f=n_{21}-n_{22}+m_2$.}
    \label{fig:EM4-1221-Proof}
\end{figure}

\begin{figure}[ht]
    \centering
    \includegraphics[scale=0.3]{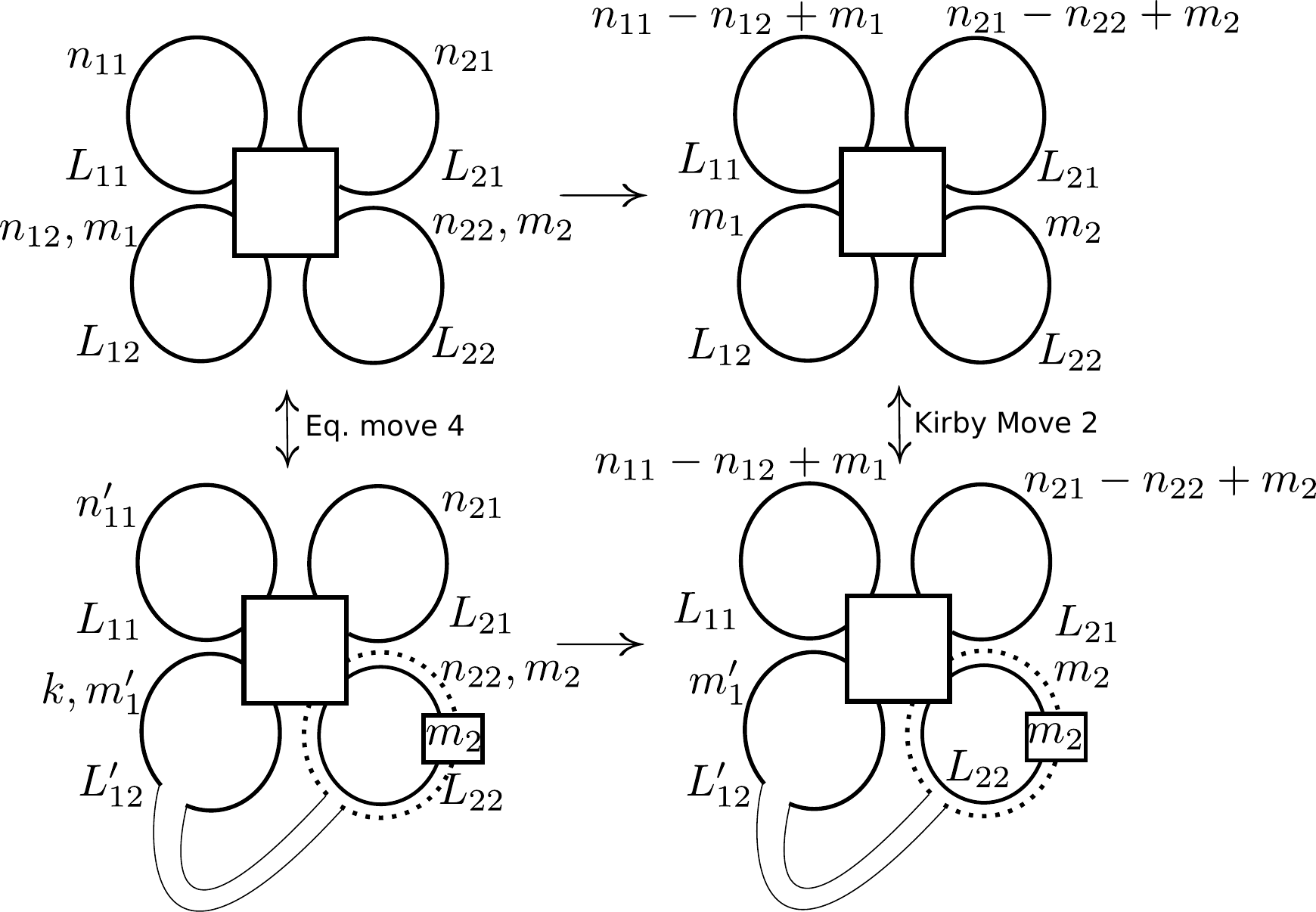}
    \caption{On the left column, the round surgery diagrams are in equivalence move 4, where $L_{12}$ is replaced with $L_{12}\#_b L_{22}(m_2)$, and on the right, their corresponding surgery diagrams in Kirby move 2.}
    \label{fig:EM4-1222-Proof}
\end{figure}
\end{proof}

\begin{lem}\label{lem: EM1n2onRSD}
Let $L=\bigcup_{i=1}^{n} (L_{i1} \cup L_{i2})$ and $L^{\prime}= \bigcup_{j=1}^{n} (L^{\prime}_{j1}\cup L^{\prime}_{j2})$ be two round surgery diagrams of joint pairs such that $L$ and $L^{\prime}$ has the same 
link components (up to knot isotopy) and linking between any two components of $L$ matches with the linking between corresponding two components of $L^{\prime}$. Assume that both of the round surgery diagrams correspond to the same Dehn surgery diagrams. Then, the round surgery diagrams $L$ and $L'$ are connected by a finite sequence of the shuffle moves of type A and type B and equivalence move 1.
\end{lem} 

\begin{proof}
We have given that $L = \bigcup_{i=1}^{n}(L_{i1}\cup L_{i2})$ is a round surgery diagram of joint pairs. Thus, by convention, we know that the link component $L_{i2}$ is the round 2-surgery knot for each $1\leq i\leq n$.
Since link components of $L$ and $L'$ are isotopic and linking between any two components of $L$ matches with the linking between corresponding two components of $L'$, we change the indexing of components of $L'$, whenever needed, such that $L^{\prime}_{ij}$ and $L_{ij}$ are knot isotopic, and $lk(L^{\prime}_{ij}, L^{\prime}_{ab}) = lk(L_{ij}, L_{ab})$ for all $1\leq a\leq n $ and $1\leq b\leq 2$.

For $1\leq a \leq n$, suppose $L^{\prime}_{a1}$ and $L^{\prime}_{a2}$ are the first pair that is not a joint pair as per the convention mentioned in Remark \ref{rem: StdPos}. 
In particular, we have the following cases. We want to convert the given pair into a joint pair as per convention, i.e., $L^{\prime}_{a2}$ is a round 2-surgery knot.
\begin{enumerate}
    \item Suppose $L^{\prime}_{k_01}\cup L^{\prime}_{k_02}$ are in joint pair but $L^{\prime}_{k_01}$ is the round 2-surgery knot.
    Then we apply the shuffle move of type A to get $L^{\prime}_{k_01}\cup L^{\prime}_{k_02}$ into a joint pair as per convention, i.e., $L^{\prime}_{k_02}$ is a round 2-surgery knot. 
    \item Suppose $L^{\prime}_{k_01}$ and $L^{\prime}_{k_02}$ are in joint pairs with different components. 
    In particular, we assume that $L^{\prime}_{k_01}$ and $L^{\prime}_{k_02}$ are in joint pair with $L^{\prime}_{i_0j_0}$ and $L^{\prime}_{i_1j_1}$, respectively, for some $1\leq i_0,i_1\leq n$ and $1\leq j_0, j_1 \leq 2$.
    We want to arrange the pair $L^{\prime}_{k_01}$ and $L^{\prime}_{k_02}$ into a joint pair as per our convention. For that, we consider the following cases. 
    In each case, we start with joint pairs as per our convention.
    \begin{enumerate}
        \item [Case 1.] The joint pairs are given as follows: $(L^{\prime}_{k_01} \cup L^{\prime}_{i_0j_0}) \cup (L^{\prime}_{k_02}\cup L^{\prime}_{i_1j_1})$.  
        We apply the type A move on the second pair to get 
        $$(L^{\prime}_{k_01} \cup L^{\prime}_{i_0j_0}) \cup ( L^{\prime}_{i_1j_1}\cup L^{\prime}_{k_02}) \xrightarrow[A]{} (L^{\prime}_{k_01} \cup L^{\prime}_{i_0j_0}) \cup ( L^{\prime}_{i_1j_1}\cup L^{\prime}_{k_02}). $$
        Then we apply the type B move on the latter to get
        $$ (L^{\prime}_{k_01} \cup L^{\prime}_{i_0j_0}) \cup ( L^{\prime}_{i_1j_1}\cup L^{\prime}_{k_02})\xrightarrow[]{B}(L^{\prime}_{k_01} \cup L^{\prime}_{k_02}) \cup ( L^{\prime}_{i_1j_1}\cup L^{\prime}_{i_0j_0}).$$ 
        Finally, we have arranged $L^{\prime}_{k_01}$ and $L^{\prime}_{k_02}$ pair in a joint pair.

        \item [Case 2.] The joint pairs are given as follows: $(L^{\prime}_{i_0j_0} \cup L^{\prime}_{k_01}) \cup (L^{\prime}_{k_02}\cup L^{\prime}_{i_1j_1})$.
        We apply the type B move to get 
        $$(L^{\prime}_{i_0j_0} \cup L^{\prime}_{k_01}) \cup (L^{\prime}_{k_02}\cup L^{\prime}_{i_1j_1})\xrightarrow[]{B}( L^{\prime}_{k_02}\cup L^{\prime}_{k_01}) \cup (L^{\prime}_{i_0j_0}\cup L^{\prime}_{i_1j_1}).$$ 
        Then we apply the type A move on the first pair to get the $k_0$ component into the joint pair as our convention, i.e.
        $$( L^{\prime}_{k_02}\cup L^{\prime}_{k_01}) \cup (L^{\prime}_{i_0j_0}\cup L^{\prime}_{i_1j_1})\xrightarrow[]{A}(L^{\prime}_{k_01}\cup L^{\prime}_{k_02}) \cup (L^{\prime}_{i_0j_0}\cup L^{\prime}_{i_1j_1}).$$
        
        \item [Case 3.] The joint pairs are given as follows: $(L^{\prime}_{k_01} \cup L^{\prime}_{i_0j_0}) \cup ( L^{\prime}_{i_1j_1}\cup L^{\prime}_{k_02})$. 
        We apply the type B move to get 
        $$(L^{\prime}_{k_01} \cup L^{\prime}_{i_0j_0}) \cup ( L^{\prime}_{i_1j_1}\cup L^{\prime}_{k_02})\xrightarrow[]{B}(L^{\prime}_{k_01} \cup L^{\prime}_{k_02}) \cup ( L^{\prime}_{i_1j_1}\cup L^{\prime}_{i_0j_0}),$$ 
        $k_0$ component into the joint pair satisfying our convention.
        
        \item [Case 4.]  The joint pairs are given as follows: $(L^{\prime}_{i_0j_0} \cup L^{\prime}_{k_01}) \cup ( L^{\prime}_{i_1j_1}\cup L^{\prime}_{k_02})$.
        We apply the type A move on the first pair to get 
        $$(L^{\prime}_{i_0j_0} \cup L^{\prime}_{k_01}) \cup ( L^{\prime}_{i_1j_1}\cup L^{\prime}_{k_02})\xrightarrow[]{A}(L^{\prime}_{k_01} \cup L^{\prime}_{i_0j_0} ) \cup ( L^{\prime}_{i_1j_1}\cup L^{\prime}_{k_02}).$$ 
        Then we apply the type B move on the latter union to get   
        $$(L^{\prime}_{k_01} \cup L^{\prime}_{i_0j_0} ) \cup ( L^{\prime}_{i_1j_1}\cup L^{\prime}_{k_02})\xrightarrow[]{B}(L^{\prime}_{k_01} \cup L^{\prime}_{k_02}) \cup ( L^{\prime}_{i_1j_1}\cup L^{\prime}_{i_0j_0}).$$
    \end{enumerate}
    \item We may repeat the above step finitely many times until we exhaust all the pairs. 
    In the end, we achieve a round surgery diagram $L' = \bigcup_{i=1}^{n}(L'_{i1}\cup L'_{i2})$ such that $L'_{i1} \cup L'_{i2}$ is a joint pair as per our convention for each $1\leq i\leq n$. 
    However, the round surgery coefficients of $L'_{ij}$ are not same as $L_{ij}$ 
    \item For each $i^{th}$-joint pair in $L$ and $L^{\prime}$, there exists an integer $k_i$ such that the round surgery coefficients on both of the joint pair agree.
    If such an integer does not exist for any joint pair, then the corresponding Dehn surgery coefficient will differ, which is a contradiction to the hypothesis that both round surgery diagrams correspond to the same Dehn surgery diagram.
    This change in the round surgery coefficient is the application of the equivalence move 1 to each of the joint pairs in the surgery diagram. 
\end{enumerate}
\end{proof}

Now, we are in a position to prove the Equivalence Theorem (Theorem \ref{thm: EqThm}).

\begin{proof}[Proof of Theprem \ref{thm: EqThm}]
By Theorem \ref{BridgeTheorem}, we convert the round surgery diagrams $R_1$ and $R_2$ to the Dehn surgery diagrams $S_1$ and $S_2$, respectively. 
Since $R_1$ and $R_2$ correspond to the same 3-manifolds, $S_1$ and $S_2$ also correspond to the same 3-manifold. 
By Kirby's theorem, we can apply the Kirby moves finitely many times to obtain $S_2$ from $S_1$.
We denote $I_{i}$ to be an intermediate surgery diagram obtained after applying a certain type of Kirby move to $I_{i-1}$, for $1\leq i\leq k$ such that $I_0=S_1$ and $I_k=S_1$. 
By Theorem \ref{BridgeTheorem}, each $I_i$ corresponds to a round surgery diagram of joint pairs, denoted by $J_i$. %We denote each intermediate round surgery diagram by $J_i$ corresponding to each $I_i$. 

By definition, the equivalence moves 3 and 4 play the role of Kirby moves of types 1 and 2 on a round surgery diagram, respectively.
Thus, each Kirby move of type 1 can be replaced with the equivalence move 3 and a Kirby move of type 2 by the equivalence move 4 in the corresponding round surgery diagrams.
Each Kirby move of type 1 increases or decreases the number of components of $I_i$ by $1$. 
Since $S_2$ has an even number of components, we apply the Kirby move of type 1 an even number, say $2p$, times. 

Each time we remove $(\pm1)$-framed unknot from Dehn surgery diagram $I_a$ via Kirby move of type 1 to obtain $I_{a+1}$, we do not change the corresponding round surgery diagram $J_a$, i.e., $J_{a+1}= J_a$, for some $1\leq a \leq k$. 
Futher, each time we add a $(\pm1)$-framed unknot to Dehn surgery diagram $I_a$ via Kirby move of type 1, we add two unknots in joint pair $U_{a1} \cup U_{a2}$ with round 2-suregry coefficient $(\pm1)$ on $U_{a2}$ in the corresponding round surgery diagram $J_a$ via equivalence move 3, repectively. 
Thus, after applying all the Kirby moves on $S_1$ to get $S_2$, the corresponding round surgery diagram $J_k$ has $2p_0$-many unlinked unknots, where the integer $p_0 \leq p$. 
These unlinked unknots are either in a joint pair among themselves or in a joint pair with other components. 
We pair those unknots, which are in joint pairs with other components than among themselves.
And apply shuffle moves of type A and B as discussed in  
cases 1 to 4 in the proof of the Lemma \ref{lem: EM1n2onRSD}.
As a result, we get $p_0$-many joint pairs on unlinked unknots. 
We delete all these pairs by equivalence move 3. 
After deletion, we get a link $J'_k$ link isotopic to $S_2$. 

Since $S_2$ is obtained by applying Lemma \ref{lem: JPtoOSD}, the round surgery diagram $R_2$ also has the same link type as $S_2$.
In particular, $J'_k$ and $R_2$ are two round surgery diagrams corresponding to the same Dehn surgery diagram $S_2$.
We apply Lemma \ref{lem: EM1n2onRSD} to get $R_2$ from $J'_k$ by applying equivalence moves 1 and 2 finitely many times.
Hence, we have obtained the round surgery diagram $R_2$ from $R_1$ by a finite sequence of the equivalence moves 1--4. 
\end{proof}

\begin{rem}
    We can think of the equivalence theorem as an analogue of the Kirby theorem for the round surgery diagrams. 
\end{rem}

\section{Applications}\label{Sec: Application}
\subsection{Round 1-surgery diagram to a Kirby diagram.}

%A round surgery diagram $L= \bigcup_{i=1}^{n}(L_{i1} \cup L_{i2})$ in $\s^3$ of a 3-manifold is called round 1-surgery diagram if on each pair $L_{i1} \cup L_{i2}$ we perform only round 1-surgery with round 1-surgery coefficient $n_{ij}$ on $L_{ij}$ for $1 \leq i\leq n$ and $1\leq j \leq 2$. 
In this subsection, we would like to convert a round 1-surgery diagram into a Kirby diagram consisting of the attaching region of a 1-handle and the attaching sphere of a 2-handle.  
We use the idea of the proof of Lemma \ref{Lem: RHasOH} for the conversion.

\begin{lem}[Round 1-surgery diagram to a Kirby diagram]\label{Lem: R1SDtoOSD}
    Suppose a link $L=L_{11} \cup L_{12}$ is round 1-surgery diagram of 3-manifold $M$ with $n_{11}$ and $n_{12}$ are round 1-surgery coefficients on $L_{11}$ and $L_{12}$, respectively. Then, it can be converted into a Kirby diagram  of the attaching region of a 1-handle and the attaching sphere of a 2-handle, such that
        \begin{enumerate}
        \item one 3-ball of the attaching sphere of the 1-handle $\s^0\times\D^3$ lies on the one component, and the other 3-ball lies on the other component of the round 1-surgery link. It makes the two components link into a 2-surgery knot, which goes over the attached 1-handle. 
         \item The two 3-balls of the attaching sphere of the 1-handle are placed so that the orientation on the 2-surgery knot agrees with the orientation on $L_{11}$ and $L_{12}$.
        \item The 2-surgery curve gets a new framing $n_{11} + n_{12} + 2 lk(L_{11}, L_{12})$.
    \end{enumerate}
\end{lem}

\begin{proof}
By Lemma \ref{Lem: RHasOH}, a round $k$-handle can be realised as a pair of $k$ and $(k+1)$-handles. 
We want to track the attaching regions of the ordinary handles under this realization. 
Since we are interested in a round 1-surgery, we take $k=1$.
Schematically, we can think of the round 1-handle $\s^1 \times \D^{1} \times \D^{2}$ as a solid cuboid ABCDHGFE. 

In the cuboid, we identify the $y$-axis with $\s^1$, the $z$-axis with $\D^2$, and the $x$-axis with $\D^1$, as shown in the Figure \ref{fig: 1RHto0-1OH}.
\begin{figure}[ht]
    \centering
    \includegraphics[scale=0.4]{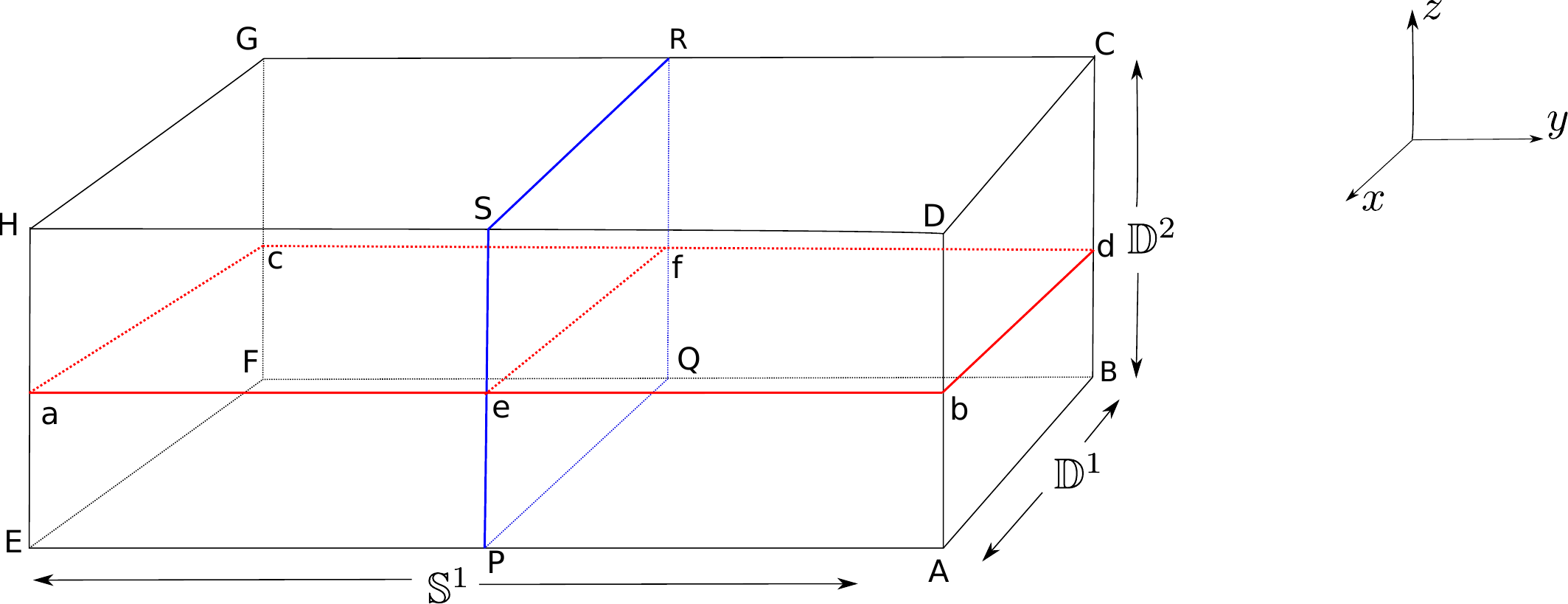}
    \caption{Round 1-handle attachment as a pair of ordinary handles of indices 1 and 2}
    \label{fig: 1RHto0-1OH}
\end{figure}
The front face $ADHE$ and back face $BCGF$ are the attaching regions of the round 1-handle. 
The circles $ab \subset ADHE$ and $cd \subset BCGF$ is the link $L_{11}\cup L_{12}$ in $\s^3$. 
In both of the circles $ab$ and $cd$, we can find small segments disjoint from the linking and knotting. 
Suppose the segments $eb \subset ab$ and $fd \subset cd$ denote those parts of the knots.
Now, we can think of $\s^1$ as a union of two intervals, i.e., $\D^1 \cup \D^1$. 
We break this cuboid into two solid cuboids, ABCDSRQP and EFGHSRQP. 
We treat the first half ABCDSRQP as a 1-handle and the second half EFGHSRQP as a 2-handle.
The squares APSD and BQRC are the attaching region of the 1-handle as it identifies with $\D^1\times \D^2 \times \s^0 = \D^3 \times \s^0$.  
Further, the attaching region of the 2-handle is identified with the cuboid EPQFGHSR, i.e. $\D^1 \times \D^1 \times \D^2= \D^2 \times \D^2$. 
The core curve of the 2-handle can be identified with $ae\cup ef \cup fc \cup ca$, where $ae$ and $cf$ are part of the different components of the round 1-surgery link.
The parts $ef$ and $ac$ go over the attached 1-handles. 
In particular, we get a Kirby diagram consisting of the attaching region of a 1-handle and the attaching sphere of a 2-handle.

In conversion, we add a 3-ball to each component such that the orientation on the new knot going through the 3-balls agrees with the orientation given on the components of the round 1-surgery link.
The sides $HS\cup SR \cup  RG \cup GH$ form the framing curve of the 2-surgery knot. 
We may consider the 3-balls such that they do not intersect with the individual framing of the components of $L$.
The framing curve $L'_{11}(n_{11})$ has linking number $n_{11}$ with $L_{11}$ and linking number $lk(L_{11}, L_{12})$ with $L_{12}$. 
Similarly, the framing curve $L'_{12}(n_{12})$ has linking number $n_{12}$ with $L_{12}$ and linking number $lk(L_{11}, L_{12})$ with $L_{12}$. 
Both of these framing curves join through the 2-handle to become a single framing curve of the new 2-surgery knot.
Thus, the new framing of the 2-surgery curve is the sum of the above individual linking numbers.
Hence, the 2-surgery knot has $n_{11} + n_{12} + 2 lk(L_{11}, L_{12})$ framing.  
\end{proof}

\begin{figure}
    \centering
    \includegraphics[scale=0.5]{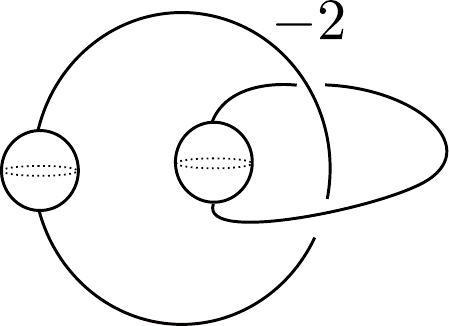}
    \caption{Kirby diagram of a 4-manifold with a circle bundle over a torus as a boundary.}
    \label{fig:KirbyDiagramFor3-torus}
\end{figure}

\begin{rem}\label{rem: R1SDtoKD}
    The above Lemma \ref{Lem: R1SDtoOSD} gives another way to construct new Kirby diagrams whose boundary is determined by the round 1-surgery diagrams.
    For example, the round 1-surgery diagram of the circle bundle over the torus given in Figure \ref{fig:3-torus} can be converted into a Kirby diagram $\mathcal{K}$ consisting of handles of indices 1 and 2 (see Figure \ref{fig:KirbyDiagramFor3-torus}).
    In particular, Kirby diagrams obtained from the above lemma can be used to distinguish 4-manifolds by showing that they have different boundary manifolds. 
\end{rem}

We have proved the following theorem in the Lemmas \ref{lem: JPtoOSD} and \ref{Lem: R1SDtoOSD}. 

\begin{thm}\label{Thm: RSDtoOSD}
Suppose $L=L_{1} \cup \ldots \cup L_{n}$ is a round surgery diagram of a closed connected oriented smooth $3$-manifold such that each $L_{i}$ is a two-component link.
\begin{enumerate}
\item[(i)] If $L_{i}$ is not in a joint pair of surgeries, then it corresponds to a Kirby diagram consisting of handles of indices 1 and 2.
\item[(ii)] And if $L_i$ is a joint pair link, then it corresponds to a Kirby diagram consisting of only 2-handles (in particular, an integral Dehn surgery diagram) as mentioned in Lemma \ref{lem: JPtoOSD}.
\end{enumerate} 
\end{thm} 

\subsection{Extension of the Taut foliation on the surgered 3-manifold.}
On a thickened torus $\mathbb{T}^2 \times [1,2]$, there are infinitely many taut foliations by annuli as a leaf.
For a given foliation $\mathcal{F}$, a leaf annulus intersects with the boundary torus $\mathbb{T}^2 \times \{i\}$ on a simple closed curve $d_i$ for $i=1,2$.
We may fix an identification of the thickened torus with $\mathbb{R}^2/\mathbb{Z}^2 \times [1,2]$. 
Using this identification, we may describe the slope $s(d_i)$ of $d_i$, which is called the {\it suture slope}.
We claim that $s(d_1)= s(d_2)$. 
Since $d_1 \cup d_2$ is the boundary of an annulus of $\mathcal{F}$, projection of this annulus on to $\mathbb{T}^2 \times \{2\}$ gives a homotopy between $d_1$ and $d_2$ in $\mathbb{T}^2\times \{ 2\}$.
Thus, simple closed curves are in the same homotopy class, which implies they have the same slope, i.e., $s(d_1)=s(d_2)$.

In particular, any taut foliation by an annulus of a thickened torus induced isotopic simple closed curves on the boundary tori. 
For more details about the taut foliation, the reader may refer to \cite{FolBook}.

We can extend a taut foliation to the resultant 3-manifold obtained after performing a round 1-surgery on a fibred link of two components if their suture slopes on the toroidal boundary are the same. 

\begin{lem}
\label{lem: Application1}
   Suppose $L$ is a fibred link of two components in $\mathbb{S}^3$. 
   Assume we perform a round 1-surgery on $L$ with the same round 1-surgery coefficient $n$ on both components. 
   Then the resultant 3-manifold $M_L$ has a taut foliation extended from the foliation of the fibred link. 
\end{lem}

\begin{proof}
We have given a two-component fibred link $L= L_{11} \cup L_{12}$. 
Suppose $N(L):= N(L_{11}) \cup N(L_{12})$, where $N(L_{1j})$ denotes the tubular neighbourhood of the  $L_{1j}$. 
Since $L$ is a fibred link, $\overline{\s^3 \setminus N(L)} $ has a codimension one foliation such that the leaf $S$ is a Seifert surface of the link.
Denote $d_j:= S\cap N(L_{1j})$ for $j=1, 2$. 
For each $j=1,2$, we may express $d_j$ with respect to the canonical coordinate curves described in the proof of Lemma \ref{Lem:R1Diagram}, 
\begin{equation}\label{Eq: SutureSlope}
    d_j= a_j \cdot m_j+ b_j \cdot l_{j}. 
\end{equation}
Since $d_j$ is a push-off of the $L_{1j}$ along the surface $S$, $b_j=1$. 
Suppose $S_j$ is a Seifert surface of $L_{1j}$. Then the knot $L_{12}$ intersects with the surface $S_1$, $lk(L_{11}, L_{12})$-many times which implies the push off $d_1$ goes $lk(L_{11}, L_{12})$-many times around the meridian $m_{1}$. A similar argument holds for $d_2$. Thus, $a_j= lk(L_{11}, L_{12})$ for $j=1,2$. 
From Equation \ref{Eq: SutureSlope}, we have 
\begin{equation} \label{Eq: SutureSlope2}
d_j = lk(L_{11}, L_{12}) \cdot m_{j} + l_{j}
\end{equation}

Suppose we perform a round 1-surgery on $L$ with round surgery coefficients $n_{1j}$ on $L_{1j}$ for $j=1,2$. 
Using the notation from the proof of the Lemma \ref{Lem:R1Diagram}, the glueing is determined by the images
\begin{eqnarray}\label{Eq: SutureSlope3}
    r_{j} &=& n \cdot m_{j} + l_{j}\\
\text{and } s_{j} &=& m_{j}, 
\end{eqnarray}
where $r_{j} = \phi(\s^1 \times \{j\} \times  \{p\})$ and $s_{j}= \phi(\{q\} \times \{j\} \times \s^1)$ for some $p, q\in \s^1$. 

%Denote the generating curves $\s^1 \times \{j\} \times \{p\}$ and $\{q\} \times \{j\}\times \s^1$ on $\mathbb{T}^2 \times \{j\}$ by $x^{j}$ and $y^j$ respectively. 

In particular, we get 
\begin{eqnarray}
\label{Eq: SutureSlope4}
    m_{j}&=& s_{j}\\
\label{Eq: SutureSlope4.1}    \text{and } l_{j}&=& r_{j}- n\cdot s_{j}.
\end{eqnarray}

By Equations \ref{Eq: SutureSlope2}, \ref{Eq: SutureSlope4} and \ref{Eq: SutureSlope4.1}, we get 
\begin{eqnarray*}
    d_j &=& lk(L_{11}, L_{12})\cdot s_{j} + r_{j} - n \cdot s_{j}\\
    &=& lk(L_{11}, L_{12}) - n)\cdot s_{j} + r_{j}. 
\end{eqnarray*}

From the observation above, we know that we can extend the foliation if the slope of $d_j$ is the same in $\mathbb{T}^2 \times [1,2]$. 
We take a taut foliation $\mathcal{F}$ on $\mathbb{T}^2 \times [1,2]$ that induce slope $lk(L_{11}, L_{12})-n $ on the boundary tori under the identification of $\mathbb{T}^2\times \{j\}$ with $\R^2/ \mathbb{Z}^2$ defined by $r_j \mapsto (1,0)$ and $s_j \mapsto (0,1)$.
Hence, the resultant 3-manifold has a taut foliation. 
\end{proof}

\begin{rem}\label{lem: Application2}
    Suppose $M$ is a disjoint union of two copies of $\s^3$, and each 3-sphere  $\s^3_j$ contains a fibred knot $K_j$ for $j=1, 2$. 
    Assume that we perform a round 1-surgery on $M$ along $L=K_1 \sqcup K_2$ with the round 1-surgery coefficients $n$ on both $K_1$ and $K_2$.
    Then the resultant 3-manifold admits a taut foliation extended from the foliations of the fibred knots.
\end{rem}

We can summarise the discussions above in Theorem~\ref{Thm: InftyTautFoln}. The theorem extends the taut foliations from a link complement to the manifold obtained by round 1-surgery along that link. %In particular, it can be thought of as a round surgery version of the Li--Roberts theorem (cf. \cite{Li2014}). 
Although in \cite{Miyoshi}, Miyoshi has proved the extensions of the foliations on the manifolds obtained by (foliated) round surgery. Our result provides explicit conditions on framings along the lines of Li--Roberts theorem (cf. \cite{Li2014}). 
%We state this observation as the following theorem. The following theorem may follow from the existing literature on taut foliations. Here, we present it from the round surgery point of view.    

\begin{thm}\label{Thm: InftyTautFoln}
    Let $M$ be a 3-manifold obtained by performing round 1-surgery on a fibred link $L_{11}\cup L_{12} \subset \s^3$ with the same round 1-surgery coefficient $n$ on $L_{11}$ and $L_{12}$. 
    Then, for each $n\in \mathbb{Z}$, $M$ admits a taut foliation $\mathcal{F}_n$.     
\end{thm}

Suppose we perform two round 1-surgeries on the given link $L_{11} \cup L_{12}$ with two round surgery coefficients $n_{11}=n_{12}=n$ and $n'_{11}=n'_{12}=m$ such that $n \neq m$.
By Theorem \ref{thm: SameClassOfR1SD}, we know that both of the round 1-surgeries correspond to the same 3-manifold $M$. And we get two taut foliations $\mathcal{F}_n$ for surgery coefficient $n$ and $\mathcal{F}_m$ for surgery coefficient $m$ on $M$ by Lemma \ref{lem: Application1}. 
We ask the following question. 
\begin{ques}\label{Qn: distinct_fol}
    For integers $n\neq m$, is there any relation between foliations $\mathcal{F}_n$ and $\mathcal{F}_m$ on $M$?
\end{ques}

The tangent bundle of the leaves of the taut foliation provides an integrable 2-plane field on $ M$. 
Due to Eliashberg--Thurston in \cite{ThursEliashConfol}, a small $C^0$-perturbation to this 2-plane field gives a tight contact structure $\xi$ on $M$ such that $(M,\xi)$ is a contact manifold.
For details about the contact structure, a reader can refer to \cite{Book-HG}. 

As a consequence of Theorem \ref{Thm: InftyTautFoln}, we get infinitely many tight contact structures on $M$, which is obtained by performing round 1-surgery on $L= L_{11} \cup L_{12}$ in $\s^3$. 
Formally, we may state it as follows. 
\begin{cor}
    \label{cor: ExsTightCS}
    Suppose $M$ is a closed connected oriented 3-manifold obtained by round 1-surgery on a fibred link $L_{11}\cup L_{12}$ in $\s^3$ with the same round 1-surgery coefficient $n$ on both components. 
    Then, for each $n\in \mathbb{Z}$, $M$ has a tight contact structure $\xi_n$ on it. 
\end{cor}

For each $n\in \mathbb{Z}$, this construction yields a tight contact structure on the fixed manifold $M$. We do not know that the resulting contact structures are distinct up to isotopy. Thus, one may ask a question similar to Question~\ref{Qn: distinct_fol}.

\bibliographystyle{plain}
\bibliography{Reference}

@article {Asimov,
    AUTHOR = {Asimov, Daniel},
     TITLE = {Round handles and non-singular {M}orse-{S}male flows},
   JOURNAL = {Ann. of Math. (2)},
  FJOURNAL = {Annals of Mathematics. Second Series},
    VOLUME = {102},
      YEAR = {1975},
    NUMBER = {1},
     PAGES = {41--54},
      ISSN = {0003-486X},
   MRCLASS = {58F10 (57D65)},
  MRNUMBER = {380883},
MRREVIEWER = {Robert\ Roussarie},
       DOI = {10.2307/1970972},
       URL = {https://doi.org/10.2307/1970972},
}

@article {Thurs,
    AUTHOR = {Thurston, W. P.},
     TITLE = {Existence of codimension-one foliations},
   JOURNAL = {Ann. of Math. (2)},
  FJOURNAL = {Annals of Mathematics. Second Series},
    VOLUME = {104},
      YEAR = {1976},
    NUMBER = {2},
     PAGES = {249--268},
      ISSN = {0003-486X},
   MRCLASS = {57D30},
  MRNUMBER = {425985},
MRREVIEWER = {D.\ B.\ Fuchs},
       DOI = {10.2307/1971047},
       URL = {https://doi.org/10.2307/1971047},
}

@article {Miyoshi,
    AUTHOR = {Miyoshi, Shigeaki},
     TITLE = {Foliated round surgery of codimension-one foliated manifolds},
   JOURNAL = {Topology},
  FJOURNAL = {Topology. An International Journal of Mathematics},
    VOLUME = {21},
      YEAR = {1982},
    NUMBER = {3},
     PAGES = {245--261},
      ISSN = {0040-9383},
   MRCLASS = {57R30 (57R90)},
  MRNUMBER = {649757},
MRREVIEWER = {Takashi\ Tsuboi},
       DOI = {10.1016/0040-9383(82)90008-8},
       URL = {https://doi.org/10.1016/0040-9383(82)90008-8},
}

@article{Li2014,
  author  = {Li, Tao and Roberts, Rachel},
  title   = {Taut foliations in knot complements},
  journal = {Pacific Journal of Mathematics},
  year    = {2014},
  volume  = {269},
  number  = {1},
  pages   = {149--168},
  doi     = {10.2140/pjm.2014.269.149}
}

@article {MFS,
    AUTHOR = {Matveev, S. V. and Fomenko, A. T. and Sharko, V. V.},
     TITLE = {Round {M}orse functions and isoenergetic surfaces of
              integrable {H}amiltonian systems},
   JOURNAL = {Mat. Sb. (N.S.)},
  FJOURNAL = {Matematicheski\u i\ Sbornik. Novaya Seriya},
    VOLUME = {135(177)},
      YEAR = {1988},
    NUMBER = {3},
     PAGES = {325--345, 414},
      ISSN = {0368-8666},
   MRCLASS = {58F05 (57M99)},
  MRNUMBER = {937644},
MRREVIEWER = {A.\ Morimoto},
       DOI = {10.1070/SM1989v063n02ABEH003276},
       URL = {https://doi.org/10.1070/SM1989v063n02ABEH003276},
}

@article {STT,
    AUTHOR = {Joubert, G. and Rosenberg, H.},
     TITLE = {Plongement du tore {$T\sb{2}$} dans la sph\`ere {$S\sb{3}$}},
   JOURNAL = {Cahiers Topologie G\'eom. Diff\'erentielle},
  FJOURNAL = {Cahiers de Topologie et G\'eom\'etrie Diff\'erentielle},
    VOLUME = {11},
      YEAR = {1969},
     PAGES = {323--328},
      ISSN = {0008-0004},
   MRCLASS = {57.20},
  MRNUMBER = {273638},
MRREVIEWER = {O.\ G.\ Harrold},
}

@article {Jiro,
    AUTHOR = {Adachi, Jiro},
     TITLE = {Contact round surgery and symplectic round handlebodies},
   JOURNAL = {Internat. J. Math.},
  FJOURNAL = {International Journal of Mathematics},
    VOLUME = {25},
      YEAR = {2014},
    NUMBER = {5},
     PAGES = {1450050, 25},
      ISSN = {0129-167X,1793-6519},
   MRCLASS = {57R17 (53D35 57R65)},
  MRNUMBER = {3215225},
MRREVIEWER = {Yildiray\ Ozan},
       DOI = {10.1142/S0129167X14500505},
       URL = {https://doi.org/10.1142/S0129167X14500505},
}

@article {Jiro19,
    AUTHOR = {Adachi, Jiro},
     TITLE = {Contact round surgery and {L}utz twists},
   JOURNAL = {Internat. J. Math.},
  FJOURNAL = {International Journal of Mathematics},
    VOLUME = {30},
      YEAR = {2019},
    NUMBER = {4},
     PAGES = {1950019, 31},
      ISSN = {0129-167X,1793-6519},
   MRCLASS = {57R17 (53D35 57R65)},
  MRNUMBER = {3950816},
MRREVIEWER = {Ferit\ \"Ozt\"urk},
       DOI = {10.1142/S0129167X19500198},
       URL = {https://doi.org/10.1142/S0129167X19500198},
}

@book {Book-HG,
    AUTHOR = {Geiges, Hansj\"{o}rg},
     TITLE = {An introduction to contact topology},
    SERIES = {Cambridge Studies in Advanced Mathematics},
    VOLUME = {109},
 PUBLISHER = {Cambridge University Press, Cambridge},
      YEAR = {2008},
     PAGES = {xvi+440},
      ISBN = {978-0-521-86585-2},
   MRCLASS = {57R17 (53D35)},
  MRNUMBER = {2397738},
MRREVIEWER = {John B. Etnyre},
       DOI = {10.1017/CBO9780511611438},
       URL = {https://doi.org/10.1017/CBO9780511611438},
}

@book {Book-Rolfsen,
    AUTHOR = {Rolfsen, Dale},
     TITLE = {Knots and links},
    SERIES = {Mathematics Lecture Series},
    VOLUME = {No. 7},
 PUBLISHER = {Publish or Perish, Inc., Berkeley, CA},
      YEAR = {1976},
     PAGES = {ix+439},
   MRCLASS = {55-01},
  MRNUMBER = {515288},
}

@book {Book-S,
    AUTHOR = {Saveliev, Nikolai},
     TITLE = {Lectures on the topology of 3-manifolds},
    SERIES = {De Gruyter Textbook},
   EDITION = {revised},
      NOTE = {An introduction to the Casson invariant},
 PUBLISHER = {Walter de Gruyter \& Co., Berlin},
      YEAR = {2012},
     PAGES = {xii+207},
      ISBN = {978-3-11-025035-0},
   MRCLASS = {57N10 (57M25 57M27)},
  MRNUMBER = {2893651},
}

@article {L,
    AUTHOR = {Lickorish, W. B. R.},
     TITLE = {A representation of orientable combinatorial {$3$}-manifolds},
   JOURNAL = {Ann. of Math. (2)},
  FJOURNAL = {Annals of Mathematics. Second Series},
    VOLUME = {76},
      YEAR = {1962},
     PAGES = {531--540},
      ISSN = {0003-486X},
   MRCLASS = {54.78},
  MRNUMBER = {151948},
MRREVIEWER = {A.\ D.\ Wallace},
       DOI = {10.2307/1970373},
       URL = {https://doi.org/10.2307/1970373},
}

@article {W,
    AUTHOR = {Wallace, Andrew H.},
     TITLE = {Modifications and cobounding manifolds},
   JOURNAL = {Canadian J. Math.},
  FJOURNAL = {Canadian Journal of Mathematics. Journal Canadien de
              Math\'{e}matiques},
    VOLUME = {12},
      YEAR = {1960},
     PAGES = {503--528},
      ISSN = {0008-414X,1496-4279},
   MRCLASS = {57.10},
  MRNUMBER = {125588},
MRREVIEWER = {Morris\ W.\ Hirsch},
       DOI = {10.4153/CJM-1960-045-7},
       URL = {https://doi.org/10.4153/CJM-1960-045-7},
}

@article {RK,
    AUTHOR = {Kirby, Robion},
     TITLE = {A calculus for framed links in {$S^{3}$}},
   JOURNAL = {Invent. Math.},
  FJOURNAL = {Inventiones Mathematicae},
    VOLUME = {45},
      YEAR = {1978},
    NUMBER = {1},
     PAGES = {35--56},
      ISSN = {0020-9910},
   MRCLASS = {57A10 (55A25)},
  MRNUMBER = {467753},
MRREVIEWER = {Roger Fenn},
       DOI = {10.1007/BF01406222},
       URL = {https://doi.org/10.1007/BF01406222},
}

@book {Book-GompfStip,
    AUTHOR = {Gompf, Robert E. and Stipsicz, Andr\'{a}s I.},
     TITLE = {{$4$}-manifolds and {K}irby calculus},
    SERIES = {Graduate Studies in Mathematics},
    VOLUME = {20},
 PUBLISHER = {American Mathematical Society, Providence, RI},
      YEAR = {1999},
     PAGES = {xvi+558},
      ISBN = {0-8218-0994-6},
   MRCLASS = {57N13 (14J80 32Q55 57-02 57R17 57R57 57R65)},
  MRNUMBER = {1707327},
MRREVIEWER = {Nikolai N. Saveliev},
       DOI = {10.1090/gsm/020},
       URL = {https://doi.org/10.1090/gsm/020},
}

@incollection {ThursEliashConfol,
    AUTHOR = {Eliashberg, Yakov M. and Thurston, William P.},
     TITLE = {Confoliations},
 BOOKTITLE = {Collected works of {W}illiam {P}. {T}hurston with commentary.
              {V}ol. {I}. {F}oliations, surfaces and differential geometry},
     PAGES = {281--351},
      NOTE = {Reprint of [1483314]},
 PUBLISHER = {Amer. Math. Soc., Providence, RI},
      YEAR = {[2022] \copyright 2022},
      ISBN = {978-1-4704-6388-5; [9781470468330]; [9781470451646]},
   MRCLASS = {53C15 (57K30 57K33 57R30)},
  MRNUMBER = {4554446},
}

@book {FolBook,
    AUTHOR = {Calegari, Danny},
     TITLE = {Foliations and the geometry of 3-manifolds},
    SERIES = {Oxford Mathematical Monographs},
 PUBLISHER = {Oxford University Press, Oxford},
      YEAR = {2007},
     PAGES = {xiv+363},
      ISBN = {978-0-19-857008-0},
   MRCLASS = {57R30 (37C85 57M50 57N10)},
  MRNUMBER = {2327361},
MRREVIEWER = {Thilo\ Kuessner},
}

\end{document}